\documentclass[11 pt]{article}
\usepackage{epsfig}
\usepackage{color}
\usepackage{subfigure}

\usepackage{bm}% the package can get heiti for all.

\usepackage{indentfirst}
\usepackage{amsmath}
\usepackage{amssymb}
\usepackage{amsthm}

\numberwithin{figure}{section}
\numberwithin{equation}{section}

 \allowdisplaybreaks

\date{}

\begin{document}
\baselineskip 15pt

\title{\bf Asymptotic solvers for ordinary differential equations with multiple frequencies}
\author{\begin{tabular}[t]{c@{\extracolsep{0em}}c}
{\large Marissa CONDON} \\
{School of Electronic Engineering, Dublin City University} \\
{\normalsize{\it E-mail: marissa.condon@dcu.ie}}\\
{\large Alfredo DEA\~{N}O} \\
{Depto.de Matem\'{a}ticas, Universidad Carlos III de Madrid} \\
{\normalsize{\it E-mail: alfredo.deanho@uc3m.es}}\\
{\large Jing GAO \thanks{Communicating author}} \\
{School of Mathematics and Statistics, Xi'an Jiaotong University} \\
{\normalsize{\it E-mail: jgao@mail.xjtu.edu.cn}}\\
{\large Arieh ISERLES} \\
 {DAMTP, Centre for Mathematical Sciences, University of Cambridge} \\
 { \normalsize {\it E-mail: A.Iserles@damtp.cam.ac.uk}}\\
\end{tabular}}

\maketitle

\begin{abstract}
We construct asymptotic expansions for ordinary differential equations
with highly oscillatory forcing terms, focussing on the case of multiple,
non-commensurate frequencies. We derive an asymptotic expansion
in inverse powers of the oscillatory parameter and use its truncation as
an exceedingly effective means to discretize the differential equation in
question. Numerical examples illustrate the effectiveness of the method.
\end{abstract}
\bigskip

{\bf{Mathematics Subject Classification:}} 65L05,65D30,42B20,42A10

\section{Introduction}
Our concern in this paper is with the numerical solution of highly oscillatory
ordinary differential equations of the form
\begin{eqnarray}
\boldsymbol{y}'(t) =
\boldsymbol{f}(\boldsymbol{y}(t)) + \sum_{m = 1}^{M} \boldsymbol{a}_m(t) e^{i \omega_m t} ,\quad t \geq 0, \quad \boldsymbol{y}(0) = \boldsymbol{y}_0 \in \mathbb{C}^d,
\end{eqnarray}
where $\boldsymbol{f}: \mathbb{C}^d \rightarrow \mathbb{C}^d$
and $\boldsymbol{a}_1, \cdots, \boldsymbol{a}_M: \mathbb{R}_+ \rightarrow \mathbb{C}^d$ are analytic and
$\omega_1, \omega_2, \cdots, \omega_M \in \mathbb{R} \setminus \{0\}$ are given frequencies. We assume that at least
some of these frequencies are large, thereby causing the solution to oscillate and rendering numerical discretization of
(1.1) by classical methods expensive and inefficient. Many phenomena in engineer and physics are described by the oscillatoryly
differential equations(Chedjou2001, Fodjouong2007, Slight2008 and so on).

A special case of (1.1) with $\omega_{2m - 1} = m \omega$, $\omega_{2m} = - m\omega$, $m = 0, 1, \cdots, \lfloor M/2 \rfloor$,
where $\omega \gg 1$, is a special case of
\begin{eqnarray}
\boldsymbol{y}'(t) = \boldsymbol{f}(\boldsymbol{y}(t)) + G(\boldsymbol{y})
\sum_{k = -\infty}^{\infty} \boldsymbol{b}_k(t) e^{i k \omega t} ,\quad t \geq 0, \quad \boldsymbol{y}(0) = \boldsymbol{y}_0 \in \mathbb{C}^d,
\end{eqnarray}
where $G : \mathbb{C}^d \times \mathbb{C}^d \rightarrow
\mathbb{C}^d$ is smooth, which has been already analysed at some
length in (Condon, Dea\~{n}o and Iserles 2010). It has been proved
that the solution of (1.2) can be expanded asymptotically in
$\omega^{-1}$,
\begin{eqnarray}
\boldsymbol{y}(t) \sim \boldsymbol{p}_{0, 0}(t) + \sum_{r = 1}^{\infty} \frac{1}{\omega^r}
\sum_{m = -\infty}^{\infty} \boldsymbol{p}_{r, m}(t) e^{i m \omega t}, \quad t \geq 0,
\end{eqnarray}
where the functions $\boldsymbol{p}_{r, m}$, which are independent of $\omega$, can be derived recursively:
$\boldsymbol{p}_{r, 0}$ by solving a non-oscillatory ODE and $\boldsymbol{p}_{r, m}$, $m \neq 0$, by recursion.

An alternative approach, based upon the Heterogeneous Multiscale Method
(E and Engquist 2003), is due to Sanz-Serna (2009). Although the theory in
(Sanz-Serna 2009) is presented for a specific equation, it can be extended in
a fairly transparent manner to (1.2) and, indeed, to (1.1). It produces the
solution in the form
\begin{eqnarray}
\boldsymbol{y}(t) \sim \sum_{m = -\infty}^{\infty} \boldsymbol{\kappa}_m(t) e^{i m \omega t},
\end{eqnarray}
where $\kappa_m(t) = O(\omega^{-1})$, $m \in \mathbb{Z}$.
\footnote[1]{In the special case considered in (Sanz-Serna 2009) it is true that $\boldsymbol{\kappa}_m(t) = O(\omega^{-|m|})$,
$m \in \mathbb{Z}$, but this does not generalise to (1.2).}
Formally, (1.3) and (1.4) are linked by
\begin{eqnarray*}
&& \boldsymbol{\kappa}_m(t) = \left\{
\begin{array}{cccccc}
\sum\limits_{r = 0}^{\infty} \frac{1}{\omega^r} \boldsymbol{p}_{r, 0}(t), \quad m = 0,\\
\sum\limits_{r = 1}^{\infty} \frac{1}{\omega^r} \boldsymbol{p}_{r, m}(t), \quad m \neq 0.
\end{array}
\right.
\end{eqnarray*}
We adopt here the approach of (Condon et al. 2010), because it allows us to
derive the expansion in a more explicit form.

The highly oscillatory term in (1.2) is periodic in $t\omega$: the main difference
with our model (1.1) is that we allow the more general setting of almost periodic
terms (Besicovitch 1932). It is justified by important applications, not least in
the modelling of nonlinear circuits (Giannini and Leuzzi 2004, Ram\'{i}rez, Su\'{a}rez,
Lizarraga and Collantes 2010).

Another difference is that we allow in (1.1) only a finite number of distinct
frequencies in the forcing term. This is intended to prevent the occurrence of
small denominators, familiar from asymptotic theory (Verhulst 1990). Note
that (Chartier, Murua and Sanz-Serna 2012) (cf. also (Chartier, Murua and Sanz-Serna 2010))
employs similar formalism - a finite number of multiple, noncommensurate
frequencies - except that it does so within the `body' of the
differential operator, rather than in the forcing term.

We commence our analysis by letting $\mathcal{U}_0 = \{1, 2, \cdots, M\}$ and $\omega_j = \kappa_j \omega$,
$j = 1, \cdots, M$, where $\omega$ is a large number which will serve as our asymptotic
parameter. Consequently, we can rewrite (1.1) in a form that emphasises the
similarities and identifies the differences with (1.2),
\begin{eqnarray}
\boldsymbol{y}'(t) = \boldsymbol{f}(\boldsymbol{y}(t)) + \sum_{m \in
\mathcal{U}_0} \boldsymbol{a}_m(t) e^{i \kappa_m \omega t} , \quad t
\geq 0, \quad \boldsymbol{y}(0) = \boldsymbol{y}_0 \in \mathbb{C}^d,
\end{eqnarray}
Section 2 is devoted to a `warm up exercise', an asymptotic expansion of a
linear version of (1.5), namely
\begin{eqnarray}
\boldsymbol{y}'(t) = \boldsymbol{A} \boldsymbol{y}(t) + \sum_{m \in
\mathcal{U}_0} \boldsymbol{a}_m(t) e^{i \kappa_m \omega t} ,\quad t
\geq 0, \quad \boldsymbol{y}(0) = \boldsymbol{y}_0 \in \mathbb{C}^d,
\end{eqnarray}
where $A$ is a $d \times d$ matrix. Of course, the solution of (1.5) can be written
explicitly, but this provides little insight into the real size of different components.
The asymptotic expansion is considerably more illuminating, as well as
hinting at the general pattern which we might expect once we turn our gaze
to the nonlinear equation (1.5).

For the ODE with the highly oscillatory forcing terms with multiple
frequencies, the asymptotic method is superior to the standard
numerical methods. With less computational expense, this asymptotic
method can obtain the higher accuracy. Especially, the asymptotic
expansion with a fixed number of terms becomes more accurate when
increasing the oscillator parameter $\omega$.

In Section 3 we will demonstrate the existence of sets
\begin{eqnarray*}
\mathcal{U}_0 \subseteq \mathcal{U}_1 \subseteq \mathcal{U}_2 \subseteq \cdots
\end{eqnarray*}
and of a mapping
$$
\sigma : \bigcup\limits_{r = 0}^{\infty} \mathcal{U}_r \rightarrow \mathbb{R}
$$
such that the solution of (1.5) can be written in the form
\begin{eqnarray}
\boldsymbol{y}(t) \sim \boldsymbol{p}_{0, 0}(t) + \sum_{r = 1}^{\infty} \frac{1}{\omega^r} \sum_{m \in \mathcal{U}_r}
\boldsymbol{p}_{r, m}(t) e^{i \sigma_m \omega t}, \quad t \geq 0.
\end{eqnarray}
As can be expected, the original parameters $\{\kappa_1, \kappa_2, \cdots, \kappa_M\}$ form a subset of
$\mathcal{U}_r$. However, we will see in the sequel that the set $\sigma(\mathcal{U}_r)$ is substantially larger
for $r \geq 3$. In the sequel we refer to elements of $\sigma(\mathcal{U}_r)$ as $\{\sigma_m : m \in \mathcal{U}_r\}$.

The functions $\boldsymbol{p}_{r, m}$, which are all independent of $\omega$, are constructed explicitly
in a recursive manner. We will demonstrate that the sets $\mathcal{U}_r$ are composed
of $n$-tuples of nonnegative integers.

In Section 4 we accompany our narrative by a number of computational
results. Setting the error functions are represented as
\begin{eqnarray*}
\boldsymbol{\epsilon}_s(t, \omega) = \boldsymbol{y}(t)
-\boldsymbol{p}_{0, 0}(t) - \sum_{r = 1}^{s} \frac{1}{\omega^r}
\sum_{m \in \mathcal{U}_r} \boldsymbol{p}_{r, m}(t) e^{i \sigma_m
\omega t},
\end{eqnarray*}
we plot the error functions in the figures to illustrate the
theoretical analysis. The expansion solvers are convergent
asymptotically. That is, for every $\varepsilon > 0$,  fixed $s$,
the bounded interval for $t$, there exists $\omega_0 > 0$ such that
for $\omega > \omega_0$, the error function
$\left|\boldsymbol{\epsilon}_s(t, \omega)\right| < \varepsilon$.
However, for increasing $s$, fixed $t$ and $\omega$, we will come to
the convergence of the expansion in future paper.

\section{The linear case}

Our concern in this section is with the linear highly oscillatory ODE (1.6),
which we recall for convenience,
\begin{eqnarray}
\boldsymbol{y}' =
\boldsymbol{A} \boldsymbol{y} + \sum_{m \in \mathcal{U}_0} \boldsymbol{a}_m(t) e^{i \kappa_m \omega t} ,\quad t \geq 0, \quad \boldsymbol{y}(0) = \boldsymbol{y}_0 \in \mathbb{C}^d,
\end{eqnarray}
Its closed-form solution can be derived at once from standard variation of
constants,
\begin{eqnarray}
\boldsymbol{y}(t) = e^{t A} \boldsymbol{y}_0 + e^{t A} \sum\limits_{m \in \mathcal{U}_0}
\int_0^t e^{- x A} \boldsymbol{a}_m(x) e^{i \kappa_m \omega x}\hbox{d} x.
\end{eqnarray}
However, the finer structure of the solution is not apparent from (2.2) without
some extra work. Each of the integrals hides an entire hierarchy of scales, and
this becomes apparent once we expand them asymptotically.

The asymptotic expansion of integrals with simple exponential operators
is well known: given $g \in C^{\infty} [0, t)$ and $|\eta| \gg 1$,
\begin{eqnarray*}
\int_0^t g(x) e^{i \eta x}\hbox{d} x \sim - \sum_{r = 1}^{\infty} \frac{1}{(-i \eta)^r}
\left[g^{(r - 1)}(t) e^{i \eta t} - g^{(r - 1)}(0)\right]
\end{eqnarray*}
(Iserles, N{\o}rsett and Olver 2006). For any $m \in \mathcal{U}_0$ we thus take
$\boldsymbol{g}(x) = e^{- x A} \boldsymbol{a}_m(x)$ and $\eta = \kappa_m \omega$ (recall
that $\kappa_m \neq 0$), therefore
\begin{eqnarray*}
\boldsymbol{y}(t) &\sim& e^{t A} \boldsymbol{y}_0\\
&& - \sum_{m \in \mathcal{U}_0} \sum_{r = 1}^{\infty} \frac{1}{(- i \kappa_m \omega)^r}
\left[e^{i \kappa_m \omega t} \sum_{\ell = 0}^{r - 1} (-1)^{r - 1 - \ell} \binom{r - 1}{\ell}
A^{r - 1 - \ell} \boldsymbol{a}_m^{(\ell)}(t)\right.\\
&& \left. -e^{t A} \sum_{\ell = 0}^{r - 1} (-1)^{r - 1 - \ell} \binom{r - 1}{\ell}
A^{r - 1 - \ell} \boldsymbol{a}_m^{(\ell)}(0)\right]\\
&& =e^{t A} \boldsymbol{y}_0 + \sum_{r = 1}^{\infty} \frac{1}{\omega^r} \sum_{m \in \mathcal{U}_0}
\left[ \frac{e^{i \kappa_m \omega t}}{(i \kappa_m)^r} \sum_{\ell = 0}^{r - 1} (-1)^{\ell} \binom{r - 1}{\ell} A^{r - 1 - \ell}\boldsymbol{a}_m^{(\ell)}(t)\right.\\
&& \left. - \frac{1}{(i \kappa_m)^r} e^{t A} \sum_{\ell = 0}^{r - 1} (-1)^{\ell} \binom{r - 1}{\ell} A^{r - 1 - \ell} \boldsymbol{a}_m^{(\ell)}(0)\right].
\end{eqnarray*}
We deduce the expansion (1.7), with $\mathcal{U}_m = \{0\} \cup \mathcal{U}_0 = \{0, 1, \cdots, M\}$,
$m \in \mathbb{N}$, and the coefficients
\begin{eqnarray*}
&& \boldsymbol{p}_{0, 0}(t) = e^{t A} \boldsymbol{y}_0,\\
&& \boldsymbol{p}_{r, 0}(t) = - e^{t A} \sum_{\ell = 0}^{r - 1} (-1)^{\ell} \binom{r - 1}{\ell} A^{r - 1 - \ell}
\sum_{m \in \mathcal{U}_0} \frac{\boldsymbol{a}_m^{(\ell)}(0)}{(i \kappa_m)^r}, \quad r \in \mathbb{N},\\
&& \boldsymbol{p}_{r, m}(t) = \frac{1}{(i \kappa_m)^r} \sum_{\ell = 0}^{r - 1} (-1)^{\ell} \binom{r - 1}{\ell} A^{r - 1 - \ell}
\boldsymbol{a}_m^{(\ell)}(t), \quad r \in \mathbb{N}, \quad m \in \mathcal{U}_0.
\end{eqnarray*}

We thus recover an expansion of the form (1.7), where $\sigma_0 = 0$ and $\sigma_m = \kappa_m$,
$m \in \mathcal{U}_0$. Note that linearity `locks' frequencies: each $\boldsymbol{p}_{r, m}$ depends just on $\boldsymbol{a}_m$,
$m \in \mathcal{U}_0$. This is no longer true in a nonlinear setting.

\section{The asymptotic expansion}

\subsection{The recurrence relations}

We are concerned with expanding asymptotically the solution of
\begin{eqnarray}
\boldsymbol{y}' = \boldsymbol{f}(\boldsymbol{y}) + \sum_{m \in
\mathcal{U}_0} \boldsymbol{a}_m(t) e^{i \kappa_m \omega t} ,\quad t
\geq 0, \quad \boldsymbol{y}(0) = \boldsymbol{y}_0 \in \mathbb{C}^d,
\end{eqnarray}
The function $f$ is analytic, thus for every $n \in \mathbb{N}$
there exists the $n$th differential of $\boldsymbol{f}$, a function
$\boldsymbol{f}_n : \overbrace{\mathbb{C}^d \times \mathbb{C}^d
\times \mathbb{C}^d}^{n \quad times} \rightarrow \mathbb{C}^d$ such
that for every sufficiently small $|t| > 0$
\begin{eqnarray*}
\boldsymbol{f}(\boldsymbol{y}_0 + t \varepsilon) = \boldsymbol{f}(\boldsymbol{y}_0) + \sum_{n = 1}^{\infty} \frac{t^n}{n !}  \boldsymbol{f}_n(\boldsymbol{y}_0) [\varepsilon, \cdots, \varepsilon].
\end{eqnarray*}
Note that $\boldsymbol{f}_n$ is linear in all its arguments in the square brackets.

We substitute (1.7) in both sides of (3.1) and expand about $\boldsymbol{p}_{0, 0}(t)$,
\begin{eqnarray*}
\boldsymbol{y}' &=& \boldsymbol{p}'_{0, 0} + \sum_{m \in \mathcal{U}_1} i \sigma_m  \boldsymbol{p}_{1, m} e^{i \sigma_m \omega t}\\
&& \quad + \sum_{r = 1}^{\infty} \frac{1}{\omega^r} \left[\sum_{m \in \mathcal{U}_r} \boldsymbol{p}'_{r, m}e^{i \sigma_m \omega t}
+ i \sigma_m \sum_{m \in \mathcal{U}_{r + 1}} \boldsymbol{p}_{r + 1, m} e^{i \sigma_m \omega t}\right]\\
&=& \boldsymbol{f}\left(\boldsymbol{p}_{0, 0} + \sum_{r = 1}^{\infty} \frac{1}{\omega^r}\sum_{m \in \mathcal{U}_r}\boldsymbol{p}_{r, m} e^{i \sigma_m \omega t}\right) + \sum_{m \in \mathcal{U}_0} \boldsymbol{a}_m e^{i \kappa_m \omega t}\\
&=& \boldsymbol{f}(\boldsymbol{p}_{0, 0}) + \sum_{n = 1}^{\infty} \frac{1}{n!} \boldsymbol{f}_n(\boldsymbol{p}_{0, 0})
\left[\sum_{\ell_1 = 1}^{\infty} \frac{1}{\omega^{\ell_1}} \sum_{k_1 \in \mathcal{U}_{\ell_1}} \boldsymbol{p}_{\ell_1, k_1} e^{i \sigma_{k_1} \omega t}, \cdots,\right.\\
&& \quad \left. \sum_{\ell_n = 1}^{\infty} \frac{1}{\omega^{\ell_n}}\sum_{k_n \in \mathcal{U}_{\ell_n}} \boldsymbol{p}_{\ell_n, k_n} e^{i \sigma_{k_n} \omega t}\right] + \sum_{m \in \mathcal{U}_0} \boldsymbol{a}_m e^{i \kappa_m \omega t}\\
&=& \boldsymbol{f}(\boldsymbol{p}_{0, 0}) + \sum_{n = 1}^{\infty} \frac{1}{n!} \sum_{\ell_1 = 1}^{\infty} \cdots \sum_{\ell_n = 1}^{\infty}
\frac{1}{\omega^{\ell_1 + \ell_2 + \cdots + \ell_n}} \sum_{k_1 \in \mathcal{U}_{\ell_1}}\cdots \sum_{k_n \in \mathcal{U}_{\ell_n}} \boldsymbol{f}_n(\boldsymbol{p}_{0, 0}) \left[\boldsymbol{p}_{\ell_1, k_1},\right.\\
&& \quad \left. \cdots,\boldsymbol{p}_{\ell_n, k_n}\right] e^{i (\sigma_{k_1} + \cdots + \sigma_{k_n}) \omega t}  + \sum_{m \in \mathcal{U}_0} \boldsymbol{a}_m e^{i \kappa_m \omega t}\\
&=& \boldsymbol{f}(\boldsymbol{p}_{0, 0}) + \sum_{n = 1}^{\infty} \frac{1}{n!} \sum_{r = n}^{\infty}
\frac{1}{\omega^r} \sum_{\boldsymbol{\ell} \in \mathbb{I}^o_{n, r}} \sum_{k_1 \in \mathcal{U}_{\ell_1}}\cdots \sum_{k_n \in \mathcal{U}_{\ell_n}} \boldsymbol{f}_n(\boldsymbol{p}_{0, 0}) \left[\boldsymbol{p}_{\ell_1, k_1},\cdots,\boldsymbol{p}_{\ell_n, k_n}\right]\\
&&  \times e^{i (\sigma_{k_1} + \cdots + \sigma_{k_n}) \omega t}  + \sum_{m \in \mathcal{U}_0} \boldsymbol{a}_m e^{i \kappa_m \omega t}\\
&=& \boldsymbol{f}(\boldsymbol{p}_{0, 0}) + \sum_{r = 1}^{\infty} \frac{1}{\omega^r} \sum_{n = 1}^{r} \frac{1}{n!}
 \sum_{\boldsymbol{\ell} \in \mathbb{I}^o_{n, r}} \sum_{k_1 \in \mathcal{U}_{\ell_1}}\cdots \sum_{k_n \in \mathcal{U}_{\ell_n}} \boldsymbol{f}_n(\boldsymbol{p}_{0, 0}) \left[\boldsymbol{p}_{\ell_1, k_1},\cdots,\boldsymbol{p}_{\ell_n, k_n}\right]\\
&&  \times e^{i (\sigma_{k_1} + \cdots + \sigma_{k_n}) \omega t}  + \sum_{m \in \mathcal{U}_0} \boldsymbol{a}_m e^{i \kappa_m \omega t}
\end{eqnarray*}
where
$$
\mathbb{I}^o_{n, r} = \{\boldsymbol{\ell} \in \mathbb{N}^n : \boldsymbol{\ell}^T \boldsymbol{1} = r\}, \quad 1 \leq n \leq r.
$$

There is a measure of redundancy in the last expression: for example there
are two terms in $\mathbb{I}^o_{2, 3}$, namely (1,2) and (2,1), but they produce identical expressions.
Consequently, we may lump them together, paying careful attention
to their multiplicity. More formally, we let
$$
\mathbb{I}_{n, r} = \{\boldsymbol{\ell} \in \mathbb{N}^n : \boldsymbol{\ell}^T \boldsymbol{1} = r, \ell_1 \leq \ell_2 \leq \cdots \leq \ell_n\}, \quad 1 \leq n \leq r,
$$
the set of ordered partitions of $r$ into $n$ natural numbers and allow $\theta_{\boldsymbol{\ell}}$ stand for
the \textit{multiplicity} of $\boldsymbol{\ell}$, i.e. the number of terms in $\mathbb{I}^o_{n, r}$ that can be brought to it
by permutations. For example, $\theta_{1, 2} = 2$, while for r = 4 there are five terms,
$$
\theta_4 = 1, \quad \theta_{1, 3} = 2, \quad \theta_{2, 2} = 1, \quad \theta_{1, 1, 2} = 3, \quad \theta_{1, 1, 1, 1} = 1.
$$
We introduce the multiplicity and obtain a more compact form for the
equation
\begin{eqnarray}
&& \boldsymbol{p}'_{0, 0} + \sum_{m \in \mathcal{U}_1} i \sigma_m  \boldsymbol{p}_{1, m} e^{i \sigma_m \omega t}\\
&& \quad + \sum_{r = 1}^{\infty} \frac{1}{\omega^r} \left[\sum_{m \in \mathcal{U}_r} \boldsymbol{p}'_{r, m}e^{i \sigma_m \omega t}
+ \sum_{m \in \mathcal{U}_{r + 1}} i \sigma_m  \boldsymbol{p}_{r + 1, m} e^{i \sigma_m \omega t}\right]\nonumber\\
&& = \boldsymbol{f}(\boldsymbol{p}_{0, 0}) + \sum_{r = 1}^{\infty} \frac{1}{\omega^r} \sum_{n = 1}^{r} \frac{1}{n!}
 \sum_{\boldsymbol{\ell} \in \mathbb{I}_{n, r}} \theta_{\boldsymbol{\ell}} \sum_{k_1 \in \mathcal{U}_{\ell_1}}\cdots \sum_{k_n \in \mathcal{U}_{\ell_n}} \boldsymbol{f}_n(\boldsymbol{p}_{0, 0}) \left[\boldsymbol{p}_{\ell_1, k_1},\cdots,\boldsymbol{p}_{\ell_n, k_n}\right]\nonumber\\
&&  \qquad \qquad \times e^{i (\sigma_{k_1} + \cdots + \sigma_{k_n}) \omega t}  + \sum_{m \in \mathcal{U}_0} \boldsymbol{a}_m e^{i \kappa_m \omega t}.\nonumber
\end{eqnarray}

We separate different powers of $\omega$ in (3.2). The outcome is
\begin{eqnarray}
\boldsymbol{p}'_{0, 0} + \sum_{m \in \mathcal{U}_1} i \sigma_m  \boldsymbol{p}_{1, m} e^{i \sigma_m \omega t}
= \boldsymbol{f}(\boldsymbol{p}_{0, 0}) + \sum_{m \in \mathcal{U}_0} \boldsymbol{a}_m e^{i \kappa_m \omega t}.
\end{eqnarray}
for $r = 0$ and
\begin{eqnarray}
&& \sum_{m \in \mathcal{U}_r} \boldsymbol{p}'_{r, m} e^{i \sigma_m \omega t}
+ \sum_{m \in \mathcal{U}_{r + 1}} i \sigma_m  \boldsymbol{p}_{r + 1, m} e^{i \sigma_m \omega t}\\
&& = \sum_{n = 1}^{r} \frac{1}{n!}\sum_{\boldsymbol{\ell} \in \mathbb{I}_{n, r}} \theta_{\boldsymbol{\ell}} \sum_{k_1 \in \mathcal{U}_{\ell_1}}\cdots \sum_{k_n \in \mathcal{U}_{\ell_n}} \boldsymbol{f}_n(\boldsymbol{p}_{0, 0}) \left[\boldsymbol{p}_{\ell_1, k_1},\cdots,\boldsymbol{p}_{\ell_n, k_n}\right] e^{i (\sigma_{k_1} + \cdots + \sigma_{k_n}) \omega t} \nonumber
\end{eqnarray}
for $r \in \mathbb{N}$.

Before we embark on detailed examination of the cases $r = 0, 1, 2$, followed
by the general case, we must impose an additional set of conditions on the
coefficients $\boldsymbol{p}_{r, m}$. Similarly to the expansion of (1.2) in (Condon et al. 2010),
we obtain non-oscillatory differential equations for the coefficients $\boldsymbol{p}_{r, 0}$, $r \in \mathbb{Z}_{+}$,
which require initial conditions. We do so by imposing the original initial
condition from (3.1) on $\boldsymbol{p}_{0, 0}$ and requiring that the terms at origin sum to zero
at every $\omega$ scale. In other words,
\begin{eqnarray}
\boldsymbol{p}_{0, 0}(0) = \boldsymbol{y}_0, \quad
\boldsymbol{p}_{r, 0}(0) = - \sum_{m \in \mathcal{U}_r \setminus \{0\}} \boldsymbol{p}_{r, m}(0), \quad r \in \mathbb{N}.
\end{eqnarray}

\subsection{The first few values of $r$}

The expansion (1.7) exhibits two distinct hierarchies of scales: both \textit{amplitudes}
$\omega^{-r}$ for $r \in \mathbb{Z}_+$ and, for each $r \in \mathbb{N}$, \textit{frequencies}
$e^{i \sigma_m \omega t}$. In (3.3) and (3.4) we
have already separated amplitudes. Next we separate frequencies.

For r = 0 (3.3) and (3.5) yield the original ODE (3.1) without a forcing
term,
$$
\boldsymbol{p}'_{0, 0} = \boldsymbol{f}(\boldsymbol{p}_{0, 0}), \quad t \geq 0, \quad \boldsymbol{p}_{0, 0}(0) = \boldsymbol{y}_0,
$$
as well as the recursions
$$
\boldsymbol{p}_{1, m} = \frac{\boldsymbol{a}_m}{i \kappa_m}, \quad m = 1, \cdots, M
$$
(recall that $\kappa_m \neq 0$). Therefore $\sigma_m = \kappa_m$, $m = 1, \cdots, M$. We set, for reasons
that will become apparent in the sequel,
$$
\mathcal{U}_1 = \{0\} \cup \mathcal{U}_0 = \{0, 1, \cdots, M\},
$$
with $\sigma_0 = 0$.

For $r = 1$ we have $\mathbb{I}_{1, 1} = \{1\}$, $\theta_1 = 1$, and (3.4) yields
\begin{eqnarray*}
\sum_{m \in \mathcal{U}_1} \boldsymbol{p}'_{1, m} e^{i \sigma_m \omega t}
+ \sum_{m \in \mathcal{U}_2} i \sigma_m  \boldsymbol{p}_{2, m} e^{i \sigma_m \omega t}
= \sum_{m \in \mathcal{U}_1} \boldsymbol{f}_1(\boldsymbol{p}_{0, 0}) \left[\boldsymbol{p}_{1, m}\right] e^{i \sigma_m \omega t}.
\end{eqnarray*}

We set
$$
\mathcal{U}_2 = \mathcal{U}_1 = \{0, 1, \cdots, M\}
$$
and (recalling that $\boldsymbol{p}_{1, m}$ are already known for $m \neq 0$)
\begin{eqnarray*}
&& \boldsymbol{p}'_{1, 0} = \boldsymbol{f}_1(\boldsymbol{p}_{0, 0})[\boldsymbol{p}_{1, 0}], \quad t \geq 0, \quad \boldsymbol{p}_{1, 0}(0) = - \sum_{m \in \mathcal{U}_1 \setminus \{0\}} \boldsymbol{p}_{1, m}(0), \\
&& \boldsymbol{p}_{2, m} = \frac{1}{i \kappa_m} \{\boldsymbol{f}_1(\boldsymbol{p}_{0, 0})[\boldsymbol{p}_{1, m}] - \boldsymbol{p}'_{1, m}\}, \quad m \in \mathcal{U}_2 \setminus \{0\}.
\end{eqnarray*}

An explanation is in order with regard to our imposition of $0 \in \mathcal{U}_1$. We
could have accounted for all the $r = 1$ terms in (3.4) without any need of the
$\boldsymbol{p}_{1, 0}$ term. However, in that case the outcome would not have been consistent
with the initial condition (3.5) and this is the rationale for the addition of this
term.

Our next case is $r = 2$. The case is not so straightforward to
deduce. Since $\mathbb{I}_{1, 2} = \{2\}$ and $\mathbb{I}_{2, 2} =
\{(1, 1)\}$, we have from (3.4)
\begin{eqnarray}
&& \sum_{m \in \mathcal{U}_2} \boldsymbol{p}'_{2, m} e^{i \sigma_m \omega t}
+ \sum_{m \in \mathcal{U}_3} i \sigma_m  \boldsymbol{p}_{3, m} e^{i \sigma_m \omega t}
= \sum_{m \in \mathcal{U}_2} \boldsymbol{f}_1(\boldsymbol{p}_{0, 0}) \left[\boldsymbol{p}_{2, m}\right] e^{i \sigma_m \omega t}\nonumber\\
&& + \frac{1}{2} \sum_{m_1 \in \mathcal{U}_1} \sum_{m_2 \in \mathcal{U}_1} \boldsymbol{f}_2(\boldsymbol{p}_{0, 0}) [\boldsymbol{p}_{1, m_1}, \boldsymbol{p}_{1, m_2}] e^{i (\sigma_{m_1} + \sigma_{m_2}) \omega t}.
\end{eqnarray}
We need to choose $\mathcal{U}_3$ to match frequencies in the above formula. The set
$\mathcal{U}_2$ accounts for the frequencies $0, \kappa_1, \kappa_2, \cdots, \kappa_M$ but we must also account for
$\kappa_i + \kappa_j$ for $i, j = 1, 2, \cdots, M$. Therefore we let
$$
\mathcal{U}_3 = \mathcal{U}_2 \cup \{(m_1, m_2) : 0 \leq m_1 \leq
m_2 \leq M\}.
$$
We note two important points. Firstly, for $i \neq j$, $\kappa_i + \kappa_j$ can be obtained
for $(j, i)$, as well as for $(i, j)$. Secondly, it might well happen that there exist
$i, j, k \in \{1, \cdots, M\}$, $i \leq j$, such that $\kappa_i + \kappa_j = \kappa_k$ {-in that case we do
not include $(i, j)$ in $\mathcal{U}_3$. This motivates the definition of the \textit{multiplicity} of
$m \in \mathcal{U}_3 \setminus \mathcal{U}_2$ (which we will generalise in the sequel to all sets $\mathcal{U}_r$). Thus,
for every $0 \leq \ell_1 \leq \ell_2 \leq M$ we let $\rho^m_{\ell_1, \ell_2}$ equal the number of cases when
$\kappa_{\pi(\ell_1)} + \kappa_{\pi(\ell_2)} = \kappa_m$, where $\pi (\boldsymbol{\ell})$ is a permutation of $\boldsymbol{\ell}$.
Likewise, we let $\rho^{m_1, m_2}_{\ell_1, \ell_2}$, where $0 \leq \ell_1 \leq \ell_2 \leq M$ and $1 \leq m_1 \leq m_2 \leq M$,
be the number of permutations such that $\kappa_{\pi(\ell_1)} + \kappa_{\pi(\ell_2)} = \kappa_{m_1} + \kappa_{m_2}$.

We can now separate frequencies. Firstly, the non-oscillatory term, corresponding
to $\sigma_0 = 0$. It yields the non-oscillatory ODE
$$
\boldsymbol{p}'_{2, 0} = \boldsymbol{f}_1(\boldsymbol{p}_{0, 0})[\boldsymbol{p}_{2, 0}] + \frac{1}{2} \mathop{\sum_{\kappa_{\ell_1} + \kappa_{\ell_2} = 0}}_{\ell_1 \leq \ell_2} \rho^0_{\ell_1, \ell_2} \boldsymbol{f}_2(\boldsymbol{p}_{0, 0}) [\boldsymbol{p}_{1, \ell_1}, \boldsymbol{p}_{1, \ell_2}],
$$
whose initial condition, according to (3.5), is
\begin{eqnarray*}
\boldsymbol{p}_{2, 0}(0) = - \sum_{m \in \mathcal{U}_2 \setminus \{0\}} \boldsymbol{p}_{2, m}(0).
\end{eqnarray*}
Secondly, we match all the terms in $\mathcal{U}_2 \setminus \{0\}$, and this results in the recurrence
\begin{eqnarray*}
i \kappa_m  \boldsymbol{p}_{3, m}
= \boldsymbol{f}_1(\boldsymbol{p}_{0, 0}) \left[\boldsymbol{p}_{2, m}\right] - \boldsymbol{p}'_{2, m}
+ \frac{1}{2} \mathop{\sum_{\kappa_{\ell_1} + \kappa_{\ell_2} = \kappa_m}}_{\ell_1 \leq \ell_2} \rho^m_{\ell_1, \ell_2} \boldsymbol{f}_2(\boldsymbol{p}_{0, 0}) [\boldsymbol{p}_{1, \ell_1}, \boldsymbol{p}_{1, \ell_2}].
\end{eqnarray*}
Finally, we match the terms in $\mathcal{U}_3 \setminus \mathcal{U}_2$. Recall that these are pairs $(m_1, m_2)$
such that $m_1 \leq m_2$ and $\kappa_{m_1} + \kappa_{m_2} \neq \sigma_j$ for $j = 0, 1, \cdots, M$. We obtain the
recurrence
\begin{eqnarray*}
i (\kappa_{m_1} + \kappa_{m_2}) \boldsymbol{p}_{3, (m_1, m_2)}
= \frac{1}{2} \mathop{\sum_{\kappa_{\ell_1} + \kappa_{\ell_2} = \kappa_{m_1} + \kappa_{m_2}}}_{\ell_1 \leq \ell_2} \rho^{m_1, m_2}_{\ell_1, \ell_2} \boldsymbol{f}_2(\boldsymbol{p}_{0, 0}) [\boldsymbol{p}_{1, \ell_1}, \boldsymbol{p}_{1, \ell_2}].
\end{eqnarray*}
(There is no danger of dividing by zero since we have ensured that $\kappa_{m_1} + \kappa_{m_2} \neq \sigma_0 = 0$.) This accounts for all the terms in (3.6).

\subsection{The general case $r \geq 1$}

We consider the `level r' equations (3.4) noting that, by induction, the sets
$\mathcal{U}_{\ell}$ are known for $\ell = 1, 2, \cdots, r$ and $0 \in \mathcal{U}_r$. Moreover, we have already
constructed all the functions $\boldsymbol{p}_{\ell, m}$s for $m \in \mathcal{U}_{\ell} \setminus \{0\}$, $\ell = 1, 2, \cdots, r$, and all the
$\boldsymbol{p}_{\ell, 0}$s for $\ell = 0, 1, \cdots, r - 1$. The current task is to construct the set $\mathcal{U}_{r + 1}$, the
functions $\boldsymbol{p}_{r + 1, m}$ for $m \in \mathcal{U}_{r + 1} \setminus \{0\}$ and the function $\boldsymbol{p}_{r, 0}$.

It will follow soon that all the terms in $\mathcal{U}_{r + 1}$ are of the form $\kappa_{j_1} + \kappa_{j_2} + \cdots + \kappa_{j_q}$, where $q \leq r$ and $j_1 \leq j_2 \leq \cdots \leq j_q$, We commence by setting $\rho^{m_1, m_2, \cdots, m_q}_{\ell_1, \ell_2, \cdots, \ell_p}$
as the number of distinct p-tuples $(\ell_1, \ell_2, \cdots, \ell_p)$, where
$$
\ell_1, \ell_2, \cdots, \ell_p, m_1, m_2, \cdots, m_q \in \{0, 1, \cdots, M\}, \quad m_1 \leq m_2 \leq \cdots \leq m_q,
$$
such that
$$
\sum_{i = 1}^{p} \kappa_{\ell_i} = \sum_{i = 1}^{q} \kappa_{m_i}.
$$

Examining the formula (3.4) we observe that the terms on the right
hand side have exponents of the form $e^{i \eta \omega t}$, where
$$
\eta = \sigma_{k_1} + \sigma_{k_2} + \cdots + \sigma_{k_n}, \quad \sigma_{k_i} \in \mathcal{U}_{\ell_i}, \quad, i = 1, \cdots, n, \quad \boldsymbol{\ell} \in \mathbb{I}_{r, n}
$$
for some $n \in \{1, 2, \cdots, r\}$. It follows at once by induction on $r$ that there exist
$q \in \{1, 2, \cdots, r\}$ and $0 \leq m_1 \leq m_2 \leq \cdots \leq m_q \leq M$ such that
$$
\eta = \sum_{i = 1}^{q} \kappa_{m_i}.
$$
It might well be that such $\eta$ can be already accounted by
$\mathcal{U}_r$, in other words that there exists $m \in
\mathcal{U}_r$ such that $\eta = \sigma_m$. Otherwise we add to
$\mathcal{U}_r$ the ordered $q$-tuple $(m_1, m_2, \cdots, m_q)$.
This process, applied to all the terms on the right of (3.4),
produces the index set $\mathcal{U}_{r + 1}$
$$
\mathcal{U}_{r + 1} = \mathcal{U}_{r} \bigcup \left\{(m_1, \cdots,
m_q): 0 \leq m_1 \leq m_2 \leq \cdots \leq m_q \leq M,  q \in {1, 2,
\cdots, r}\right\}.
$$

We impose natural partial ordering on $\mathcal{U}_{r}$: first the singletons in lexicographic
ordering, then the pairs in lexicographic ordering, then the triplets
etc. This defines a relation $m_1 \preceq m_2$ for all $m_1, m_2 \in \mathcal{U}_r$. We let
$$
\mathcal{W}^n_{r, m} = \left\{(\boldsymbol{\ell}, \boldsymbol{k}):
k_i \in \mathcal{U}_{\ell_i}, \boldsymbol{\ell} \in \mathbb{I}_{n,
r}, \sum_{i = 1}^{n} \sigma_{k_i} = \sum_{i = 1}^{q} \sigma_{m_i},
k_1 \preceq \cdots \preceq k_n, m_1 \leq \cdots \leq m_q\right\},
$$
where $m \in \mathcal{U}_r$ and $n \in \{1, 2, \cdots, r\}$.

Let us commence our construction of recurrence relations by considering
$m \in \mathcal{U}_r \setminus \{0\}$. In that case we have
\begin{eqnarray}
&& i \sigma_m  \boldsymbol{p}_{r + 1, m} = - \boldsymbol{p}'_{r, m} \nonumber\\
&&+ \sum_{n = 1}^{r} \frac{1}{n!} \sum_{\boldsymbol{\ell} \in \mathbb{I}_{n, r}} \theta_{\boldsymbol{\ell}} \sum_{(\boldsymbol{\ell}, \boldsymbol{k}) \in \mathcal{W}^n_{r, m}} \rho^m_{\boldsymbol{k}} \boldsymbol{f}_n(\boldsymbol{p}_{0, 0}) \left[\boldsymbol{p}_{\ell_1, k_1}, \cdots, \boldsymbol{p}_{\ell_n, k_n}\right].
\end{eqnarray}
Next we consider $m \in \mathcal{U}_{r + 1} \setminus \mathcal{U}_{r}$. Now the first sum on the left of (3.4)
disappears and the outcome is
\begin{eqnarray}
i \sigma_m  \boldsymbol{p}_{r + 1, m}
= \sum_{n = 1}^{r} \frac{1}{n!} \sum_{\boldsymbol{\ell} \in \mathbb{I}_{n, r}} \theta_{\boldsymbol{\ell}}
\sum_{(\boldsymbol{\ell}, \boldsymbol{k}) \in
\mathcal{W}^n_{r, m}} \rho^m_{\boldsymbol{k}} \boldsymbol{f}_n(\boldsymbol{p}_{0, 0}) \left[\boldsymbol{p}_{\ell_1, k_1}, \cdots, \boldsymbol{p}_{\ell_n, k_n}\right].\nonumber\\
\end{eqnarray}
Finally we cater for the case $m = 0$: now the recurrence is a non-oscillatory
ODE,
\begin{eqnarray}
&& \boldsymbol{p}'_{r, 0} = \sum_{n = 1}^{r} \frac{1}{n!} \sum_{\boldsymbol{\ell} \in \mathbb{I}_{n, r}} \theta_{\boldsymbol{\ell}} \sum_{(\boldsymbol{\ell}, \boldsymbol{k}) \in \mathcal{W}^n_{r, 0}} \rho^0_{\boldsymbol{k}} \boldsymbol{f}_n(\boldsymbol{p}_{0, 0}) \left[\boldsymbol{p}_{\ell_1, k_1}, \cdots, \boldsymbol{p}_{\ell_n, k_n}\right], \quad t \geq 0,\nonumber\\
&& \\
&& \boldsymbol{p}_{r, 0}(0) = - \sum_{m \in \mathcal{U}_r \setminus \{0\}} \boldsymbol{p}_{r, m}(0).\nonumber
\end{eqnarray}

For example, in the case $r = 3$ we have
$$
\mathbb{I}_{1, 3} = \{3\}, \quad \mathbb{I}_{2, 3} = \{(1, 2)\}, \quad \mathbb{I}_{3, 3} = \{(1, 1, 1)\}, \quad \theta_3 = \theta_{1, 1, 1}= 1, \quad \theta_{1, 2} = 2,
$$
while
$$
\mathcal{U}_3 = \{0, 1, \cdots, M\} \cup \{(m_1, m_2) : m_1 \leq m_2, \kappa_{m_1} + \kappa_{m_2} \neq \kappa_m, \quad \forall m = 0, \cdots, M\}.
$$
We thus deduce from (3.7) that
\begin{eqnarray*}
&& i \kappa_m  \boldsymbol{p}_{4, m} = - \boldsymbol{p}'_{3, m} + \boldsymbol{f}_1(\boldsymbol{p}_{0, 0})[\boldsymbol{p}_{3, m}]
+ \mathop{\sum_{\kappa_{j_1} + \kappa_{j_2} = \kappa_{m}}}_{j_1 \leq j_2} \rho^{m}_{j_1, j_2} \boldsymbol{f}_2(\boldsymbol{p}_{0, 0}) [\boldsymbol{p}_{1, j_1}, \boldsymbol{p}_{2, j_2}]\\
&& \qquad \qquad + \frac{1}{6} \mathop{\sum_{\kappa_{j_1} + \kappa_{j_2} + \kappa_{j_3} = \kappa_{m}}}_{j_1 \leq j_2 \leq j_3} \rho^{m}_{j_1, j_2, j_3} \boldsymbol{f}_3(\boldsymbol{p}_{0, 0}) [\boldsymbol{p}_{1, j_1}, \boldsymbol{p}_{1, j_2}, \boldsymbol{p}_{1, j_3}]
\end{eqnarray*}
for all $m \in \mathcal{U}_2 \setminus \{0\}$ and
\begin{eqnarray*}
i (\kappa_{m_1} + \kappa_{m_2}) \boldsymbol{p}_{4, (m_1, m_2)} &=& - \boldsymbol{p}'_{3, (m_1, m_2)} + \boldsymbol{f}_1(\boldsymbol{p}_{0, 0})[\boldsymbol{p}_{3, (m_1, m_2)}]\\
&& + \mathop{\sum_{\kappa_{j_1} + \kappa_{j_2} = \kappa_{m_1} + \kappa_{m_2}}}_{j_1 \leq j_2} \rho^{m_1, m_2}_{j_1, j_2} \boldsymbol{f}_2(\boldsymbol{p}_{0, 0}) [\boldsymbol{p}_{1, j_1}, \boldsymbol{p}_{2, j_2}]\\
&& + \frac{1}{6} \mathop{\sum_{\kappa_{j_1} + \kappa_{j_2} + \kappa_{j_3} = \kappa_{m_1} + \kappa_{m_2}}}_{j_1 \leq j_2 \leq j_3} \rho^{m_1, m_2}_{j_1, j_2, j_3} \boldsymbol{f}_3(\boldsymbol{p}_{0, 0}) [\boldsymbol{p}_{1, j_1}, \boldsymbol{p}_{1, j_2}, \boldsymbol{p}_{1, j_3}]
\end{eqnarray*}
for $(m_1, m_2) \in \mathcal{U}_3 \setminus \mathcal{U}_2$. Next we use (3.8):
\begin{eqnarray*}
&& i (\kappa_{m_1} + \kappa_{m_2} + \kappa_{m_3}) \boldsymbol{p}_{4, (m_1, m_2, m_3)} \\
&=& \frac{1}{6} \mathop{\sum_{\kappa_{j_1} + \kappa_{j_2} + \kappa_{j_3} = \kappa_{m_1} + \kappa_{m_2} + \kappa_{m_3}}}_{j_1 \leq j_2 \leq j_3} \rho^{m_1, m_2, m_3}_{j_1, j_2, j_3} \boldsymbol{f}_3(\boldsymbol{p}_{0, 0}) [\boldsymbol{p}_{1, j_1}, \boldsymbol{p}_{1, j_2}, \boldsymbol{p}_{1, j_3}]
\end{eqnarray*}
for all $1 \leq m_1 \leq m_2 \leq m_3 \leq M$ such that $\kappa_{m_1} + \kappa_{m_2} + \kappa_{m_3} \neq \sigma_m$ for $m \in \mathcal{U}_3$.

Finally, we invoke (3.9) to derive a non-oscillatory ODE for $\boldsymbol{p}_{3, 0}$, namely
\begin{eqnarray*}
\boldsymbol{p}'_{3, 0} &=&  \boldsymbol{f}_1(\boldsymbol{p}_{0, 0})[\boldsymbol{p}_{3, 0}]
+ \sum_{\kappa_{j_1} + \kappa_{j_2} = 0} \rho^{0}_{j_1, j_2} \boldsymbol{f}_2(\boldsymbol{p}_{0, 0}) [\boldsymbol{p}_{1, j_1}, \boldsymbol{p}_{1, j_2}]\\
&& \qquad \qquad + \frac{1}{6} \mathop{\sum_{\kappa_{j_1} + \kappa_{j_2} + \kappa_{j_3} = 0}}_{j_1 \leq j_2 \leq j_3} \rho^{0}_{j_1, j_2, j_3} \boldsymbol{f}_3(\boldsymbol{p}_{0, 0}) [\boldsymbol{p}_{1, j_1}, \boldsymbol{p}_{1, j_2}, \boldsymbol{p}_{1, j_3}],\\
&& \boldsymbol{p}_{3, 0}(0) = - \sum_{m \in \mathcal{U}_3 \setminus \{0\}} \boldsymbol{p}_{3, m}(0).
\end{eqnarray*}

\subsection{A worked-out example}

Let $M = 3$ and
$$
\kappa_1 = 1, \quad \kappa_2 = \sqrt{2}, \quad \kappa_3 = - 1 - \sqrt{2}.
$$
Therefore $\mathcal{U}_1 = \mathcal{U}_2 = \{0, 1, 2, 3\}$,
$$
\sigma_0 = 0, \quad \sigma_1 = 1, \quad \sigma_2 = \sqrt{2}, \quad \sigma_3 = - 1 - \sqrt{2}
$$
and
$$
\rho^m_k = \delta_{k, m}, \quad k, m = 0, 1, 2, 3.
$$

Consequently, $\boldsymbol{p}'_{0, 0} = \boldsymbol{f}(\boldsymbol{p}_{0, 0})$, $\boldsymbol{p}_{0, 0}(0) = \boldsymbol{y}(0)$ and
$$
\boldsymbol{p}_{1, 1} = \frac{\boldsymbol{a}_1}{i}, \quad \boldsymbol{p}_{1, 2} = \frac{\boldsymbol{a}_2}{\sqrt{2} i}, \quad \boldsymbol{p}_{1, 3} = - \frac{\boldsymbol{a}_3}{(1 + \sqrt{2}) i}.
$$

We commence with $r = 1$. The ODE is now
$$
\boldsymbol{p}'_{1, 0} = \boldsymbol{f}_1(\boldsymbol{p}_{0, 0})[\boldsymbol{p}_{1, 0}], \quad t \geq 0, \quad
\boldsymbol{p}_{1, 0}(0) = - \boldsymbol{p}_{1, 1}(0) - \boldsymbol{p}_{1, 2}(0)- \boldsymbol{p}_{1, 3}(0),
$$
while the recurrences are
\begin{eqnarray*}
&& \boldsymbol{p}_{2, 1} = \frac{1}{i} \{\boldsymbol{f}_1(\boldsymbol{p}_{0, 0})[\boldsymbol{p}_{1, 1}] - \boldsymbol{p}'_{1, 1}\}, \qquad
\boldsymbol{p}_{2, 2} = \frac{1}{\sqrt{2} i} \{\boldsymbol{f}_1(\boldsymbol{p}_{0, 0})[\boldsymbol{p}_{1, 2}] - \boldsymbol{p}'_{1, 2}\}\\
&& \boldsymbol{p}_{2, 3} = - \frac{1}{(1 + \sqrt{2})i} \{\boldsymbol{f}_1(\boldsymbol{p}_{0, 0})[\boldsymbol{p}_{1, 3}] - \boldsymbol{p}'_{1, 3}\}.
\end{eqnarray*}

Note that the first two "levels" of the expansion are
$$
\boldsymbol{y}(t) \approx \boldsymbol{p}_{0, 0} + \frac{1}{\omega} \left[\boldsymbol{p}_{1, 0} + \boldsymbol{p}_{1, 1}e^{i \omega t}
+ \boldsymbol{p}_{1, 2}e^{i \sqrt{2} \omega t} + \boldsymbol{p}_{1, 3}e^{-i (1 + \sqrt{2}) \omega t}\right].
$$

Next to $\mathcal{U}_3 = \{0, 1, 2, 3, (1, 1), (1, 2), (1, 3), (2, 2), (2, 3), (3, 3)\}$, with $\sigma_{1, 1} = 2$,
$\sigma_{1, 2} = 1 + \sqrt{2}$, $\sigma_{1, 3} = - \sqrt{2}$, $\sigma_{2, 2} = 2 \sqrt{2}$, $\sigma_{2, 3} = -1$ and $\sigma_{3, 3} = - 2 - 2\sqrt{2}$.

The only way to obtain $\sigma_0 = 0$ using two terms from $\mathcal{U}_1$ is $0 + 0$, therefore
$\rho^0_{0, 0} = 1$, otherwise $\rho^0_{\ell_1, \ell_2} = 0$. However, to obtain $\sigma_m$ for $m = 1, 2, 3$ we have
two options: $0 + m$ and $m + 0$. Therefore $\rho^m_{0, m} = 2$, otherwise $\rho^m_{\ell_1, \ell_2} = 0$. For
$\rho^{m_1, m_2}_{\ell_1, \ell_2}$
we note that $\rho^{1, 1}_{1, 1} = \rho^{2, 2}_{2, 2} = \rho^{3, 3}_{3, 3} = 1$, $\rho^{1, 2}_{1, 2} = \rho^{1, 3}_{1, 3} = \rho^{2, 3}_{2, 3} = 2$, otherwise the coefficient is zero. Therefore
\begin{eqnarray*}
&& \boldsymbol{p}'_{2, 0} = \boldsymbol{f}_1(\boldsymbol{p}_{0, 0})[\boldsymbol{p}_{2, 0}] + \frac{1}{2}\boldsymbol{f}_2(\boldsymbol{p}_{0, 0})[\boldsymbol{p}_{1, 0}, \boldsymbol{p}_{1, 0}], \\
&& \boldsymbol{p}_{2, 0}(0) = - \boldsymbol{p}_{2, 1}(0) - \boldsymbol{p}_{2, 2}(0)- \boldsymbol{p}_{2, 3}(0),
\end{eqnarray*}
and the $O(\omega^{-2})$ terms are
$$
\frac{1}{\omega^2} [\boldsymbol{p}_{2, 0} + \boldsymbol{p}_{2, 1}e^{i \omega t} + \boldsymbol{p}_{2, 2}e^{i \sqrt{2} \omega t} + \boldsymbol{p}_{2, 3}e^{-i (1 + \sqrt{2}) \omega t}].
$$
Moreover,
\begin{eqnarray*}
&&\boldsymbol{p}_{3, 1} = \frac{1}{i} \{\boldsymbol{f}_1(\boldsymbol{p}_{0, 0})[\boldsymbol{p}_{2, 1}] - \boldsymbol{p}'_{2, 1} + \boldsymbol{f}_2(\boldsymbol{p}_{0, 0})[\boldsymbol{p}_{1, 0}, \boldsymbol{p}_{1, 1}]\},\\
&&\boldsymbol{p}_{3, 2} = \frac{1}{\sqrt{2}i} \{\boldsymbol{f}_1(\boldsymbol{p}_{0, 0})[\boldsymbol{p}_{2, 2}] - \boldsymbol{p}'_{2, 2} + \boldsymbol{f}_2(\boldsymbol{p}_{0, 0})[\boldsymbol{p}_{1, 0}, \boldsymbol{p}_{1, 2}]\},\\
&&\boldsymbol{p}_{3, 3} = - \frac{1}{(1 + \sqrt{2})i} \{\boldsymbol{f}_1(\boldsymbol{p}_{0, 0})[\boldsymbol{p}_{2, 3}] - \boldsymbol{p}'_{2, 3} + \boldsymbol{f}_2(\boldsymbol{p}_{0, 0})[\boldsymbol{p}_{1, 0}, \boldsymbol{p}_{1, 3}]\},\\
&&\boldsymbol{p}_{3, (1, 1)} = \frac{1}{4 i}\boldsymbol{f}_2(\boldsymbol{p}_{0, 0})[\boldsymbol{p}_{1, 1}, \boldsymbol{p}_{1, 1}],\\
&&\boldsymbol{p}_{3, (1, 2)} = \frac{1}{(1 + \sqrt{2}) i}\boldsymbol{f}_2(\boldsymbol{p}_{0, 0})[\boldsymbol{p}_{1, 1}, \boldsymbol{p}_{1, 2}],\\
&&\boldsymbol{p}_{3, (1, 3)} = - \frac{1}{\sqrt{2} i}\boldsymbol{f}_2(\boldsymbol{p}_{0, 0})[\boldsymbol{p}_{1, 1}, \boldsymbol{p}_{1, 3}],\\
&&\boldsymbol{p}_{3, (2, 2)} = \frac{1}{4 \sqrt{2} i}\boldsymbol{f}_2(\boldsymbol{p}_{0, 0})[\boldsymbol{p}_{1, 2}, \boldsymbol{p}_{1, 2}],\\
&&\boldsymbol{p}_{3, (2, 3)} = - \frac{1}{i}\boldsymbol{f}_2(\boldsymbol{p}_{0, 0})[\boldsymbol{p}_{1, 2}, \boldsymbol{p}_{1, 3}],\\
&&\boldsymbol{p}_{3, (3, 3)} = - \frac{1}{4 (1 + \sqrt{2}) i}\boldsymbol{f}_2(\boldsymbol{p}_{0, 0})[\boldsymbol{p}_{1, 3}, \boldsymbol{p}_{1, 3}].
\end{eqnarray*}

The number of terms in $\mathcal{U}_4$ is significantly larger: we need to add to $\mathcal{U}_3$
a further nine terms - altogether we have 19 terms, which are displayed in
Table 1. All the nonzero $\rho^m_{\ell_1, \ell_2, \ell_3}$
are displayed there as well. Note that there
is no $(1, 2, 3)$ term, because $\kappa_1 + \kappa_2 + \kappa_3 = 0$, hence it is counted together with
zero.

In particular,
\begin{eqnarray*}
&& \boldsymbol{p}'_{3, 0} = \boldsymbol{f}(\boldsymbol{p}_{0, 0})[\boldsymbol{p}_{3, 0}] + \boldsymbol{f}_2(\boldsymbol{p}_{0, 0})[\boldsymbol{p}_{1, 0}, \boldsymbol{p}_{2, 0}] + \frac{1}{6}\boldsymbol{f}_3(\boldsymbol{p}_{0, 0})[\boldsymbol{p}_{1, 0}, \boldsymbol{p}_{1, 0}, \boldsymbol{p}_{1, 0}], \\
&& \boldsymbol{p}_{3, 0}(0) = - \sum_{j = 1}^{3}\boldsymbol{p}_{3, j}(0) - \sum_{j = 1}^{3} \sum_{\ell=j}^{3}\boldsymbol{p}_{3, (j, \ell)}(0).
\end{eqnarray*}
We conclude that the $O(\omega^{-3})$ terms in the asymptotic expansion of $\boldsymbol{y}$ are
\begin{eqnarray*}
&& \frac{1}{\omega^3} \left[\boldsymbol{p}_{3, 0} + \boldsymbol{p}_{3, 1}e^{i \omega t} + \boldsymbol{p}_{3, 2}e^{i \sqrt{2} \omega t} + \boldsymbol{p}_{3, 3}e^{-i (1 + \sqrt{2}) \omega t} + \boldsymbol{p}_{3, (1, 1)}e^{i 2 \omega t}\right.\\
&& \left. + \boldsymbol{p}_{3, (1, 2)}e^{i (1 + \sqrt{2}) \omega t} + \boldsymbol{p}_{3, (1, 3)}e^{-i \sqrt{2} \omega t}
+ \boldsymbol{p}_{3, (2, 2)} e^{i 2\sqrt{2} \omega t} + \boldsymbol{p}_{3, (2, 3)} e^{-i \omega t}\right.\\
&& \left. + \boldsymbol{p}_{3, (3, 3)} e^{-i (2 + 2\sqrt{2}) \omega t}\right].
\end{eqnarray*}

All this can, of course, be carried forward to larger values of $r$.

\begin{table}[h]
\centering
\caption{Ordered elements of $\mathcal{U}_4$}\label{tab:table}
\begin{tabular}{ l | l | l}
 $m$ & $\sigma_m$ & $\rho^m_{\ell}$ \\ \hline
 $0$ & $\sigma_0 = 0$ & $\rho^0_{0,0,0} = 1$,$\rho^0_{1,2,3} = 2$,$\rho^0_{2,1,3} = 2$,$\rho^0_{3,1,2} = 2$,\\
 $1$ & $\sigma_1 = 1$ & $\rho^1_{0,0,1} = 3$,\\
 $2$ & $\sigma_2 = \sqrt{2}$ & $\rho^2_{0,0,2} = 3$,\\
 $3$ & $\sigma_3 = -1-\sqrt{2}$ & $\rho^3_{0,0,3} = 3$,\\
 $(1, 1)$ & $\sigma_{1, 1} = 2$ & $\rho^{1, 1}_{0,1,1} = 3$,\\
 $(1, 2)$ & $\sigma_{1, 2} = 1 + \sqrt{2}$ & $\rho^{1, 2}_{0,1,2} = 6$,\\
 $(1, 3)$ & $\sigma_{1, 3} = -\sqrt{2}$ & $\rho^{1, 3}_{0,1,3} = 6$,\\
 $(2, 2)$ & $\sigma_{2, 2} = 2\sqrt{2}$ & $\rho^{2, 2}_{0,2,2} = 3$,\\
 $(2, 3)$ & $\sigma_{2, 3} = -1$ & $\rho^{2, 3}_{0,2,3} = 6$,\\
 $(3, 3)$ & $\sigma_{3, 3} = - 2 - 2 \sqrt{2}$ & $\rho^{3, 3}_{0,3,3} = 3$,\\
 $(1, 1, 1)$ & $\sigma_{1, 1, 1} = 3$ & $\rho^{1, 1, 1}_{1, 1, 1} = 1$,\\
 $(1, 1, 2)$ & $\sigma_{1, 1, 2} = 2 + \sqrt{2}$ & $\rho^{1, 1, 2}_{1, 1, 2} = 3$,\\
 $(1, 1, 3)$ & $\sigma_{1, 1, 3} = 1 - \sqrt{2}$ & $\rho^{1, 1, 3}_{1, 1, 3} = 3$,\\
 $(1, 2, 2)$ & $\sigma_{1, 2, 2} = 1 + 2\sqrt{2}$ & $\rho^{1, 2, 2}_{1, 2, 2} = 3$,\\
 $(1, 3, 3)$ & $\sigma_{1, 3, 3} = - 1 - 2\sqrt{2}$ & $\rho^{1, 3, 3}_{1, 3, 3} = 3$,\\
 $(2, 2, 2)$ & $\sigma_{2, 2, 2} = 3\sqrt{2}$ & $\rho^{2, 2, 2}_{2, 2, 2} = 1$,\\
 $(2, 2, 3)$ & $\sigma_{2, 2, 3} = -1 +\sqrt{2}$ & $\rho^{2, 2, 3}_{2, 2, 3} = 3$,\\
 $(2, 3, 3)$ & $\sigma_{2, 3, 3} = -2 -\sqrt{2}$ & $\rho^{2, 3, 3}_{2, 3, 3} = 3$,\\
 $(3, 3, 3)$ & $\sigma_{3, 3, 3} = -3 -3\sqrt{2}$ & $\rho^{3, 3, 3}_{3, 3, 3} = 1$,\\
\end{tabular}
\end{table}

\subsection{Two non-commensurate frequencies}

An interesting special case is $M = 2$, where, without loss of generality, $\kappa_1 \neq 0$
is rational and $\kappa_2$ irrational: this means that the only integer solution to
$m_1 \kappa_1 + m_2 \kappa_2 = 0$ is $m1 = m2 = 0$. This simplifies the argument a great deal.

Simple calculation now confirms that
\begin{eqnarray*}
&& \boldsymbol{p}'_{0, 0} = \boldsymbol{f}(\boldsymbol{p}_{0, 0}), \quad  \boldsymbol{p}_{0, 0}(0) = \boldsymbol{y}_0,\\
&& \qquad \qquad \boldsymbol{p}_{1, m} = \frac{\boldsymbol{a}_m}{i \kappa_m}, \quad m  = 1, 2,\\
&& \boldsymbol{p}'_{1, 0} = \boldsymbol{f}_1(\boldsymbol{p}_{0, 0})[\boldsymbol{p}_{1, 0}], \quad t \geq 0, \quad \boldsymbol{p}_{1, 0}(0) = - \boldsymbol{p}_{1, 1}(0) - \boldsymbol{p}_{1, 2}(0),\\
&& \qquad \qquad \boldsymbol{p}_{2, m} = \frac{1}{i \kappa_m} \{\boldsymbol{f}_1(\boldsymbol{p}_{0, 0})[\boldsymbol{p}_{1, m}] - \boldsymbol{p}'_{1, m}\}, \quad m  = 1, 2,\\
&& \boldsymbol{p}'_{2, 0} = \boldsymbol{f}_1(\boldsymbol{p}_{0, 0})[\boldsymbol{p}_{2, 0}] + \frac{1}{2}\boldsymbol{f}_2(\boldsymbol{p}_{0, 0})[\boldsymbol{p}_{1, 0}, \boldsymbol{p}_{1, 0}], \quad \boldsymbol{p}_{2, 0}(0) = - \boldsymbol{p}_{2, 1}(0) - \boldsymbol{p}_{2, 2}(0),\\
&& \quad \boldsymbol{p}_{3, m} = \frac{1}{i \kappa_m} \{\boldsymbol{f}_1(\boldsymbol{p}_{0, 0})[\boldsymbol{p}_{2, m}] - \boldsymbol{p}'_{2, m} + \boldsymbol{f}_2(\boldsymbol{p}_{0, 0})[\boldsymbol{p}_{1, 0}, \boldsymbol{p}_{1, m}]\}, \quad m  = 1, 2,\\
&& \quad \boldsymbol{p}_{3, (m, m)} = \frac{1}{4 i \kappa_m} \boldsymbol{f}_2(\boldsymbol{p}_{0, 0})[\boldsymbol{p}_{1, m}, \boldsymbol{p}_{1, m}], \quad m  = 1, 2,\\
&& \quad \boldsymbol{p}_{3, (1, 2)} = \frac{1}{i (\kappa_1 + \kappa_2)} \boldsymbol{f}_2(\boldsymbol{p}_{0, 0})[\boldsymbol{p}_{1, 1}, \boldsymbol{p}_{1, 2}],\\
&& \boldsymbol{p}'_{3, 0} = \boldsymbol{f}_1(\boldsymbol{p}_{0, 0})[\boldsymbol{p}_{3, 0}] + \boldsymbol{f}_2(\boldsymbol{p}_{0, 0})[\boldsymbol{p}_{1, 0}, \boldsymbol{p}_{2, 0}] + \frac{1}{6}\boldsymbol{f}_3(\boldsymbol{p}_{0, 0})[\boldsymbol{p}_{1, 0}, \boldsymbol{p}_{1, 0}, \boldsymbol{p}_{1, 0}], \\
&& \quad \boldsymbol{p}_{3, 0}(0) = - \boldsymbol{p}_{3, 1}(0) - \boldsymbol{p}_{3, 2}(0) - \boldsymbol{p}_{3, (1, 1)}(0) - \boldsymbol{p}_{3, (1, 2)}(0) - \boldsymbol{p}_{3, (2, 2)}(0).
\end{eqnarray*}

The next `generation' is
\begin{eqnarray*}
&&\boldsymbol{p}_{4, m} = \frac{1}{i \kappa_m} \left\{\boldsymbol{f}_1(\boldsymbol{p}_{0, 0})[\boldsymbol{p}_{3, m}] - \boldsymbol{p}'_{3, m} + \boldsymbol{f}_2(\boldsymbol{p}_{0, 0})[\boldsymbol{p}_{1, 0}, \boldsymbol{p}_{2, m}]\right.\\
&& \quad \left. + \boldsymbol{f}_2(\boldsymbol{p}_{0, 0})[\boldsymbol{p}_{1, m}, \boldsymbol{p}_{2, 0}] + \frac{1}{6} \boldsymbol{f}_3(\boldsymbol{p}_{0, 0})[\boldsymbol{p}_{1, 0}, \boldsymbol{p}_{1, 0}, \boldsymbol{p}_{1, 0}]\right\}, \quad m  = 1, 2,\\
&&\boldsymbol{p}_{4, (m, m)} = \frac{1}{2 i \kappa_m} \left\{\boldsymbol{f}_1(\boldsymbol{p}_{0, 0})[\boldsymbol{p}_{3, (m, m)}] + \boldsymbol{f}_2(\boldsymbol{p}_{0, 0})[\boldsymbol{p}_{1, m}, \boldsymbol{p}_{2, m}]\right.\\
&& \quad \left. + \frac{1}{2} \boldsymbol{f}_3(\boldsymbol{p}_{0, 0})[\boldsymbol{p}_{1, 0}, \boldsymbol{p}_{1, m}, \boldsymbol{p}_{1, m}]\right\}, \quad m  = 1, 2,\\
&&\boldsymbol{p}_{4, (1, 2)} = \frac{1}{i (\kappa_1 + \kappa_2)} \left\{\boldsymbol{f}_1(\boldsymbol{p}_{0, 0})[\boldsymbol{p}_{3, (1, 2)}] + \boldsymbol{f}_2(\boldsymbol{p}_{0, 0})[\boldsymbol{p}_{1, 1}, \boldsymbol{p}_{2, 2}]\right.\\
&& \quad \left. + \boldsymbol{f}_2(\boldsymbol{p}_{0, 0})[\boldsymbol{p}_{1, 2}, \boldsymbol{p}_{2, 1}] + \boldsymbol{f}_3(\boldsymbol{p}_{0, 0})[\boldsymbol{p}_{1, 0}, \boldsymbol{p}_{1, 1}, \boldsymbol{p}_{1, 2}]\right\},\\
&&\boldsymbol{p}_{4, (m, m, m)} = \frac{1}{18 i \kappa_m} \boldsymbol{f}_3(\boldsymbol{p}_{0, 0})[\boldsymbol{p}_{1, m}, \boldsymbol{p}_{1, m}, \boldsymbol{p}_{1, m}], \quad m = 1, 2,\\
&&\boldsymbol{p}_{4, (1, 1, 2)} = \frac{1}{2 i (2 \kappa_1 + \kappa_2)} \boldsymbol{f}_3(\boldsymbol{p}_{0, 0})[\boldsymbol{p}_{1, 1}, \boldsymbol{p}_{1, 1}, \boldsymbol{p}_{1, 2}], \\
&&\boldsymbol{p}_{4, (1, 2, 2)} = \frac{1}{2 i (\kappa_1 + 2 \kappa_2)} \boldsymbol{f}_3(\boldsymbol{p}_{0, 0})[\boldsymbol{p}_{1, 1}, \boldsymbol{p}_{1, 2}, \boldsymbol{p}_{1, 2}]
\end{eqnarray*}
and
\begin{eqnarray*}
&& \boldsymbol{p}'_{4, 0} = \boldsymbol{f}_1(\boldsymbol{p}_{0, 0})[\boldsymbol{p}_{4, 0}] + \boldsymbol{f}_2(\boldsymbol{p}_{0, 0})[\boldsymbol{p}_{1, 0}, \boldsymbol{p}_{3, 0}] + \frac{1}{2}\boldsymbol{f}_2(\boldsymbol{p}_{0, 0})[\boldsymbol{p}_{2, 0}, \boldsymbol{p}_{2, 0}]\\
&& + \frac{1}{2}\boldsymbol{f}_3(\boldsymbol{p}_{0, 0})[\boldsymbol{p}_{1, 0}, \boldsymbol{p}_{1, 0}, \boldsymbol{p}_{2, 0}] + \frac{1}{24}
\boldsymbol{f}_4(\boldsymbol{p}_{0, 0})[\boldsymbol{p}_{1, 0}, \boldsymbol{p}_{1, 0}, \boldsymbol{p}_{1, 0}, \boldsymbol{p}_{1, 0}],\\
&& \quad \boldsymbol{p}_{4, 0}(0) = - \boldsymbol{p}_{4, 1}(0) - \boldsymbol{p}_{4, 2}(0) - \boldsymbol{p}_{4, (1, 1)}(0) - \boldsymbol{p}_{4, (1, 2)}(0) - \boldsymbol{p}_{4, (2, 2)}(0)\\
&& \quad - \boldsymbol{p}_{4, (1, 1, 1)}(0) - \boldsymbol{p}_{4, (1, 1, 2)}(0) - \boldsymbol{p}_{4, (1, 2, 2)}(0) - \boldsymbol{p}_{4, (2, 2, 2)}(0).
\end{eqnarray*}
Greater, but not
insurmountable effort is required to develop a general asymptotic expansion
in this case. However, to all intents and purposes, expanding up to $r = 4$ is
sufficient to derive an exceedingly accurate solution.

\subsection{Comments}

\textbf{Comment 1}: As we increase levels $r$, we are increasingly likely to encounter
the well-known phenomenon of small denominators (Verhulst 1990): unless
$\kappa_m = c \psi_m$, $m = 1, \cdots, M$, where all the  $\psi_m$s are rational (in which case,
replacing $\omega$ by its product with the least common denominator of the  $\psi_m$s,
we are back to the case (1.2) of frequencies being integer multiples of $\omega$) and
positive, the set of all finite-length linear combinations of the $\kappa_m$s is dense in
$\mathbb{R}$ (Besicovitch 1932). In particular, we can approach $0$ arbitrarily close by
such linear combinations. This is different from $\kappa_m$s summing up exactly to
zero: as demonstrated in Subsection 3.4, we can deal with the latter problem
but not with the denominators in (3.7) or (3.8) becoming arbitrarily small in
magnitude. Like with other averaging techniques, there is no simple remedy
to this phenomenon. This restricts the range of r at which the asymptotic
expansion is effective. Having said so, and bearing in mind that the truncation
of (1.7) to $r \leq R$ yields an error of $O(\omega^{-R-1})$ and, $|\omega|$ being large, we are
likely to obtain very high accuracy before small denominators kick in. Hence,
this phenomenon has little practical implications.

Incidentally, this is precisely the reason for the requirement that, unlike in
(1.2), the number of initial frequencies is finite. Otherwise, we could have encountered
small denominators already for $r = 3$ and this would have definitely
placed genuine restrictions on the applicability of our approach.

\textbf{Comment 2}: There is an alternative to our expansion. We may decide that
the $\kappa_m$s are symbols, rather than specific numbers. Not being assigned specific
values, it is meaningless to talk about sets $\mathcal{W}^n_{r, m}$ because $\kappa_{\ell_1} + \kappa_{\ell_2} = \kappa_m$, say,
has no meaning (except when $\ell_1 = 0$, $\ell_2 = m$). Of course, in that case we may
have several distinct terms which correspond to the same frequency, once we
allocate specific values to the $\kappa_m$s, but the \textit{quid pro quo} is considerable simplification
and no multiplicities (which depend on specific values of $\kappa_m$s, hence
need be re-evaluated each time we have new frequencies). Unfortunately, this
approach has another, more critical, shortcoming. We must identify all linear
combinations of $\kappa_m$s that sum up to zero, because they require an altogether
different treatment.

It would have been possible to proceed differently, by separating all the
terms in of the form $\sum_{j = 1}^{s} \kappa_{\ell_j}$ into two subsets: those that sum to zero and
those that are nonzero. We do not need to specify which is which - this becomes
apparent only once values are allocated - just to remember that the nonzero
sums give rise to new frequencies, with coefficients derived by recursion, while
zero sums are lumped into a differential equation for a non-oscillatory term.
While this is certainly feasible, it seems that the current approach is probably
simpler and more transparent.

\section{Numerical experiments}

In the current section we present two examples that illustrate the construction
of our expansions and demonstrate the effectiveness of our approach. In
each case we compare the pointwise error incurred by a truncated expansion
(1.7) with either the exact solution or the Maple routine \texttt{rkf45} with exceedingly
high error tolerance, using 20 significant decimal digits. Specifically, we
measure the components of
\begin{eqnarray*}
\boldsymbol{y}(t) -\boldsymbol{p}_{0, 0}(t) - \sum_{r = 1}^{s} \frac{1}{\omega^r} \sum_{m \in \mathcal{U}_r}
\boldsymbol{p}_{r, m}(t) e^{i \sigma_m \omega t}
\end{eqnarray*}
for different values of $s$.

\subsection{A linear example}

We consider the equation
\begin{eqnarray}
\ddot{x} + \frac{3}{5} \dot{x} + \frac{21}{5} x = t e^{\sqrt{2} i \omega t}  + t^2 e^{- (1 + \sqrt{2}) i \omega t}, \quad t
\geq 0, \quad x(0) = \dot{x}(0) = 0.5.
\end{eqnarray}

Letting $\boldsymbol{y} = [x, \dot{x}]^T$, we reformulate (4.1) as the system
\begin{eqnarray*}
\boldsymbol{y}' = \left(\begin{array}{cc}
 0 & 1 \\
 -\frac{21}{5} & -\frac{3}{5}
\end{array}\right) \boldsymbol{y} + \left(\begin{array}{cc}
0 \\
1
\end{array}\right) \left(t e^{\sqrt{2} i \omega t}  + t^2 e^{- (1 + \sqrt{2}) i \omega t}\right),\quad t \geq 0, \quad
\boldsymbol{y}(0) = \left[\begin{array}{cc}
\frac{1}{2} \\
\frac{1}{2}
\end{array}\right].
\end{eqnarray*}
This being a linear equation, the exact solution and its asymptotic expansion
are available explicitly using the theory from Section 2.

Figures 4.1 and 4.2 display the real part of the error functions in
computing $x$ and $\dot{x}$ , respectively, for $s$ between $0$ and
$4$ within $t \in [0, 5]$. It is clear that each time we increase
$s$, the error indeed decreases substantially, in line with our
theory. \footnote[2] {The fact that Figs 4.1a and 4.1b are identical
is a fluke-anyway, it is evident from the results on Page 5.}

\begin{figure}[htbp]
\centering
{\label{fig:4.1}\includegraphics[width=4cm,height=4cm]{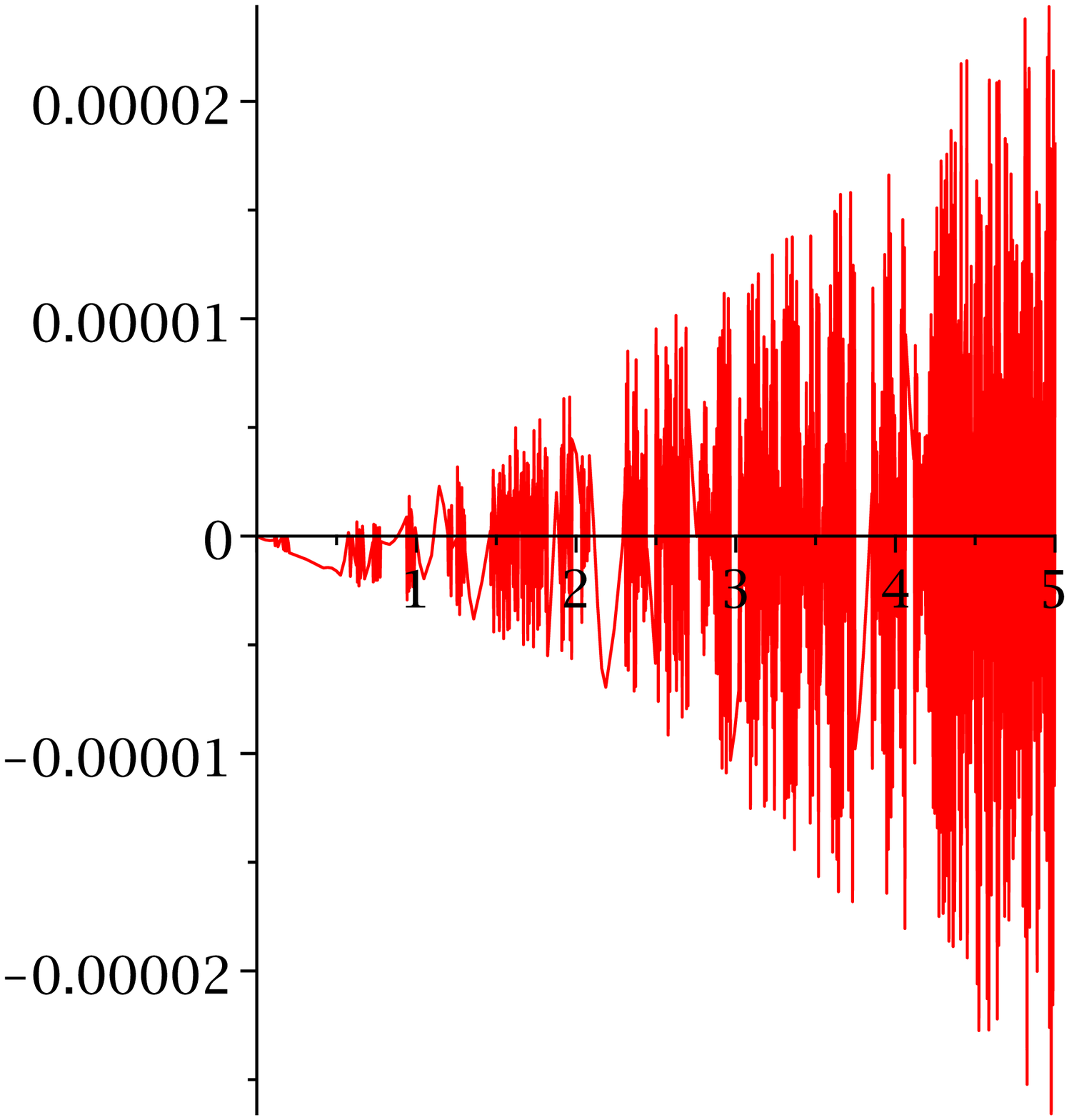}}\qquad\qquad\qquad
{\label{fig:4.1}\includegraphics[width=4cm,height=4cm]{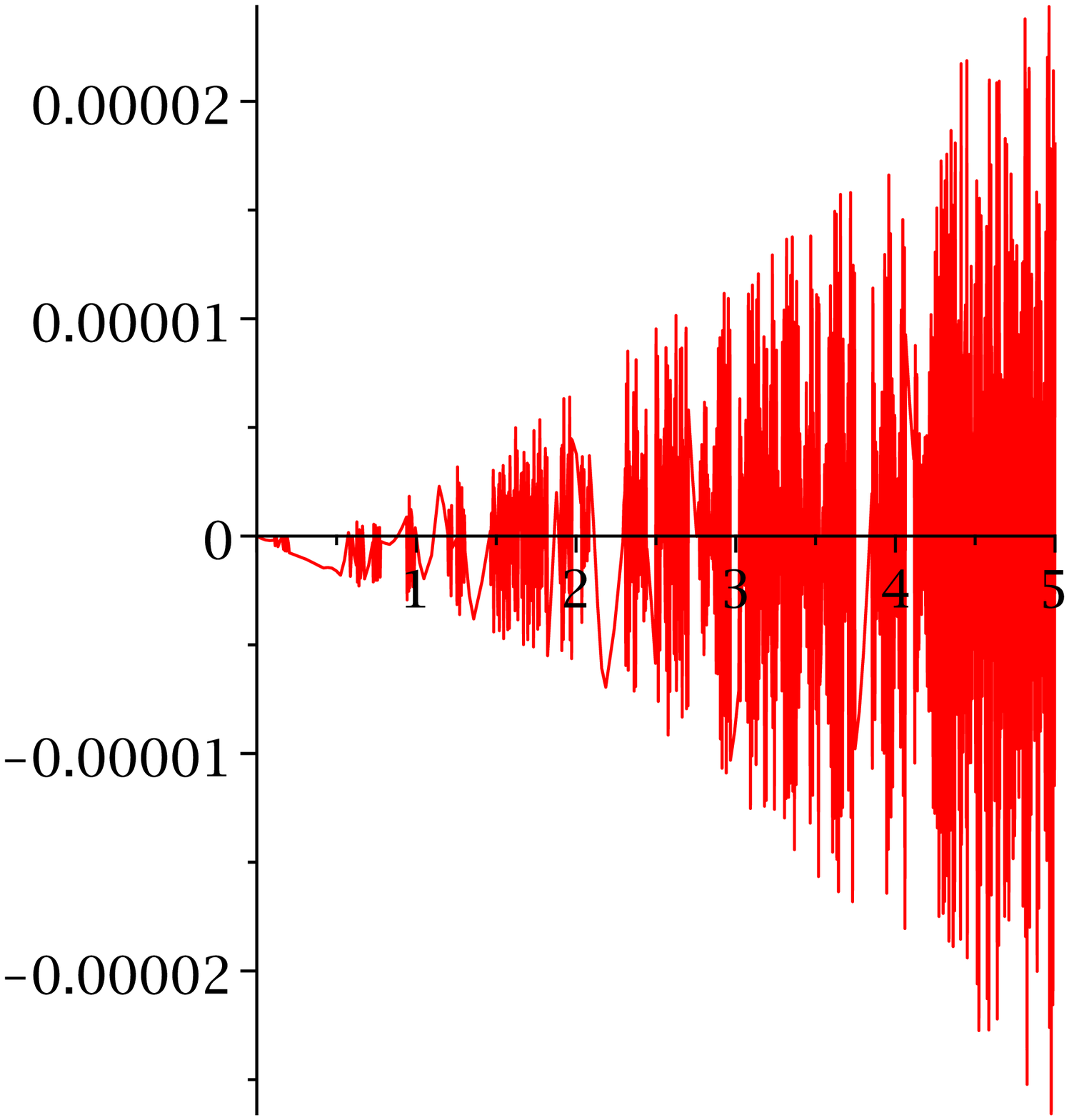}}\\
{\label{fig:4.1}\includegraphics[width=4cm,height=4cm]{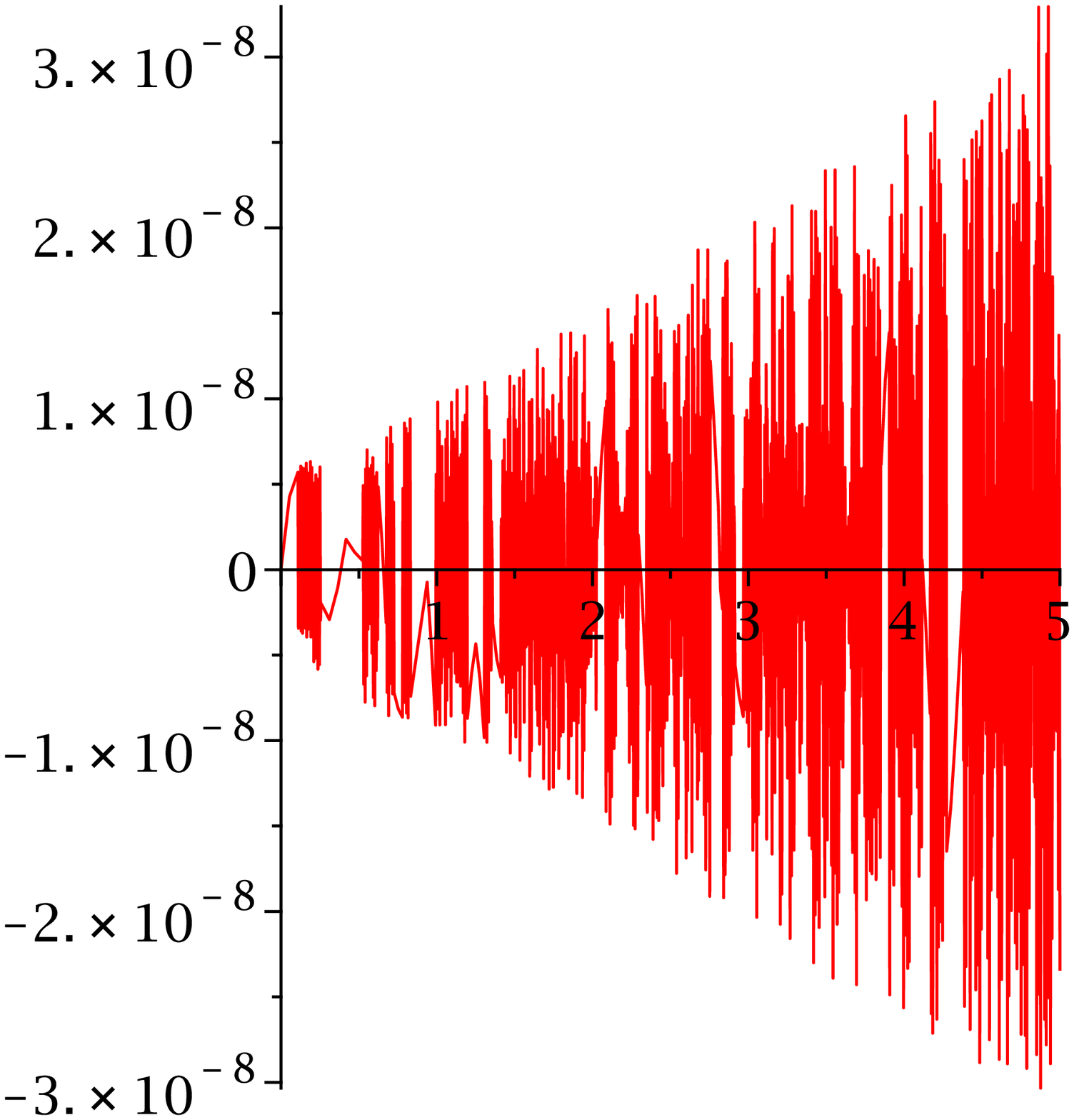}}\qquad\qquad\qquad
{\label{fig:4.1}\includegraphics[width=4cm,height=4cm]{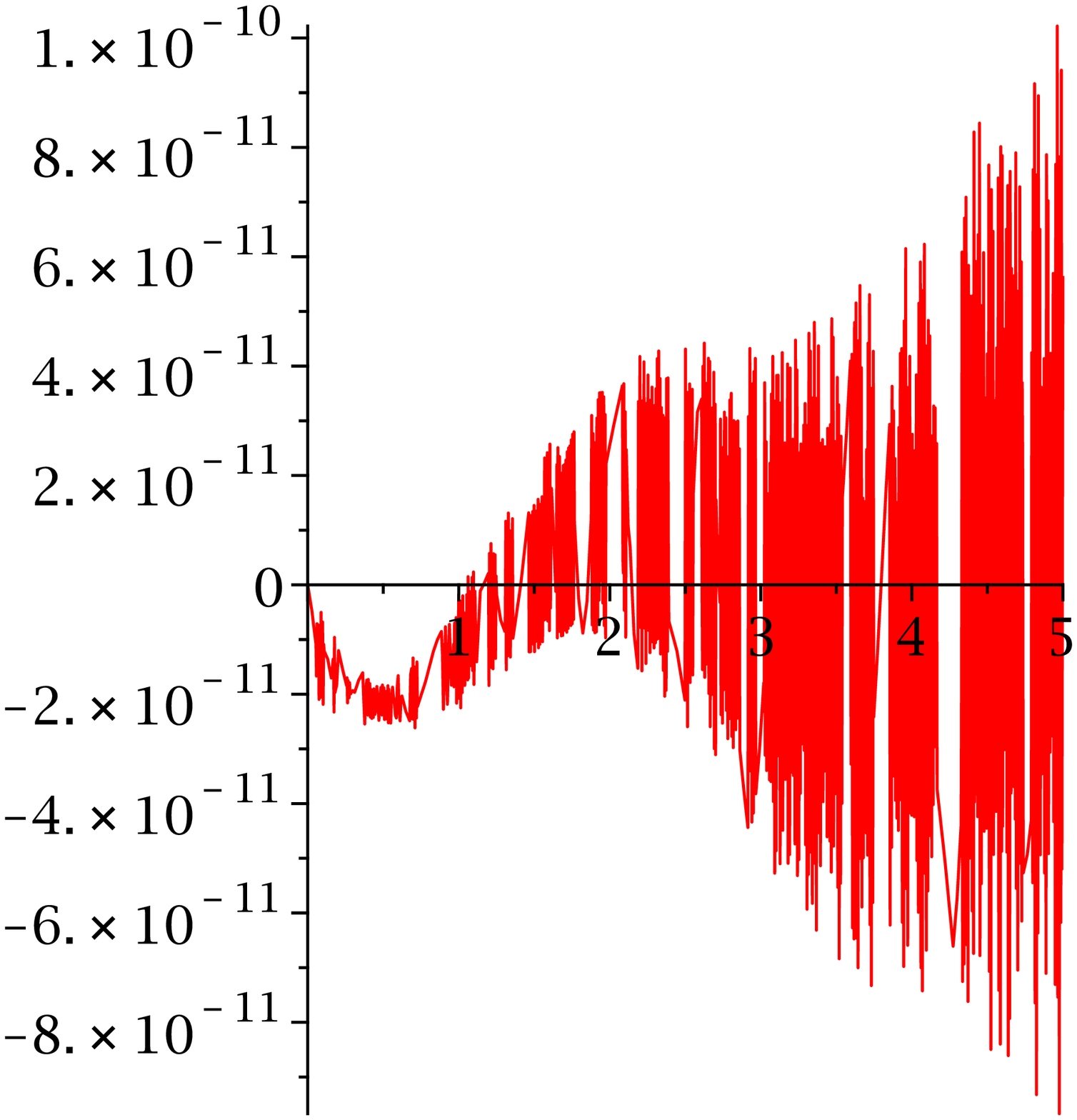}}\\
{\label{fig:4.1}\includegraphics[width=4cm,height=4cm]{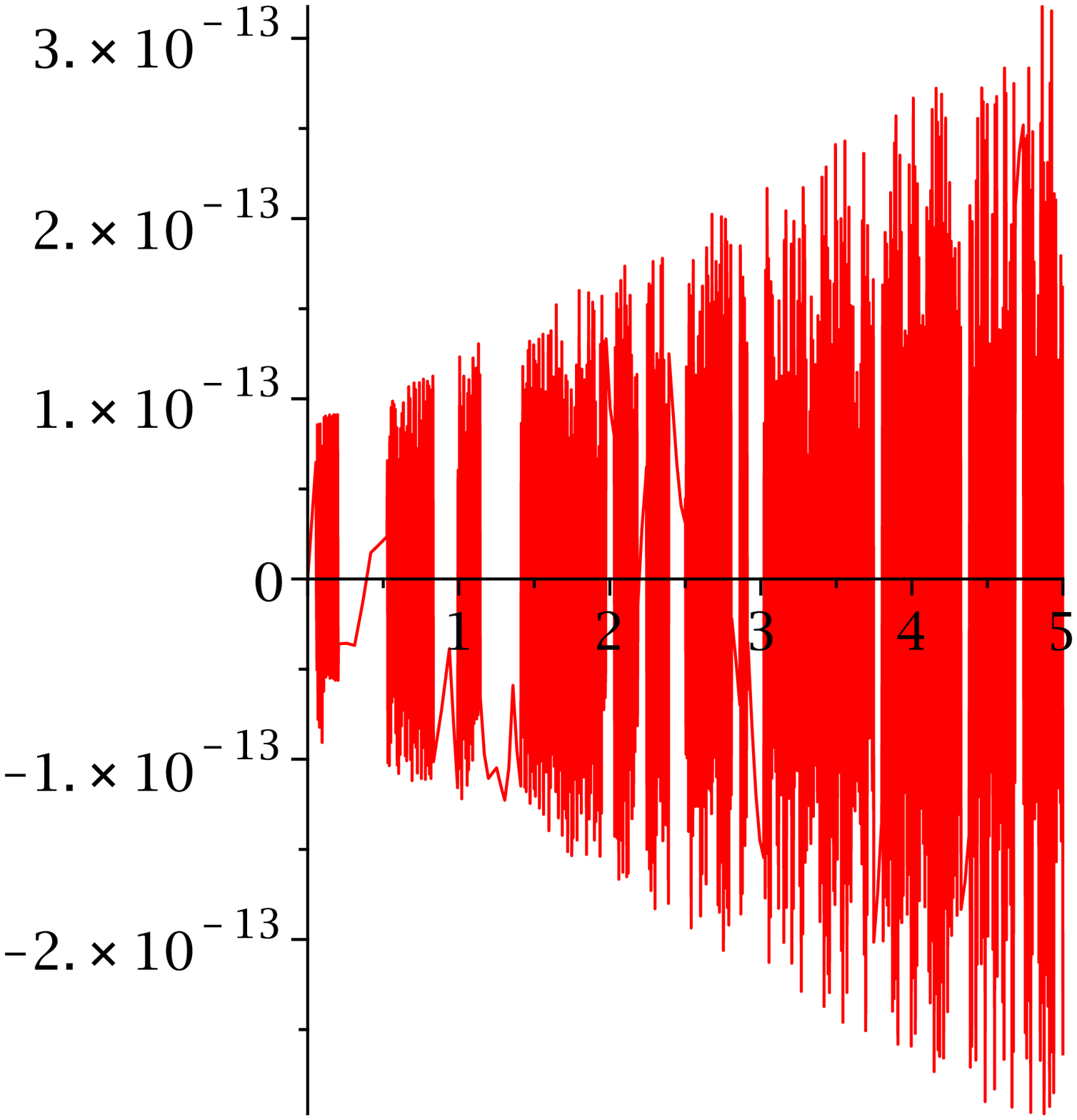}}
\caption{The real part of the error in $x$ committed by the
asymptotic expansion, as applied to the linear system (4.1) with
$\omega = 500$ for $s = 0, 1, 2, 3, 4$ (from top left onwards).}
\end{figure}

\begin{figure}[htbp]
\centering
{\label{fig:4.2}\includegraphics[width=4cm,height=4cm]{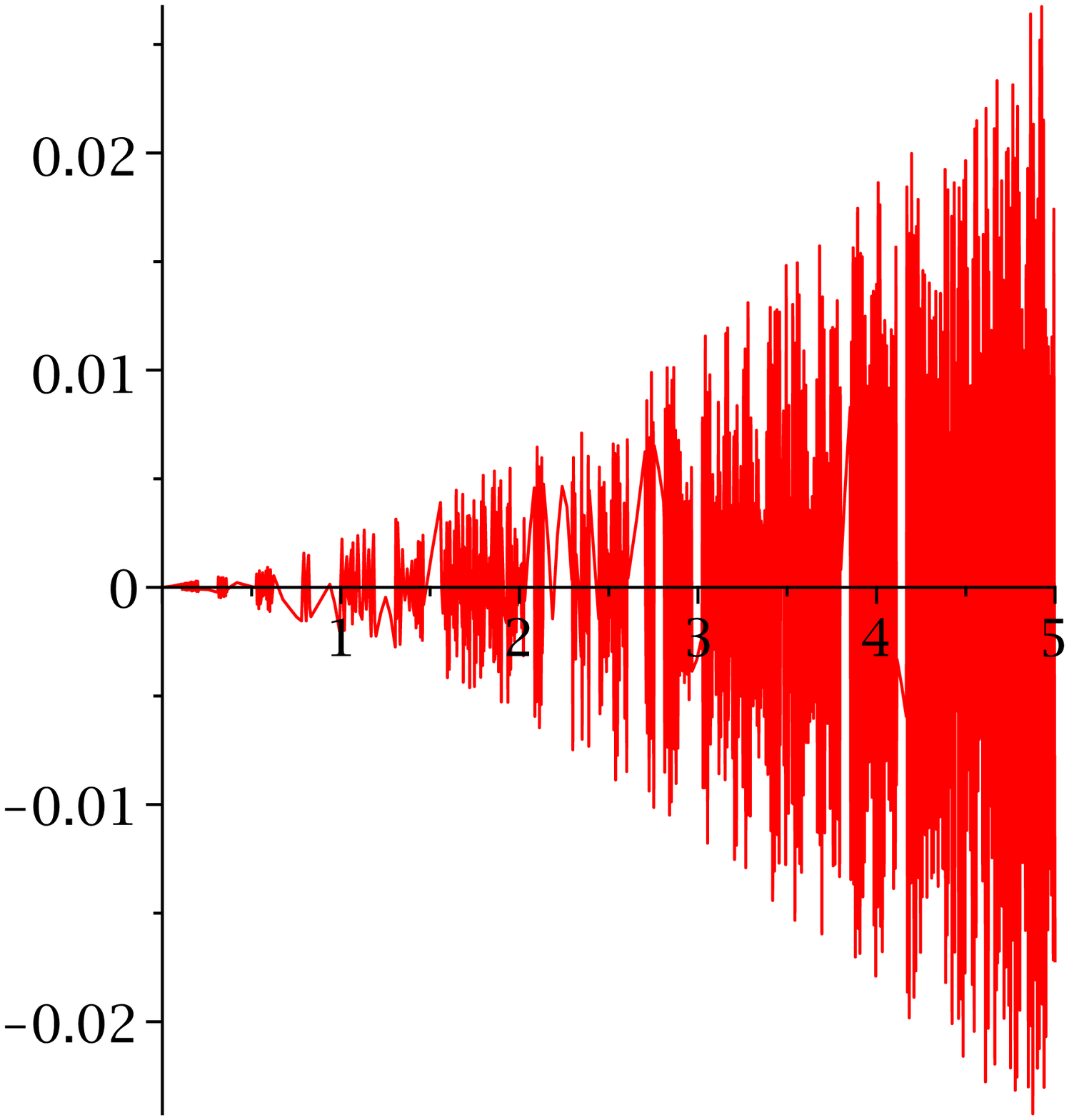}}\qquad\qquad\qquad
{\label{fig:4.2}\includegraphics[width=4cm,height=4cm]{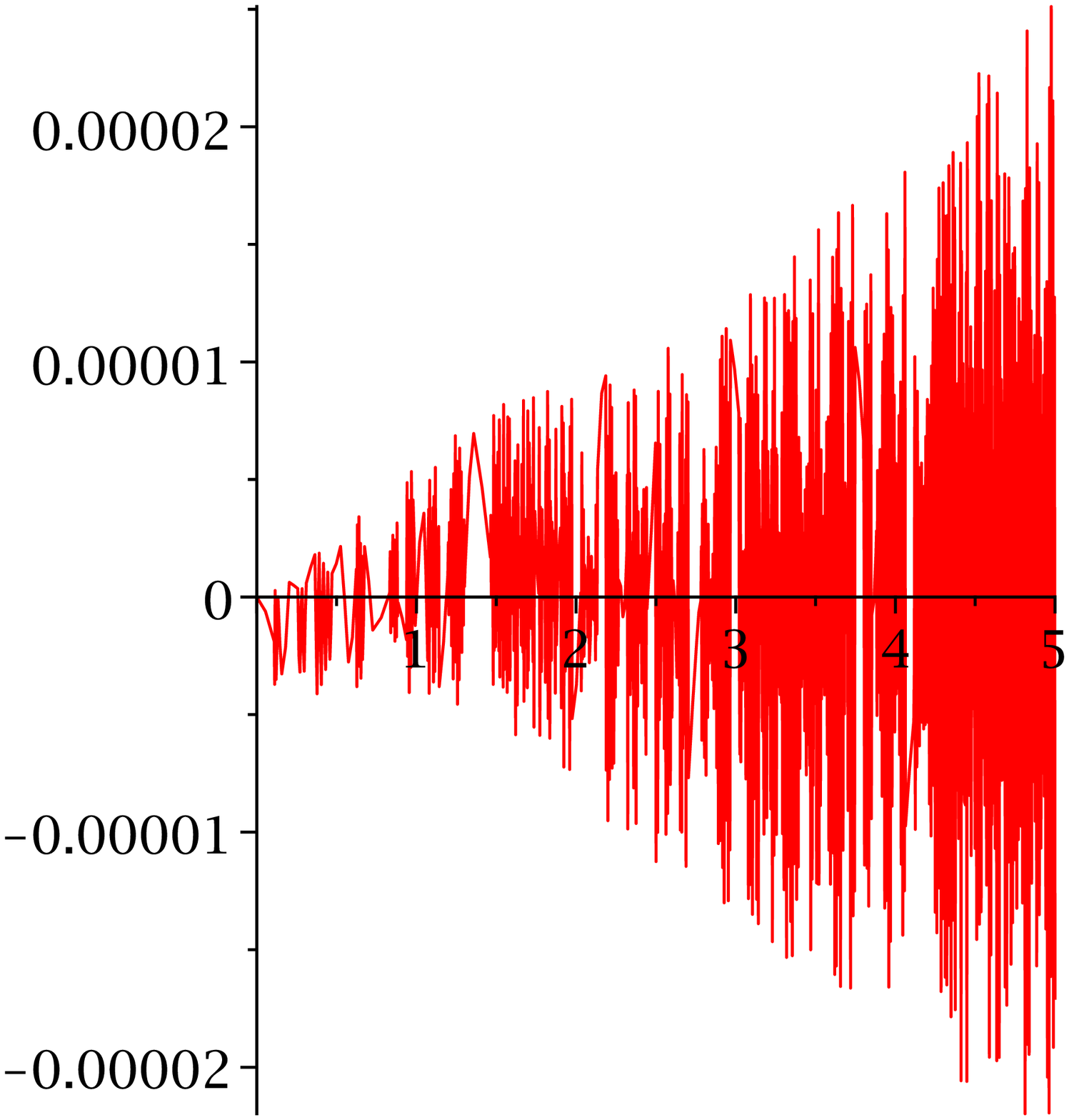}}\\
{\label{fig:4.2}\includegraphics[width=4cm,height=4cm]{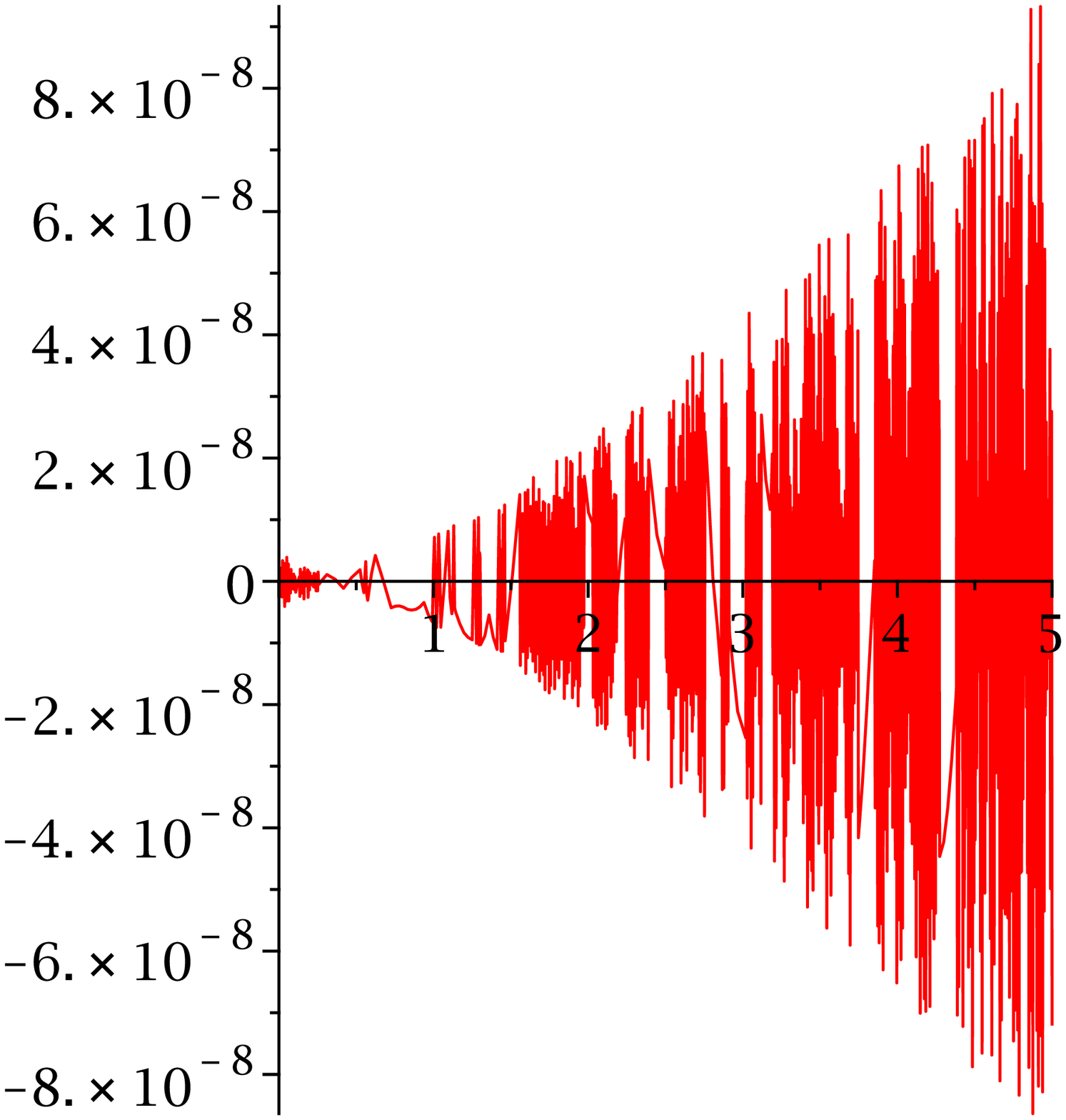}}\qquad\qquad\qquad
{\label{fig:4.2}\includegraphics[width=4cm,height=4cm]{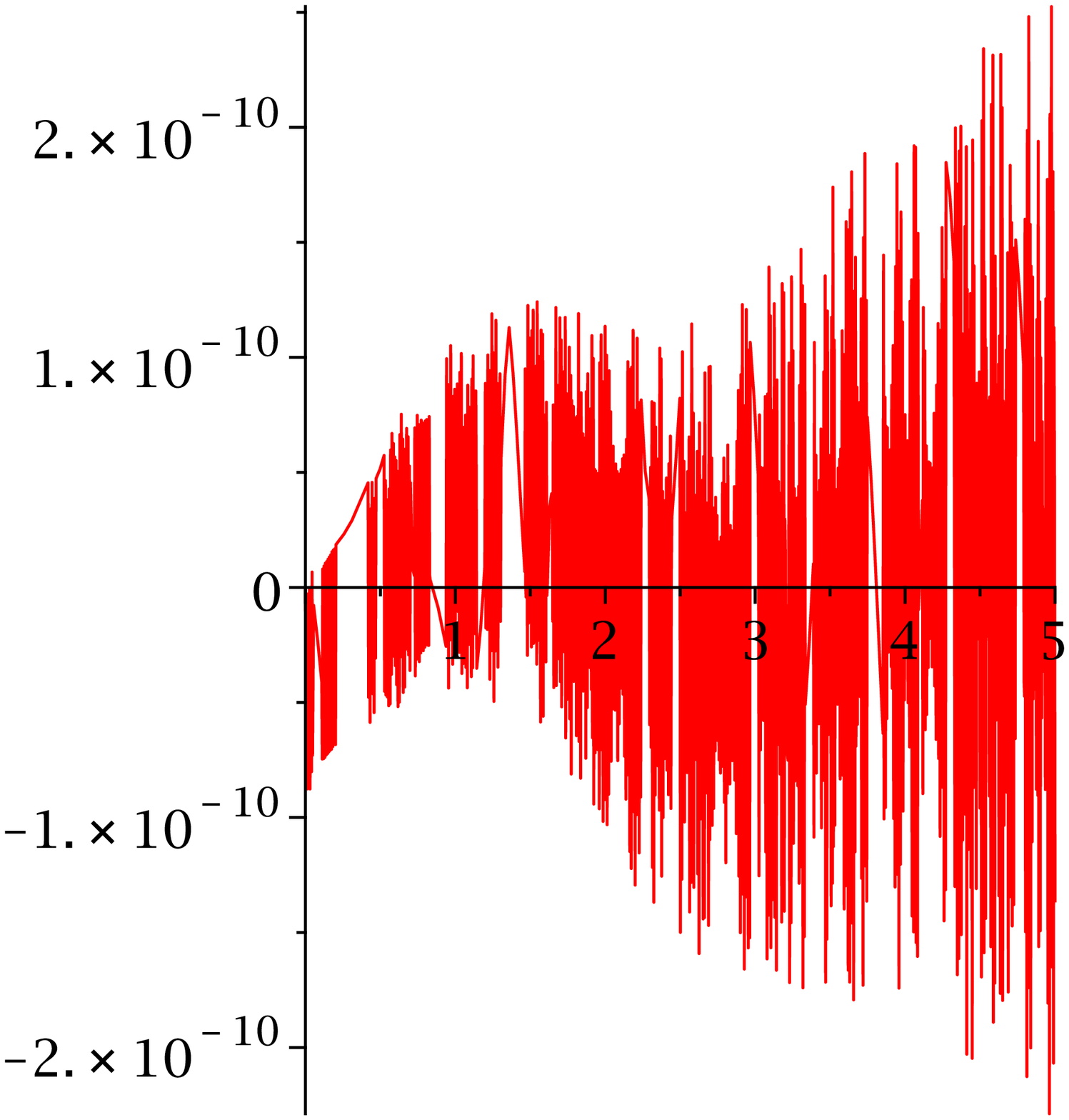}}\\
{\label{fig:4.2}\includegraphics[width=4cm,height=4cm]{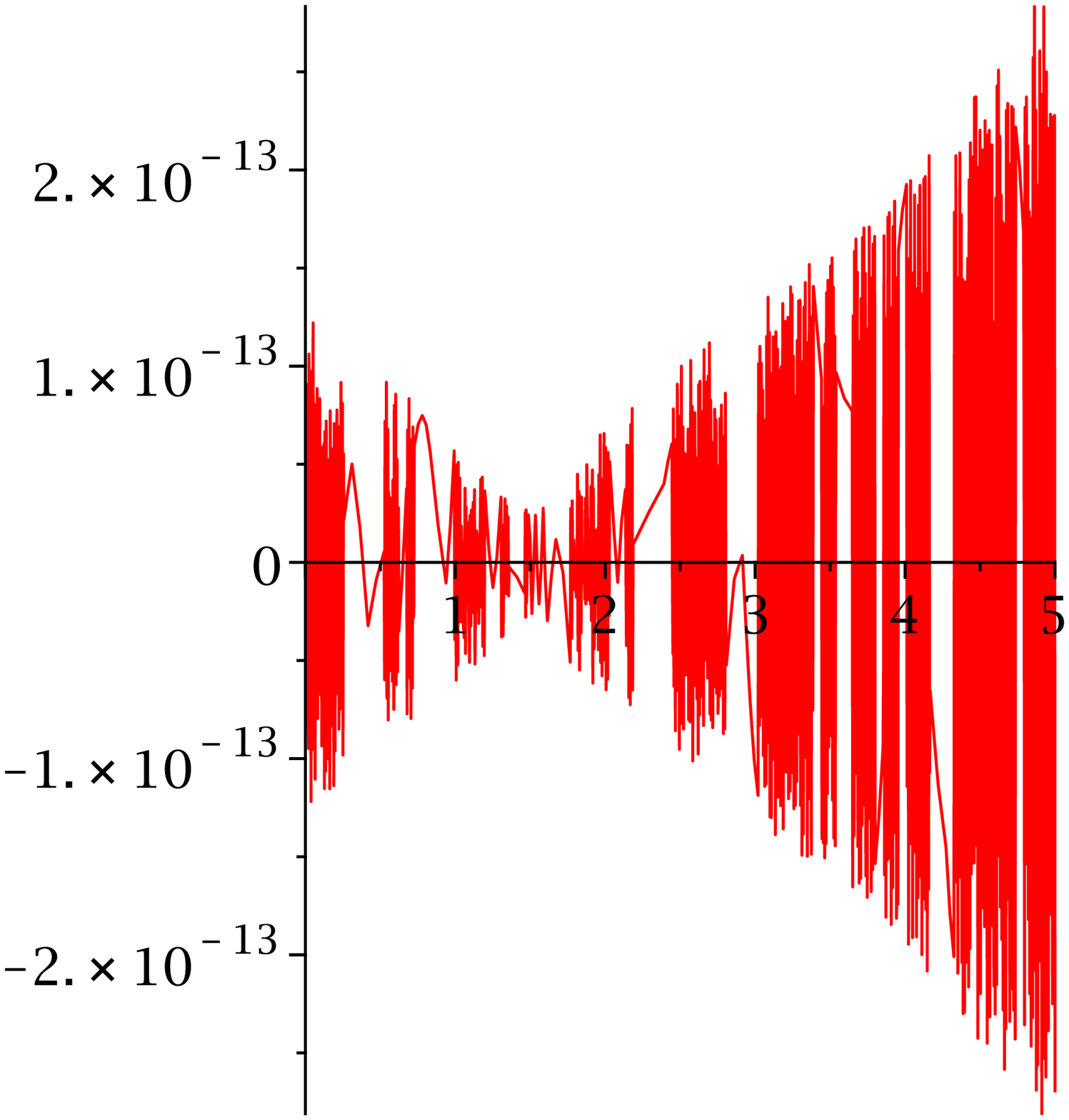}}
\caption{The real part of the error in $\dot{x}$ committed by the asymptotic
expansion, as applied to the linear system (4.1) with $\omega = 500$ for $s = 0, 1, 2, 3, 4$
(from top left onwards). It is shown that the error decreases once $s$ is increased from $0$ to $1$.}
\end{figure}

Identical information is reported in Figs 4.3 and 4.4 for frequency
$\omega = 5000$ within $t \in [0, 5]$. A comparison with the two
previous figures emphasises the important point that the efficiency
of the asymptotic - numerical method grows with $\omega$, while the
cost is to all intents and purposes identical. Indeed, wishing to
produce similar error to our method with $s = 4$, the Maple routine
\texttt{rkf45} needs be applied with absolute and relative error
tolerances of $10^{-13}$ and $10^{-18}$ respectively.
\footnote[3]{Such error tolerances are impossible in Matlab, which
explains our use of Maple.} Although the method is robust enough to
produce correct magnitude of global errors, this comes at a steep
price. Thus, while our method takes less than one second to compute
the solution and requires $\approx 17.5$ kbytes of storage,
\texttt{rkf45} takes $\approx 2740$ seconds to compute the solution
for $\omega = 500$ and requires $\approx 10^7$ kbytes. This
increases to $\approx 4533$ seconds and $\approx 1.6 \times 10^7$
kbytes for $\omega = 5000$.

\begin{figure}[htbp]
\centering
{\label{fig:4.3}\includegraphics[width=4cm,height=4cm]{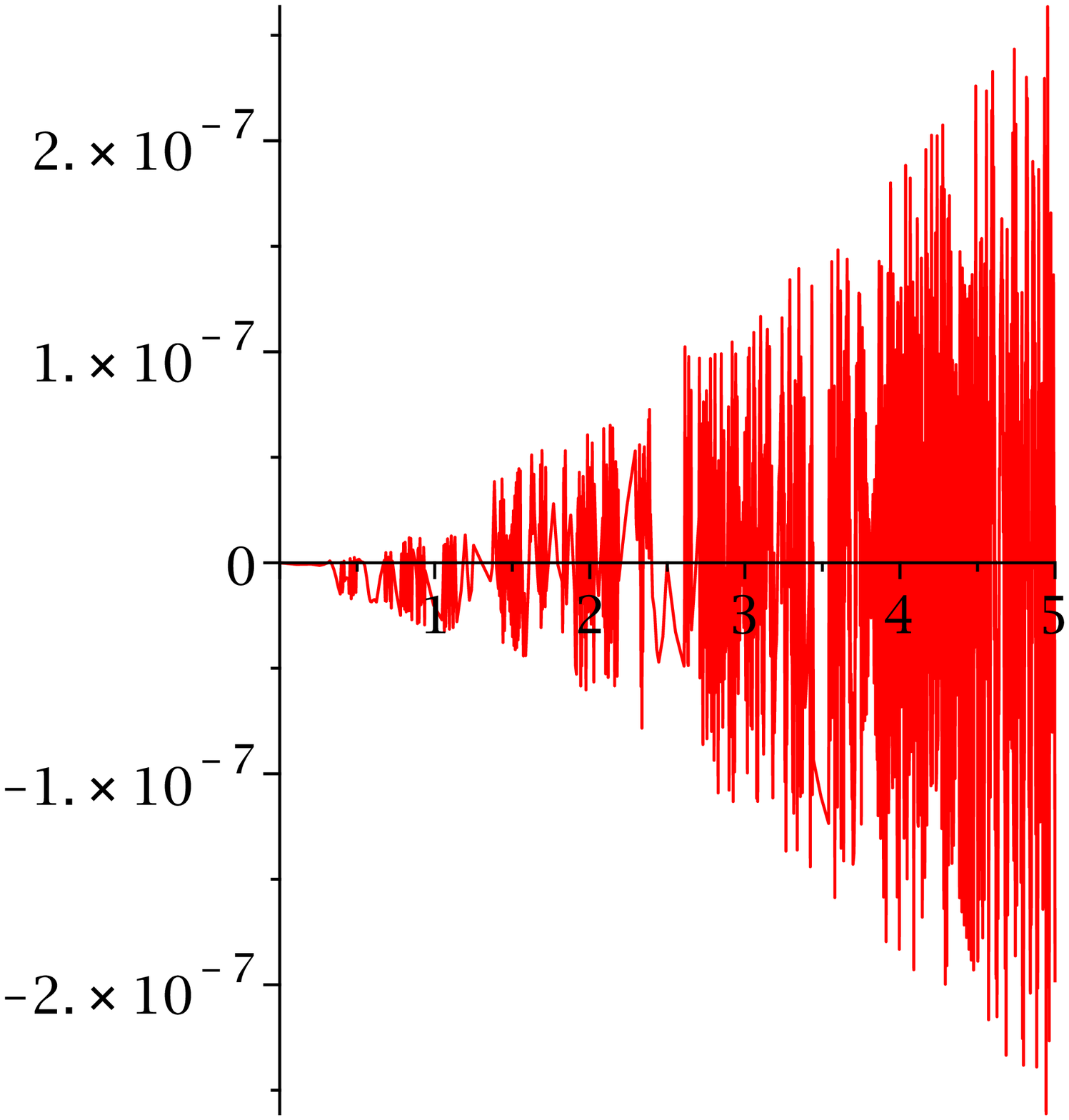}}\qquad\qquad\qquad
{\label{fig:4.3}\includegraphics[width=4cm,height=4cm]{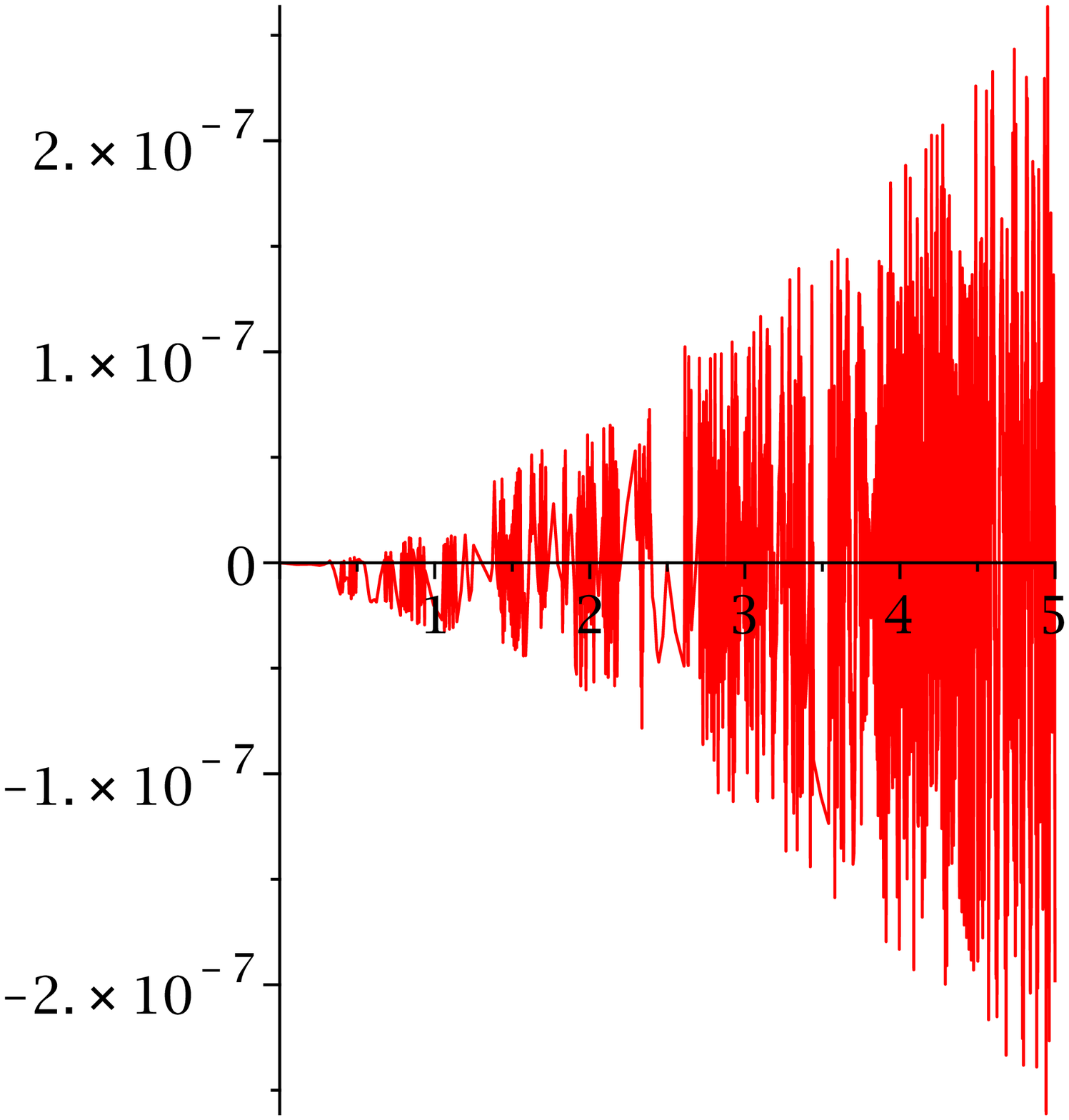}}\\
{\label{fig:4.3}\includegraphics[width=4cm,height=4cm]{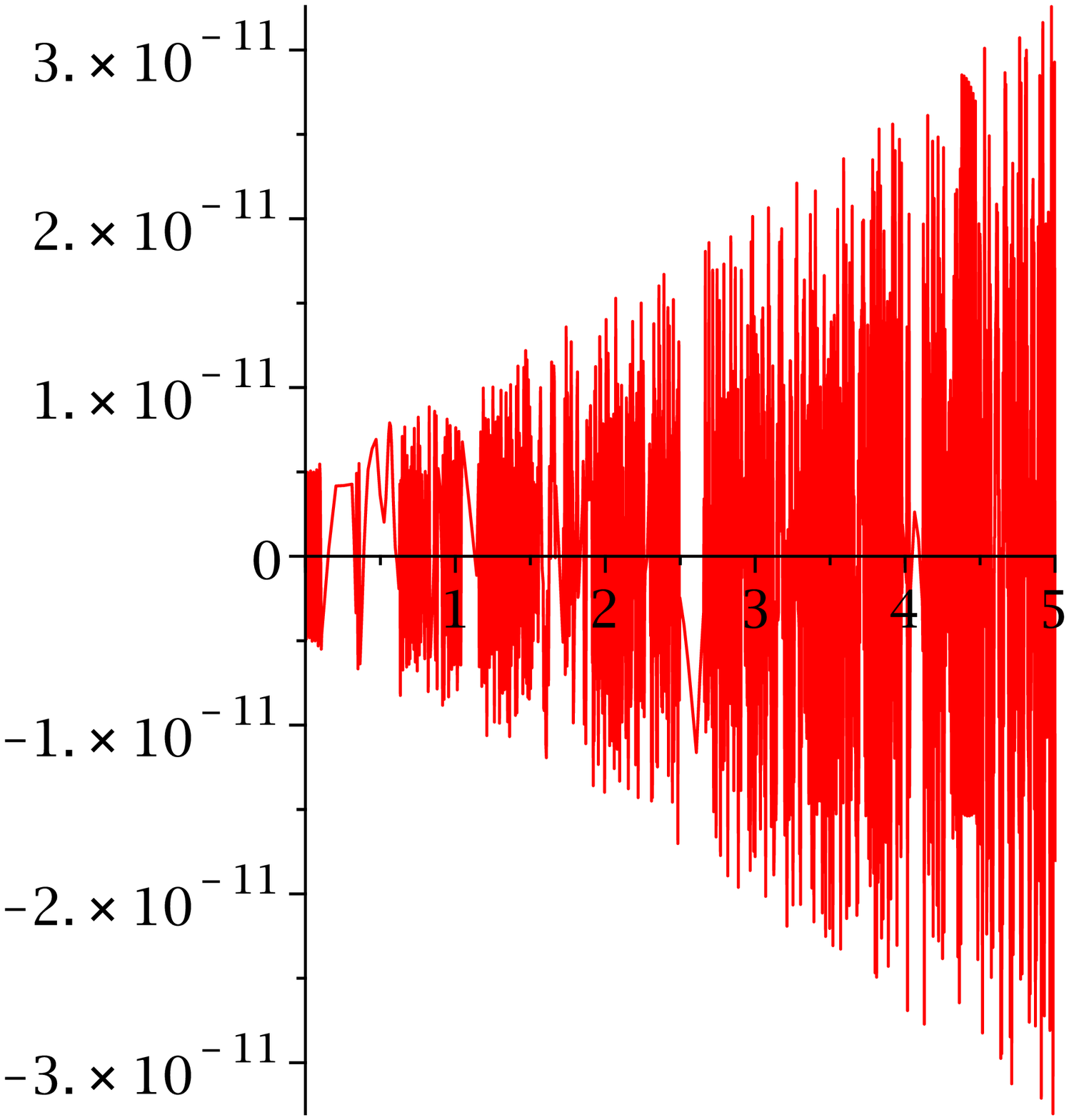}}\qquad\qquad\qquad
{\label{fig:4.3}\includegraphics[width=4cm,height=4cm]{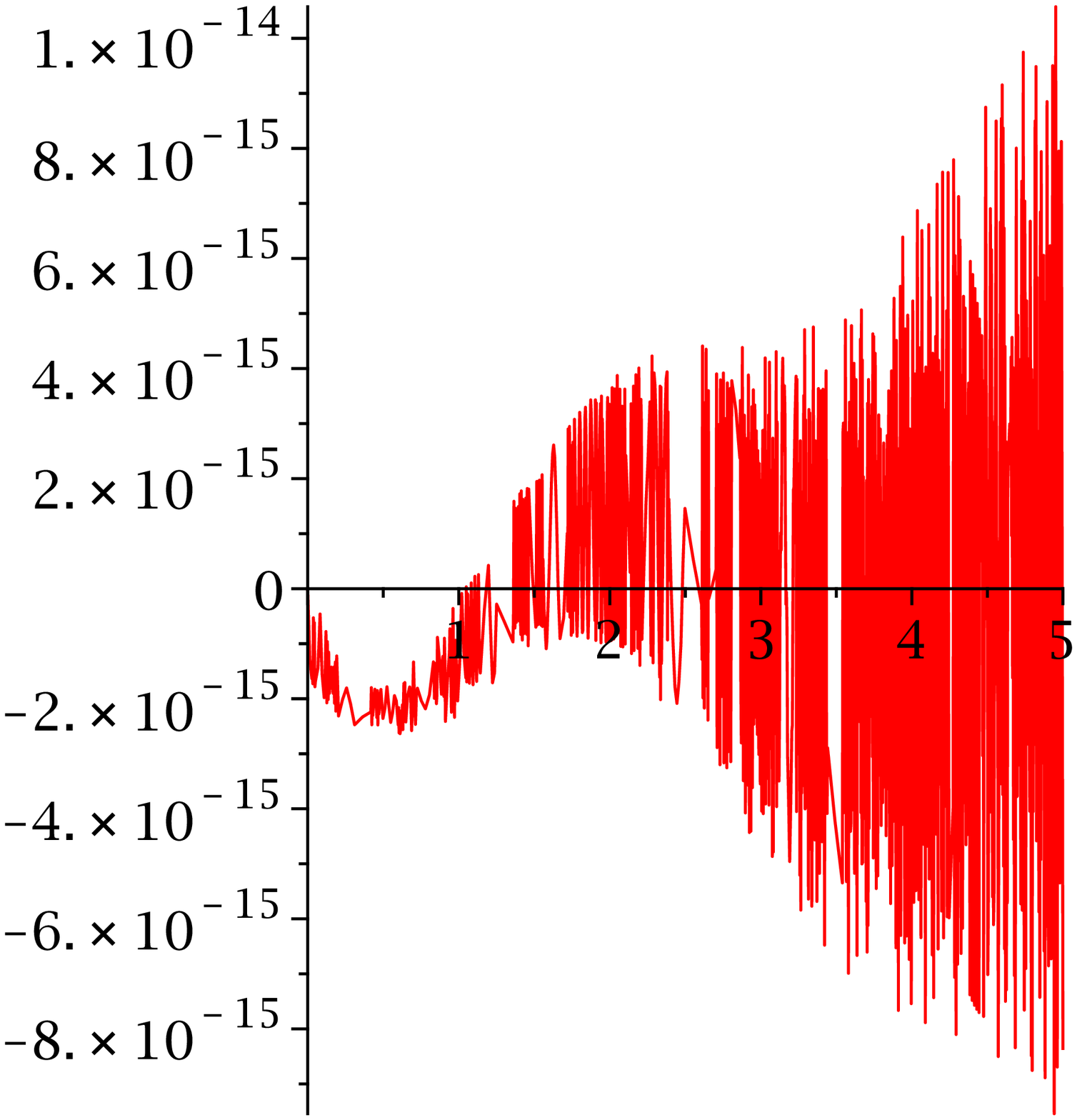}}\\
{\label{fig:4.3}\includegraphics[width=4cm,height=4cm]{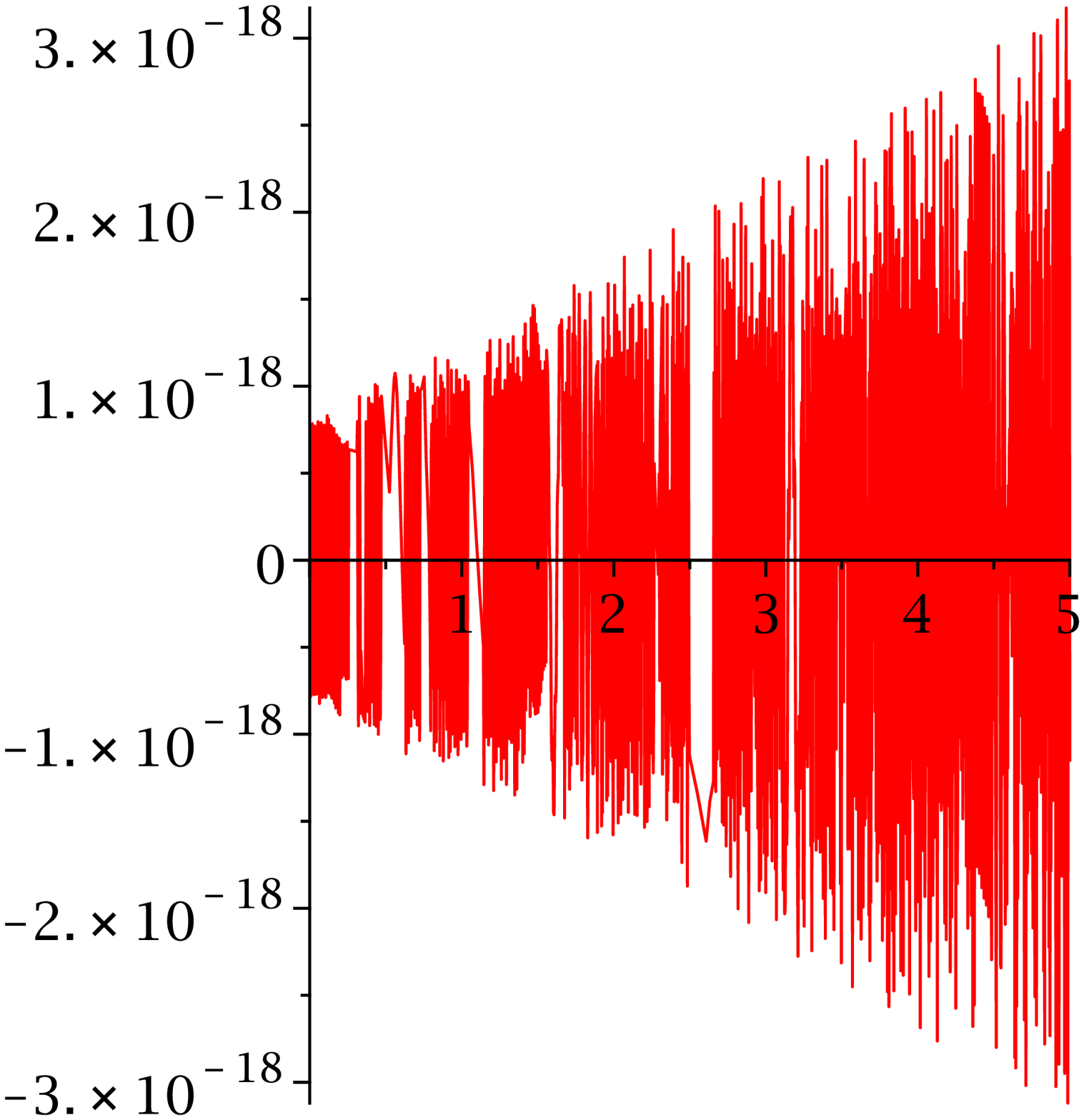}}
\caption{The real part of the error in $x$ committed by the asymptotic expansion,
as applied to the linear system (4.1) with $\omega = 5000$ for $s = 0, 1, 2, 3, 4$
(from top left onwards).}
\end{figure}

\begin{figure}[htbp]
\centering
{\label{fig:re}\includegraphics[width=4cm,height=4cm]{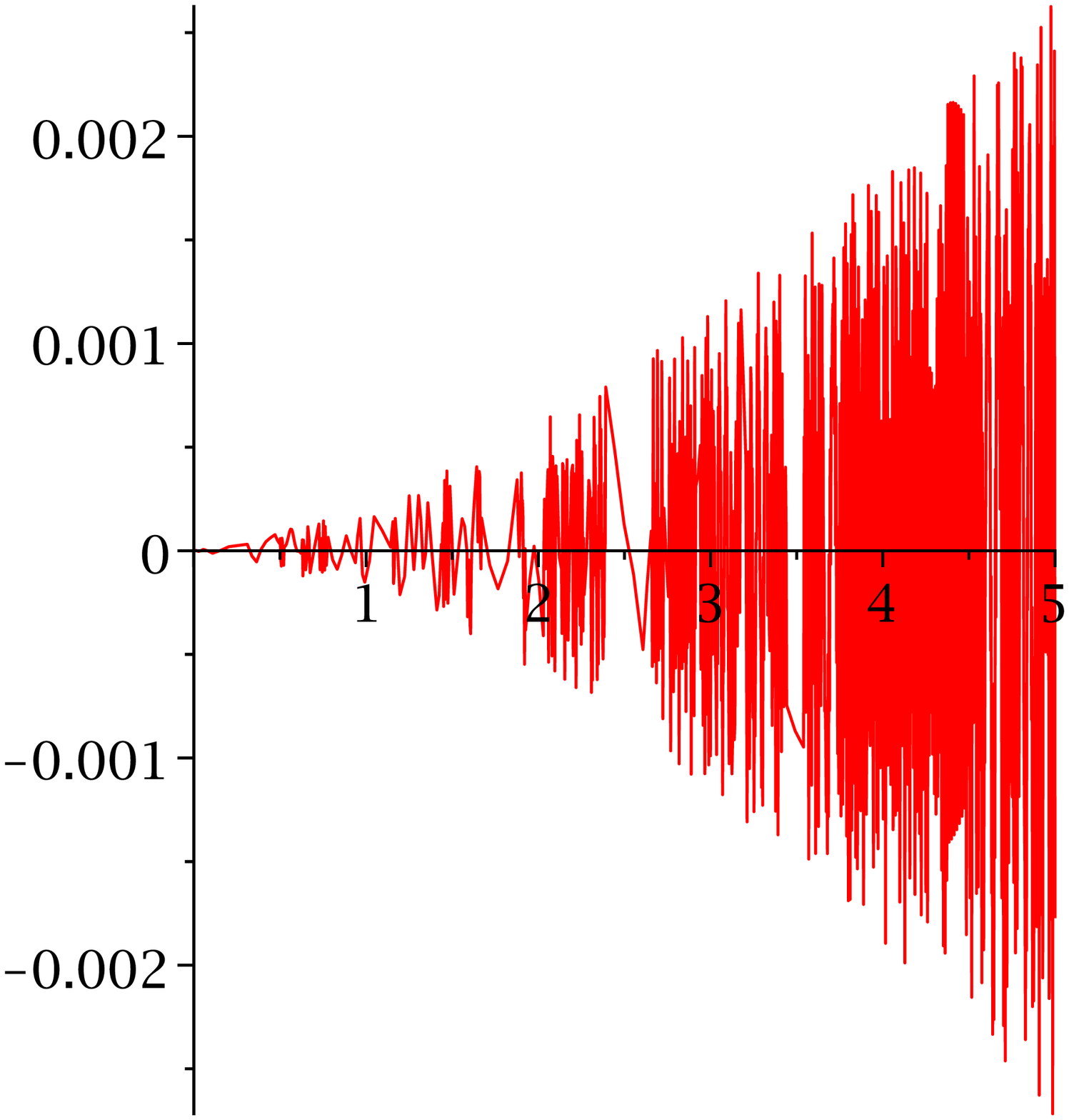}}\qquad\qquad\qquad
{\label{fig:re}\includegraphics[width=4cm,height=4cm]{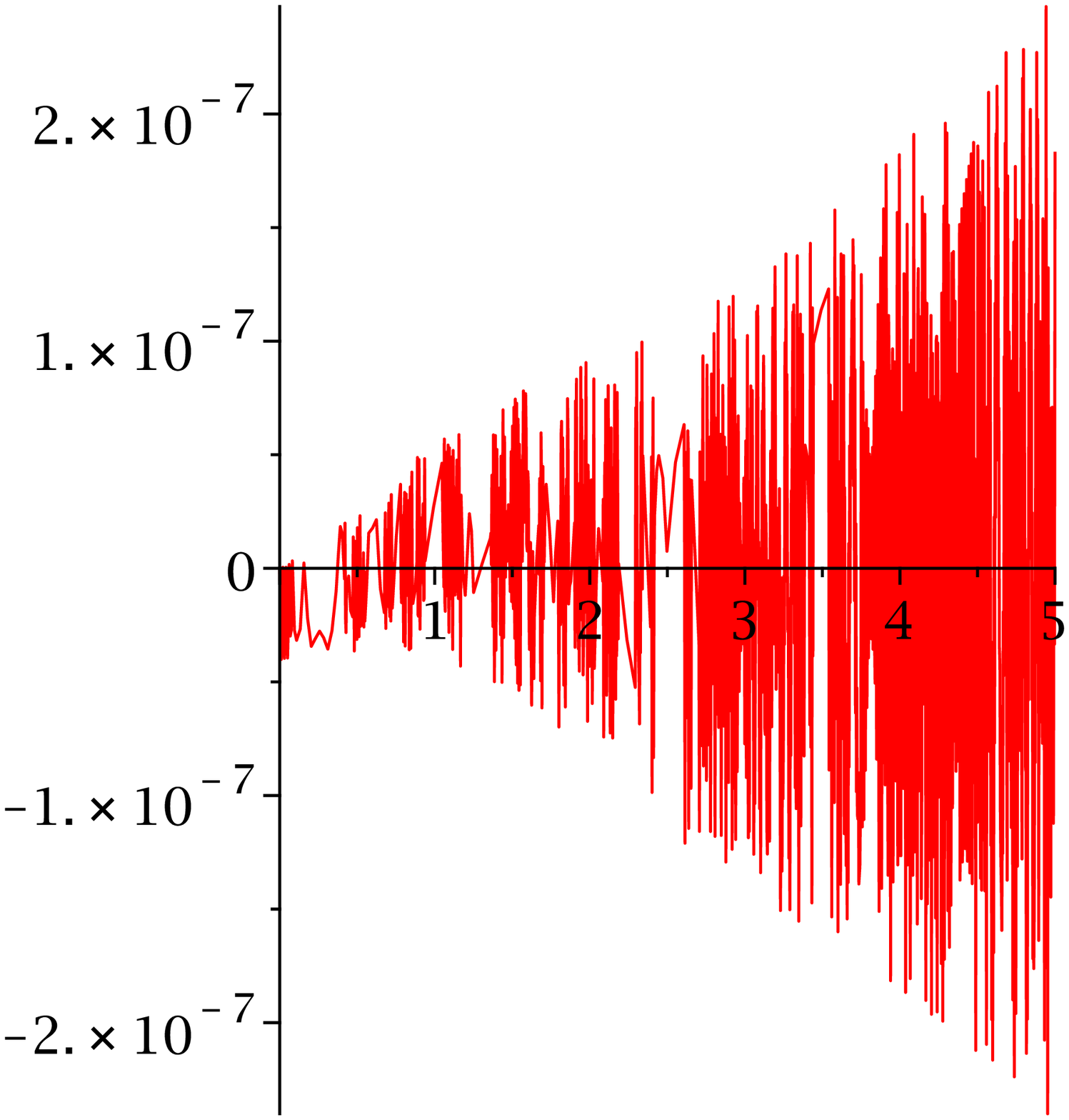}}\\
{\label{fig:re}\includegraphics[width=4cm,height=4cm]{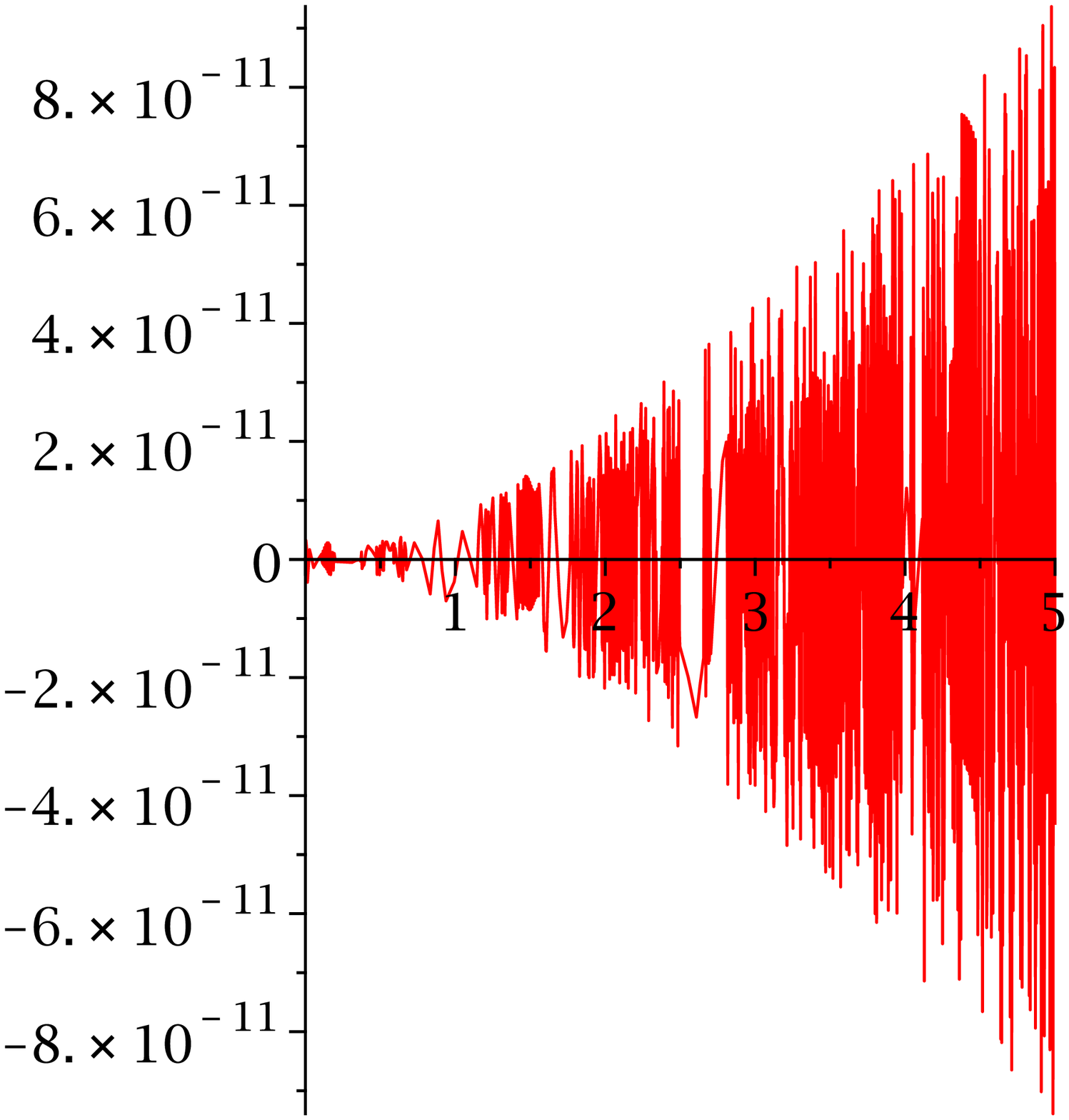}}\qquad\qquad\qquad
{\label{fig:re}\includegraphics[width=4cm,height=4cm]{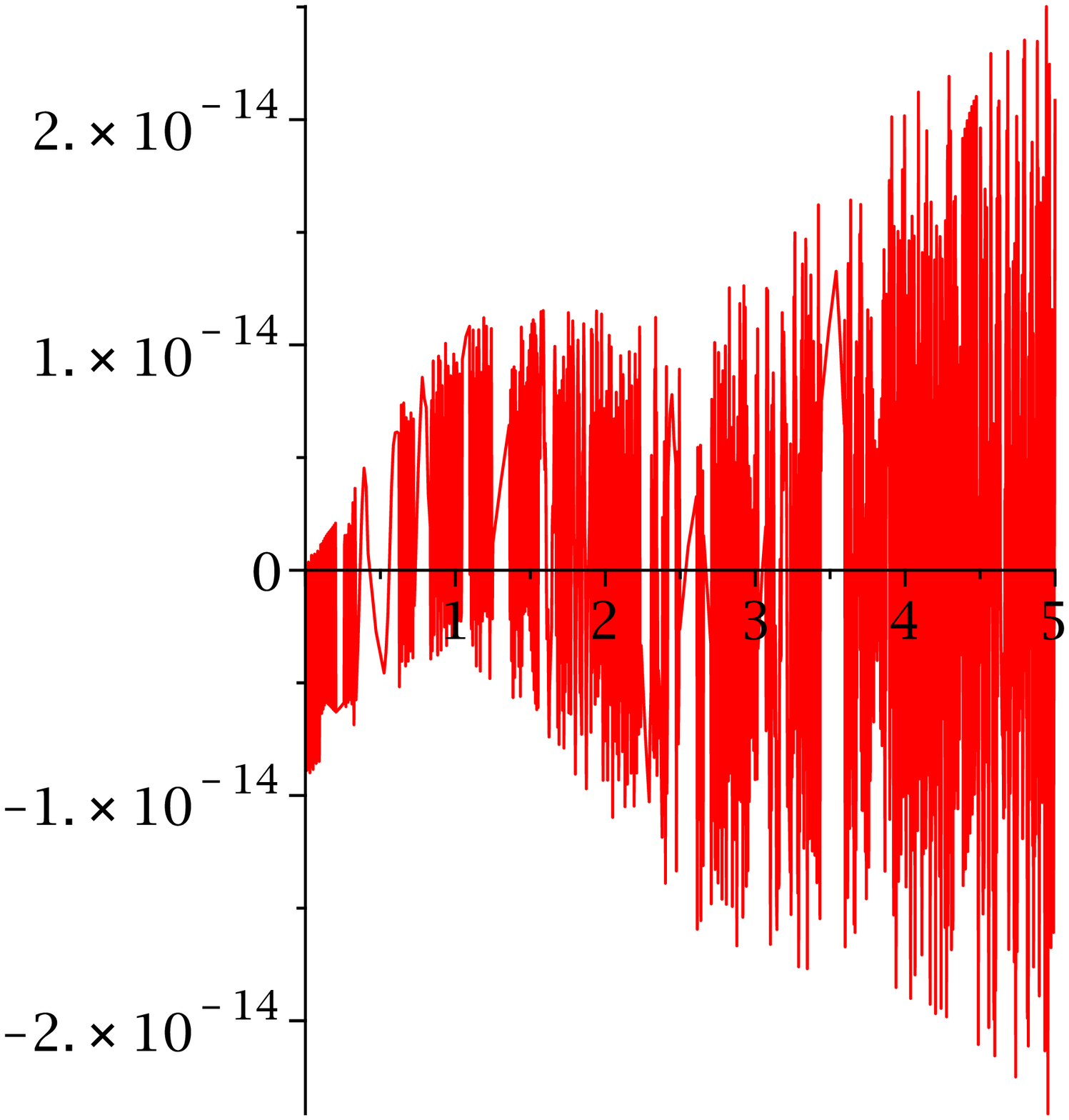}}\\
{\label{fig:re}\includegraphics[width=4cm,height=4cm]{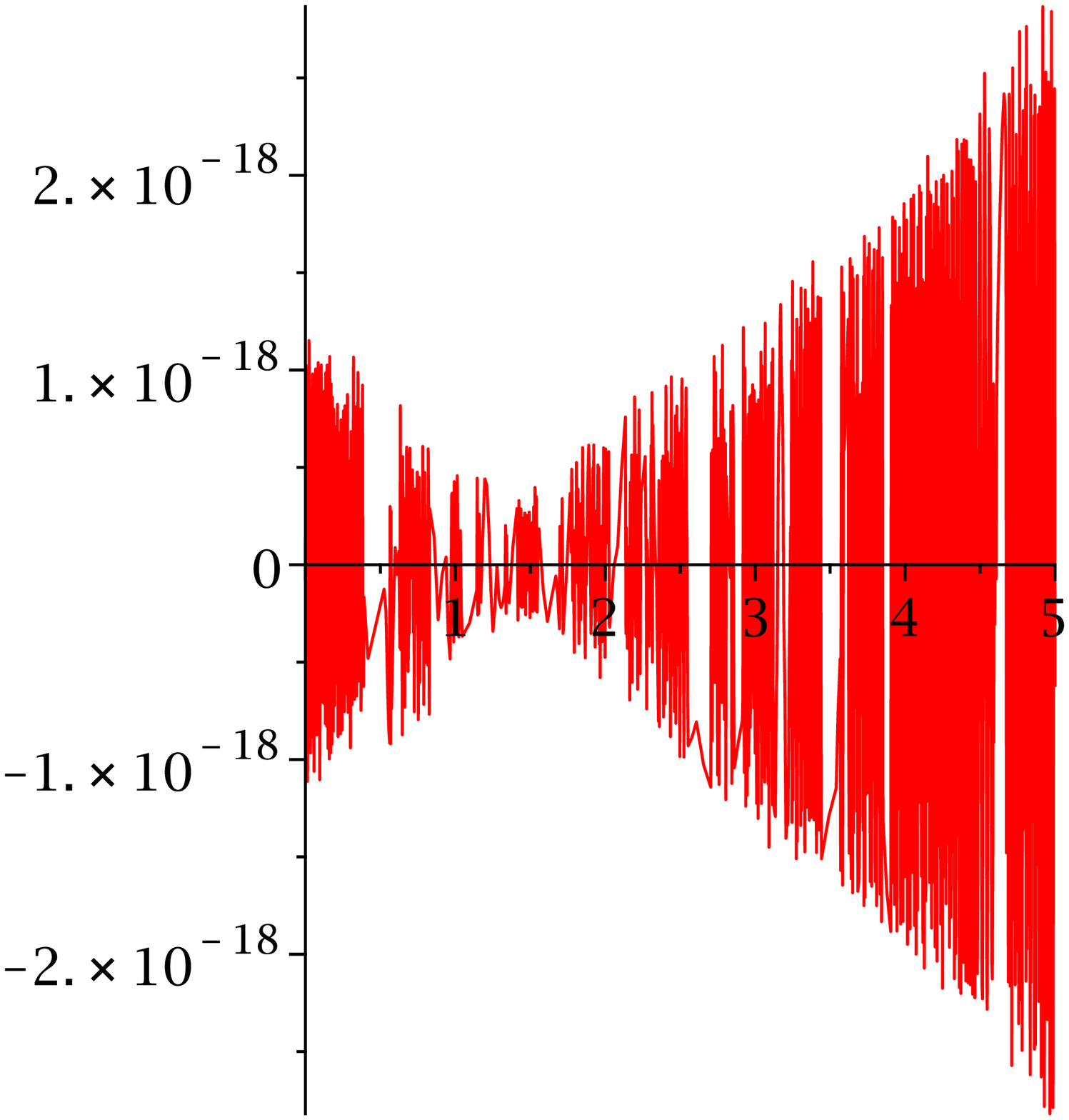}}
\caption{The real part of the error in $\dot{x}$ committed by the asymptotic expansion,
as applied to the linear system (4.1) with $\omega = 5000$ for $s = 0, 1, 2, 3, 4$
(from top left onwards).}
\end{figure}

\subsection{A nonlinear example in Memristor circuits}

In this subsection, the method in this paper is developed for
Memristor circuits subject to high-frequency signals. Consider the
following differential equation governing a circuit with two
memristors similar to that given in BoCheng2011,
\begin{eqnarray}
\left\{
\begin{array}{llllll}
y'_1(t) = y_3(t)\\
y'_2(t) = \frac{y_4(t) - y_3(t)}{1 + e(1 + 3y_2^2(t))}\\
y'_3(t) = a y_3(t) (d - (1 + 3 y_1^2(t))) - \frac{a(y_3(t) - y_4(t))(1 + 3y_2^2(t))}{1 + e(1 + 3y_2^2(t))}\\
y'_4(t) = \frac{(y_3(t) - y_4(t))(1 + 3 y_2^2(t))}{1 + e(1 + 3y_2^2(t))} + y_5(t)\\
y'_5(t) = -b y_4(t) - cy_5(t) + s(t), \quad t \geq 0
\end{array}\right.
\end{eqnarray}
with the initial conditions
$$
\left(y_1(0), y_2(0), y_3(0), y_4(0), y_5(0)\right)^T = (c_1, c_2,
0, 10^{-4}, 0)^T=y_0
$$
where $a = 8$, $b = 10$, $c = 0$, $d = 2$, $e = 0.1$, $c_1 = -0.8$,
$c_2 = - 0.4$. The unknown functions are $y_1(t), y_2(t), y_3(t),
y_4(t), y_5(t)$. The forcing term is $s(t) = \frac{A b}{2i} e^{i
\kappa_1 \omega t} - \frac{A b}{2i} e^{i \kappa_2 \omega t} +
\frac{A b}{2i} e^{i \kappa_3 \omega t} - \frac{A b}{2i} e^{i
\kappa_4 \omega t}$, in which $A = 0.1$, $\kappa_1 = 1$, $\kappa_1 =
-1$, $\kappa_3 = \sqrt{2}$, $\kappa_4 = -\sqrt{2}$ and $\omega$ is
our oscillatory parameter. The circuit figure is shown in Fig. 4.5
where the corresponding parameter relationship between Fig. 4.5 and
the equation (4.2) are
\begin{eqnarray*}
&& \phi_1 = y_1, \quad \phi_2 = y_2, \quad v_3 = y_3, \quad v_4 =
y_4, \quad i_5 = y_5, \quad C_2 = 1, \quad a = 1/C_1, \\
&& W_1 = 1 + 3 y_1^2,\quad W_2 = 1 + 3 y_2^2, \quad b = 1/L, \quad c
= r/L, \quad d = G, \quad e = R.
\end{eqnarray*}

\begin{figure}[htbp]
\centering
{\label{fig:4.5}\includegraphics[width=13cm,height=6cm]{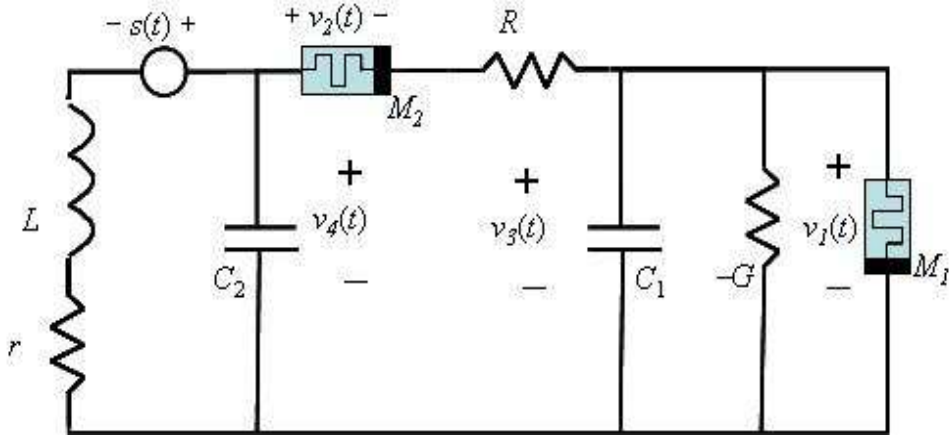}}\\
\caption{The Memristor circuits.}
\end{figure}

This circuit equation can be written in vector form
\begin{eqnarray*}
\boldsymbol{y}'(t) = \boldsymbol{f}(\boldsymbol{y}) +
\boldsymbol{a}_1 (t) e^{i \kappa_1 \omega t} + \boldsymbol{a}_2 (t)
e^{i \kappa_2 \omega t} + \boldsymbol{a}_3 (t) e^{i \kappa_3 \omega
t} + \boldsymbol{a}_4 (t) e^{i \kappa_4 \omega t}, \quad t\geq 0,
\end{eqnarray*}
where $\boldsymbol{y}(t) = \left(y_1(t), y_2(t), y_3(t), y_4(t),
y_5(t)\right)^T$,
\begin{eqnarray*}
\boldsymbol{f}(\boldsymbol{y}) = \left(
\begin{array}{llllll}
y_3\\
\frac{y_4 - y_3}{1 + e(1 + 3y_2^2)}\\
a y_3 (d - (1 + 3 y_1^2)) - \frac{a(y_3 - y_4)(1 + 3y_2^2)}{1 + e(1
+ 3y_2^2)}\\
\frac{(y_3 - y_4)(1 + 3 y_2^2)}{1 + e(1 + 3y_2^2)} + y_5\\
-b y_4 - cy_5
\end{array}\right)
\end{eqnarray*}
and
\begin{eqnarray*}
\boldsymbol{a}_1 (t) = \left(
\begin{array}{llllll}
0\\
0\\
0\\
0\\
\frac{A b}{2i}
\end{array}\right),
\boldsymbol{a}_2 (t) = \left(
\begin{array}{llllll}
0\\
0\\
0\\
0\\
-\frac{A b}{2i}
\end{array}\right),
\boldsymbol{a}_3 (t) = \left(
\begin{array}{llllll}
0\\
0\\
0\\
0\\
\frac{A b}{2i}
\end{array}\right),
\boldsymbol{a}_4 (t) = \left(
\begin{array}{llllll}
0\\
0\\
0\\
0\\
-\frac{A b}{2i}
\end{array}\right),
\end{eqnarray*}
or in a more compact form as \begin{eqnarray} \boldsymbol{y}'(t) +
\boldsymbol{f}(\boldsymbol{y}(t)) \ = \sum_{m = 1}^{4}
\boldsymbol{a}_m(t) e^{i \kappa_m \omega t} = \sum_{m \in
\mathcal{U}_0} \boldsymbol{a}_m(t) e^{i \kappa_m \omega t} , \quad t
\geq 0.\nonumber\\
\end{eqnarray}

where  $\mathcal{U}_0 = \{1, 2, \cdots, 4\}$ is an initial set and
$\omega_j = \kappa_j \omega$, $j = 1, 2, \cdots$.

Then the asymptotic method is developed for this
type of equation.

The formation of the terms in the asymptotic method for the given
memristor system shall now be described.

\subsubsection{The zeroth terms}

Denote $(\boldsymbol{p})_j$ as the $j$-th element of the vector
$\boldsymbol{p}$. When $r = 0$, set $\mathcal{U}_1 = \{0, 1, 2, 3,
4\}$. Then the zeroth term $\boldsymbol{p}_{0, 0}(t)$ obeys
$$
\boldsymbol{p}'_{0, 0} = \boldsymbol{f}(\boldsymbol{p}_{0, 0}),
\quad t \geq 0, \quad \boldsymbol{p}_{0, 0}(0) = \boldsymbol{y}_0 =
(c_1, c_2, 0, 10^{-4}, 0)^T,
$$
where
\begin{eqnarray*}
\boldsymbol{f}(\boldsymbol{p}_{0, 0}) = \left(
\begin{array}{llllll}
(\boldsymbol{p}_{0, 0})_3\\
\frac{(\boldsymbol{p}_{0, 0})_4 - (\boldsymbol{p}_{0, 0})_3}{1 + e(1 + 3(\boldsymbol{p}_{0, 0})_2^2)}\\
a (\boldsymbol{p}_{0, 0})_3 (d - (1 + 3 (\boldsymbol{p}_{0,
0})_1^2)) - \frac{a((\boldsymbol{p}_{0, 0})_3 - (\boldsymbol{p}_{0,
0})_4)(1 + 3(\boldsymbol{p}_{0, 0})_2^2)}{1 + e(1
+ 3(\boldsymbol{p}_{0, 0})_2^2)}\\
\frac{((\boldsymbol{p}_{0, 0})_3 - (\boldsymbol{p}_{0, 0})_4)
(1 + 3 (\boldsymbol{p}_{0, 0})_2^2)}{1 + e(1 + 3(\boldsymbol{p}_{0, 0})_2^2)} + (\boldsymbol{p}_{0, 0})_5\\
-b (\boldsymbol{p}_{0, 0})_4 - c(\boldsymbol{p}_{0, 0})_5
\end{array}\right)
\end{eqnarray*}

In addition, the recursions enable the determination of
$\boldsymbol{p}_{1, m}(t)$, $m \neq 0$,
\begin{eqnarray*}
&& \boldsymbol{p}_{1, 1} (t) = \frac{1}{i \kappa_1}\left(
\begin{array}{llllll}
0\\
0\\
0\\
0\\
\frac{A b}{2i}
\end{array}\right), \quad
\boldsymbol{p}_{1, 2} (t) = \frac{1}{i \kappa_2} \left(
\begin{array}{llllll}
0\\
0\\
0\\
0\\
-\frac{A b}{2i}
\end{array}\right),\nonumber\\
&& \boldsymbol{p}_{1, 3} (t) = \frac{1}{i \kappa_3}\left(
\begin{array}{llllll}
0\\
0\\
0\\
0\\
\frac{A b}{2i}
\end{array}\right), \quad
\boldsymbol{p}_{1, 4} (t) = \frac{1}{i \kappa_4}\left(
\begin{array}{llllll}
0\\
0\\
0\\
0\\
-\frac{A b}{2i}
\end{array}\right).
\end{eqnarray*}
The corresponding derivatives are $\boldsymbol{p}'_{1, 1}(t) =
\boldsymbol{p}'_{1, 2}(t) = \boldsymbol{p}'_{1, 3}(t) =
\boldsymbol{p}'_{1, 4}(t) = \boldsymbol{0}$.

\subsubsection{The $r=1$ terms}

For $r = 1$, set $\mathcal{U}_2 = \mathcal{U}_1 = \{0, 1, 2, 3,
4\}$. This yields
\begin{eqnarray*}
&& \boldsymbol{p}'_{1, 0} = \boldsymbol{f}_1(\boldsymbol{p}_{0,
0})[\boldsymbol{p}_{1, 0}], \quad t \geq 0, \\
&& \boldsymbol{p}_{1, 0}(0) = - \sum_{m \in \mathcal{U}_1\setminus
\{0\}}\boldsymbol{p}_{1, m}(0) = \left(
\begin{array}{llllll}
0 \\
0 \\
0\\
0 \\
\frac{Ab}{2}(2 + \sqrt{2})\end{array}\right),\\
&& \boldsymbol{p}_{2, m} = \frac{1}{i \kappa_m}
\{\boldsymbol{f}_1(\boldsymbol{p}_{0, 0})[\boldsymbol{p}_{1, m}] -
\boldsymbol{p}'_{1, m}\} = \frac{1}{i \kappa_m}
\boldsymbol{f}_1(\boldsymbol{p}_{0, 0})[\boldsymbol{p}_{1, m}],
\quad m \in \mathcal{U}_2 \setminus \{0\},\nonumber
\end{eqnarray*}
where
\begin{eqnarray*}
&& \boldsymbol{p}_{2, 1} (t) = \frac{1}{(i \kappa_1)^2}\left(
\begin{array}{llllll}
0\\
0\\
0\\
\frac{Ab}{2 i}\\
\frac{- c A b}{2i}
\end{array}\right), \quad
\boldsymbol{p}_{2, 2} (t) = \frac{1}{(i \kappa_2)^2} \left(
\begin{array}{llllll}
0\\
0\\
0\\
-\frac{Ab}{2 i}\\
\frac{c A b}{2i}
\end{array}\right),\nonumber\\
&& \boldsymbol{p}_{2, 3} (t) = \frac{1}{(i \kappa_3)^2}\left(
\begin{array}{llllll}
0\\
0\\
0\\
\frac{Ab}{2 i}\\
\frac{-c A b}{2i}
\end{array}\right), \quad
\boldsymbol{p}_{2, 4} (t) = \frac{1}{(i \kappa_4)^2}\left(
\begin{array}{llllll}
0\\
0\\
0\\
-\frac{Ab}{2 i}\\
\frac{c A b}{2i}
\end{array}\right).
\end{eqnarray*}

\subsubsection{when $r=2$}

The $r=2$ layer is the first layer in which additional frequencies
must be considered and we set
$$\mathcal{U}_3 = \{0, 1, 2, 3, 4, (1,
1), (1, 3), (1, 4), (2, 2), (2, 3), (2, 4), (3, 3), (4, 4)\}.
$$  Note
that the (1,2) term and (3,4) terms are not present as addition of
these frequencies would result in zero which is present in the set.

We will first consider the $\boldsymbol{p}_{2, 0}(t)$ term. Since
$\rho^0_{1, 1} = 1$, $\rho^0_{1, 2} = 2$ and $\rho^0_{3, 4} = 2$,
the term $\boldsymbol{p}_{2, 0}$ satisfies
\begin{eqnarray*}
&& \boldsymbol{p}'_{2, 0} = \boldsymbol{f}_1(\boldsymbol{p}_{0,
0})[\boldsymbol{p}_{2, 0}] + \frac{1}{2}
\mathop{\sum_{\kappa_{\ell_1} + \kappa_{\ell_2} = 0}}_{\ell_1 \leq
\ell_2} \rho^0_{\ell_1, \ell_2} \boldsymbol{f}_2(\boldsymbol{p}_{0,
0})[\boldsymbol{p}_{1, \ell_1}, \boldsymbol{p}_{1, \ell_2}]\\
&& = \boldsymbol{f}_1(\boldsymbol{p}_{0, 0})[\boldsymbol{p}_{2, 0}]
+ \frac{1}{2} \boldsymbol{f}_2(\boldsymbol{p}_{0,
0})[\boldsymbol{p}_{1, 0}, \boldsymbol{p}_{1, 0}]\\
&&\qquad  + \boldsymbol{f}_2(\boldsymbol{p}_{0,
0})[\boldsymbol{p}_{1, 1}, \boldsymbol{p}_{1, 2}] +
\boldsymbol{f}_2(\boldsymbol{p}_{0, 0})[\boldsymbol{p}_{1, 3},
\boldsymbol{p}_{1, 4}]
\end{eqnarray*}
with the initial condition
\begin{eqnarray*}
\boldsymbol{p}_{2, 0}(0) = - \sum_{m \in \mathcal{U}_2 \setminus
\{0\}} \boldsymbol{p}_{2, m}(0) = \boldsymbol{0},
\end{eqnarray*}
where
\begin{eqnarray*}
\boldsymbol{f}_2(\boldsymbol{p}_{0, 0})[\boldsymbol{p}_{1, m},
\boldsymbol{p}_{1, k}] = \left(
\begin{array}{llllll}
\boldsymbol{p}_{1, m}^T M_1(\boldsymbol{p}_{0, 0}) \boldsymbol{p}_{1, k}\\
\boldsymbol{p}_{1, m}^T M_2(\boldsymbol{p}_{0, 0}) \boldsymbol{p}_{1, k}\\
\boldsymbol{p}_{1, m}^T M_3(\boldsymbol{p}_{0, 0}) \boldsymbol{p}_{1, k}\\
\boldsymbol{p}_{1, m}^T M_4(\boldsymbol{p}_{0, 0}) \boldsymbol{p}_{1, k}\\
\boldsymbol{p}_{1, m}^T M_5(\boldsymbol{p}_{0, 0})
\boldsymbol{p}_{1, k}
\end{array}\right)
\end{eqnarray*}
and $M_j(\boldsymbol{y})$ is the $5 \times 5$ dimensional Jacobian
matrix evaluated at the vector function $\boldsymbol{p}_{0, 0}(t)$
\begin{eqnarray*}
M_j(\boldsymbol{p}_{0, 0}) = \left(
\begin{array}{llllll}
\frac{\partial^2 f_j}{\partial y_1 \partial y_1} & \frac{\partial^2
f_j}{\partial y_1 \partial y_2} & \frac{\partial^2 f_j}{\partial y_1
\partial y_3} & \frac{\partial^2 f_j}{\partial y_1 \partial y_4} &
\frac{\partial^2 f_j}{\partial y_1 \partial y_5}\\
\frac{\partial^2 f_j}{\partial y_2 \partial y_1} & \frac{\partial^2
f_j}{\partial y_2 \partial y_2} & \frac{\partial^2 f_j}{\partial y_2
\partial y_3} & \frac{\partial^2 f_j}{\partial y_2 \partial y_4} &
\frac{\partial^2 f_j}{\partial y_2 \partial y_5}\\
\frac{\partial^2 f_j}{\partial y_3 \partial y_1} & \frac{\partial^2
f_j}{\partial y_3 \partial y_2} & \frac{\partial^2 f_j}{\partial y_3
\partial y_3} & \frac{\partial^2 f_j}{\partial y_3 \partial y_4} &
\frac{\partial^2 f_j}{\partial y_3 \partial y_5}\\
\frac{\partial^2 f_j}{\partial y_4 \partial y_1} & \frac{\partial^2
f_j}{\partial y_4 \partial y_2} & \frac{\partial^2 f_j}{\partial y_4
\partial y_3} & \frac{\partial^2 f_j}{\partial y_4 \partial y_4} &
\frac{\partial^2 f_j}{\partial y_4 \partial y_5}\\
\frac{\partial^2 f_j}{\partial y_5 \partial y_1} & \frac{\partial^2
f_j}{\partial y_5 \partial y_2} & \frac{\partial^2 f_j}{\partial y_5
\partial y_3} & \frac{\partial^2 f_j}{\partial y_5 \partial y_4} &
\frac{\partial^2 f_j}{\partial y_5 \partial y_5}
\end{array}\right)_{\boldsymbol{y} = \boldsymbol{p}_{0, 0}(t)}.
\end{eqnarray*}

Furthermore, due to the fact that the first four elements of
$\boldsymbol{p}_{1, m}(t)$, $m \neq 0$, are zero, the nonoscillatory
equation for $\boldsymbol{p}_{2, 0}(t)$ simplifies to
$$
\boldsymbol{p}'_{2, 0} = \boldsymbol{f}_1(\boldsymbol{p}_{0,
0})[\boldsymbol{p}_{2, 0}] + \frac{1}{2}
\boldsymbol{f}_2(\boldsymbol{p}_{0, 0})[\boldsymbol{p}_{1, 0},
\boldsymbol{p}_{1, 0}], \quad \boldsymbol{p}_{2, 0}(0) =
\boldsymbol{0}, \quad t \geq 0.
$$

Now consider the set $\mathcal{U}_3$ and the recursions for
$\boldsymbol{p}_{3, m}$. First, we match all of the terms in
$\mathcal{U}_2 \setminus \{0\}$ as these are also in
$\mathcal{U}_3$,
\begin{eqnarray*}
&& i \kappa_m \boldsymbol{p}_{3, m} =
\boldsymbol{f}_1(\boldsymbol{p}_{0, 0}) [\boldsymbol{p}_{2, m}] -
\boldsymbol{p}'_{2, m} + \frac{1}{2} \mathop{\sum_{\kappa_{\ell_1} +
\kappa_{\ell_2}= \kappa_m}}_{\ell_1 \leq \ell_2} \rho^m_{\ell_1,
\ell_2} \boldsymbol{f}_2(\boldsymbol{p}_{0, 0})
[\boldsymbol{p}_{1, \ell_1}, \boldsymbol{p}_{1, \ell_2}] \nonumber\\
&& = \boldsymbol{f}_1(\boldsymbol{p}_{0, 0}) [\boldsymbol{p}_{2,
m}].
\end{eqnarray*}
We then match to the remainder of the elements in $\mathcal{U}_3$.
\begin{eqnarray*}
&& i (\kappa_{m_1} + \kappa_{m_2}) \boldsymbol{p}_{3, (m_1, m_2)} =
\frac{1}{2} \mathop{\sum_{\kappa_{\ell_1} + \kappa_{\ell_2}=
\kappa_{m_1} + \kappa_{m_2}}}_{\ell_1 \leq \ell_2} \rho^{m_1,
m_2}_{\ell_1, \ell_2} \boldsymbol{f}_2(\boldsymbol{p}_{0, 0})
[\boldsymbol{p}_{1, \ell_1}, \boldsymbol{p}_{1, \ell_2}] =
0.\nonumber\\
\end{eqnarray*}

Because of the nature of the elements of  $\boldsymbol{p}_{1, m}$
and  $\boldsymbol{p}_{2, m}$, the non-zero terms of
$\boldsymbol{p}_{3, m}(t)$, $m \neq 0$, are
\begin{eqnarray*}
&& \boldsymbol{p}_{3, 1} (t) = \frac{1}{(i \kappa_1)^3}\left(
\begin{array}{llllll}
0\\
\frac{1}{1 + e\left(1 + 3 (\boldsymbol{p}_{0, 0})^2_2\right)} \frac{A b}{2 i}\\
\frac{a\left(1 + 3 (\boldsymbol{p}_{0, 0})^2_2\right)}
{1 + e\left(1 + 3 (\boldsymbol{p}_{0, 0})^2_2\right)} \frac{A b}{2 i}\\
\frac{ -\left(1 + 3 (\boldsymbol{p}_{0, 0})^2_2\right)}
{1 + e\left(1 + 3 (\boldsymbol{p}_{0, 0})^2_2\right)} \frac{A b}{2 i} + \frac{-c A b}{2 i}\\
-\frac{A b^2}{2 i} +\frac{c^2 A b}{2 i}
\end{array}\right), \\
&& \boldsymbol{p}_{3, 2} (t) = \frac{1}{(i \kappa_2)^3}\left(
\begin{array}{llllll}
0\\
\frac{1}{1 + e\left(1 + 3 (\boldsymbol{p}_{0, 0})^2_2\right)} \frac{-A b}{2 i}\\
\frac{a\left(1 + 3 (\boldsymbol{p}_{0, 0})^2_2\right)}
{1 + e\left(1 + 3 (\boldsymbol{p}_{0, 0})^2_2\right)} \frac{-A b}{2 i}\\
\frac{ -\left(1 + 3 (\boldsymbol{p}_{0, 0})^2_2\right)}
{1 + e\left(1 + 3 (\boldsymbol{p}_{0, 0})^2_2\right)} \frac{-A b}{2 i} + \frac{c A b}{2 i}\\
\frac{A b^2}{2 i} - \frac{c^2 A b}{2 i}
\end{array}\right), \\
&& \boldsymbol{p}_{3, 3} (t) = \frac{1}{(i \kappa_3)^3}\left(
\begin{array}{llllll}
0\\
\frac{1}{1 + e\left(1 + 3 (\boldsymbol{p}_{0, 0})^2_2\right)} \frac{A b}{2 i}\\
\frac{a\left(1 + 3 (\boldsymbol{p}_{0, 0})^2_2\right)}
{1 + e\left(1 + 3 (\boldsymbol{p}_{0, 0})^2_2\right)} \frac{A b}{2 i}\\
\frac{ -\left(1 + 3 (\boldsymbol{p}_{0, 0})^2_2\right)}
{1 + e\left(1 + 3 (\boldsymbol{p}_{0, 0})^2_2\right)} \frac{A b}{2 i} + \frac{-c A b}{2 i}\\
-\frac{A b^2}{2 i} + \frac{c^2 A b}{2 i}
\end{array}\right), \\
&& \boldsymbol{p}_{3, 4} (t) = \frac{1}{(i \kappa_4)^3}\left(
\begin{array}{llllll}
0\\
\frac{1}{1 + e\left(1 + 3 (\boldsymbol{p}_{0, 0})^2_2\right)} \frac{-A b}{2 i}\\
\frac{a\left(1 + 3 (\boldsymbol{p}_{0, 0})^2_2\right)}
{1 + e\left(1 + 3 (\boldsymbol{p}_{0, 0})^2_2\right)} \frac{-A b}{2 i}\\
\frac{ -\left(1 + 3 (\boldsymbol{p}_{0, 0})^2_2\right)}
{1 + e\left(1 + 3 (\boldsymbol{p}_{0, 0})^2_2\right)} \frac{-A b}{2 i} + \frac{c A b}{2 i}\\
\frac{A b^2}{2 i} - \frac{c^2 A b}{2 i}.
\end{array}\right).
\end{eqnarray*}

\subsubsection{The $r=3$ terms}

When $r = 3$, we note that
$\rho^0_{0, 0} = 1$, $\rho^0_{1, 2} = 2$, $\rho^0_{3, 4} = 2$,
$\rho^0_{0, 0, 0} = 1$, $\rho^0_{0, 1, 2} = 6$, $\rho^0_{0, 3, 4} =
6$.  Hence, the equation for $\boldsymbol{p}_{3, 0}$ is
\begin{eqnarray*}
&& \boldsymbol{p}'_{3, 0} = \boldsymbol{f}_1(\boldsymbol{p}_{0,
0})[\boldsymbol{p}_{3, 0}] + \mathop{\sum_{\kappa_{\ell_1} +
\kappa_{\ell_2} = 0}}_{\ell_1 \leq \ell_2} \rho^0_{\ell_1, \ell_2}
\boldsymbol{f}_2(\boldsymbol{p}_{0,
0})[\boldsymbol{p}_{1, \ell_1}, \boldsymbol{p}_{2, \ell_2}]\nonumber\\
&& + \frac{1}{6} \mathop{\sum_{\kappa_{\ell_1} + \kappa_{\ell_2} +
\kappa_{\ell_3} = 0}}_{\ell_1 \leq \ell_2 \leq \ell_3}
\rho^0_{\ell_1, \ell_2, \ell_3} \boldsymbol{f}_3(\boldsymbol{p}_{0,
0})[\boldsymbol{p}_{1, \ell_1}, \boldsymbol{p}_{1, \ell_2}, \boldsymbol{p}_{1, \ell_3}]\nonumber\\
&& = \boldsymbol{f}_1(\boldsymbol{p}_{0, 0})[\boldsymbol{p}_{3, 0}]
+ \boldsymbol{f}_2(\boldsymbol{p}_{0, 0})[\boldsymbol{p}_{1, 0},
\boldsymbol{p}_{2, 0}] + 2\boldsymbol{f}_2(\boldsymbol{p}_{0,
0})[\boldsymbol{p}_{1, 1},
\boldsymbol{p}_{2, 2}]\nonumber\\
&&  + 2 \boldsymbol{f}_2(\boldsymbol{p}_{0, 0})[\boldsymbol{p}_{1,
3}, \boldsymbol{p}_{2, 4}] + \frac{1}{6}
\boldsymbol{f}_3(\boldsymbol{p}_{0, 0})[\boldsymbol{p}_{1, 0},
\boldsymbol{p}_{1, 0}, \boldsymbol{p}_{1,
0}]\nonumber\\
&&+ \boldsymbol{f}_3(\boldsymbol{p}_{0, 0})[\boldsymbol{p}_{1, 0},
\boldsymbol{p}_{1, 1}, \boldsymbol{p}_{1, 2}] +
\boldsymbol{f}_3(\boldsymbol{p}_{0, 0})[\boldsymbol{p}_{1, 0},
\boldsymbol{p}_{1, 3}, \boldsymbol{p}_{1, 4}]\nonumber\\
&& = \boldsymbol{f}_1(\boldsymbol{p}_{0, 0})[\boldsymbol{p}_{3, 0}]
+ \boldsymbol{f}_2(\boldsymbol{p}_{0, 0})[\boldsymbol{p}_{1, 0},
\boldsymbol{p}_{2, 0}] + \frac{1}{6}
\boldsymbol{f}_3(\boldsymbol{p}_{0, 0})[\boldsymbol{p}_{1, 0},
\boldsymbol{p}_{1, 0}, \boldsymbol{p}_{1, 0}]\nonumber\\
&&+ \boldsymbol{f}_3(\boldsymbol{p}_{0, 0})[\boldsymbol{p}_{1, 0},
\boldsymbol{p}_{1, 1}, \boldsymbol{p}_{1, 2}] +
\boldsymbol{f}_3(\boldsymbol{p}_{0, 0})[\boldsymbol{p}_{1, 0},
\boldsymbol{p}_{1, 3}, \boldsymbol{p}_{1, 4}].\nonumber\\
\end{eqnarray*}
with
\begin{eqnarray*}
\boldsymbol{p}_{3, 0}(0) = - \sum_{m \in \mathcal{U}_3 \setminus
\{0\}} \boldsymbol{p}_{3, m}(0) = \left(
\begin{array}{llllll}
0\\
-\frac{A b}{1 + e (1 + 3 c^2_2)} \left(1 + \frac{1}{2 \sqrt{2}}\right)\\
-\frac{A b a (1 + 3 c^2_2)}{1 + e (1 + 3 c^2_2)} \left(1 + \frac{1}{2 \sqrt{2}}\right)\\
\left[c A b + \frac{A b(1 + 3 c^2_2)}{1 + e (1 + 3 c^2_2)}\right] \left(1 + \frac{1}{2 \sqrt{2}}\right)\\
\left(A b^2 - c^2 A b\right)\left(1 + \frac{1}{2 \sqrt{2}}\right)\\
\end{array}\right),
\end{eqnarray*}
where
\begin{eqnarray*}
&& (\boldsymbol{f}_3(\boldsymbol{p}_{0, 0}))_j[\boldsymbol{p}_{1,
m_1}, \boldsymbol{p}_{1, m_2}, \boldsymbol{p}_{1, m_3}]\nonumber\\
&& = \sum_{k_1 = 1}^{5} \sum_{k_2 = 1}^{5} \sum_{k_3 = 1}^{5}
\frac{\partial^3 \boldsymbol{f}_j}{\partial y_{k_1} \partial y_{k_2}
\partial y_{k_3} }|_{\boldsymbol{p}_{0,
0}}(\boldsymbol{p}_{1, m_1})_{k_1}(\boldsymbol{p}_{1,
m_2})_{k_2}(\boldsymbol{p}_{1,
m_3})_{k_3}.\nonumber\\
\end{eqnarray*}

Therefore, the asymptotic expansion including terms up to $r = 3$ is
\begin{eqnarray*}
&& \boldsymbol{y}(t) \sim \boldsymbol{p}_{0, 0}(t) +
\frac{1}{\omega} \left[\boldsymbol{p}_{1, 0}(t) + \boldsymbol{p}_{1,
1}(t) e^{i \kappa_1 \omega t} + \boldsymbol{p}_{1, 2}(t) e^{i
\kappa_2 \omega t}\right.\nonumber\\
&& \left. \qquad \qquad\qquad\qquad + \boldsymbol{p}_{1, 3}(t) e^{i
\kappa_3 \omega t} +
\boldsymbol{p}_{1, 4}(t) e^{i \kappa_4 \omega t}\right]\nonumber\\
&& + \frac{1}{\omega^2} \left[\boldsymbol{p}_{2, 0}(t) +
\boldsymbol{p}_{2, 1}(t) e^{i \kappa_1 \omega t} +
\boldsymbol{p}_{2, 2}(t) e^{i \kappa_2 \omega t} +
\boldsymbol{p}_{2, 3}(t) e^{i \kappa_3 \omega t} +
\boldsymbol{p}_{2, 4}(t) e^{i \kappa_4 \omega t}\right]\nonumber\\
&& + \frac{1}{\omega^3} \left[\boldsymbol{p}_{3, 0}(t) +
\boldsymbol{p}_{3, 1}(t) e^{i \kappa_1 \omega t} +
\boldsymbol{p}_{3, 2}(t) e^{i \kappa_2 \omega t} +
\boldsymbol{p}_{3, 3}(t) e^{i \kappa_3 \omega t} +
\boldsymbol{p}_{3, 4}(t) e^{i \kappa_4 \omega t}\right].
\end{eqnarray*}

\subsubsection{Numerical experiments}

The nonlinear Memristor circuits do not have a known analytical
solution, we come to a reference solution, the Maple routine
\textit{rkf45} with the accuracy tolerance $AbsErr=10^{-10}$ and
$RelErr=10^{-10}$. The terms $\boldsymbol{p}_{1,0}$,
$\boldsymbol{p}_{2,0}$ and $\boldsymbol{p}_{3,0}$ satisfy the
non-oscillatory ODEs which is solved by the Maple routine
\textit{rkf45} with $AbsErr=10^{-10}$ and $RelErr=10^{-10}$.

Figures 4.6 to 4.10 show the error functions for $y_1$, $y_2$,
$y_3$, $y_4$ and $y_5$ for the truncated parameter $s = 0, 1, 2, 3$
within $t \in [0, 3]$ when the oscillatory parameter is $\omega =
100$ and $\omega = 1000$. The error is seen to greatly reduce with
an increasing number of $r$ levels. Furthermore, with increasing the
oscillatory parameter, the error of the asymptotic method decreases
rapidly, a very important virtue of the method.

\begin{figure}[htbp]
\centering
{\label{fig:4.6}\includegraphics[width=6cm,height=5cm]{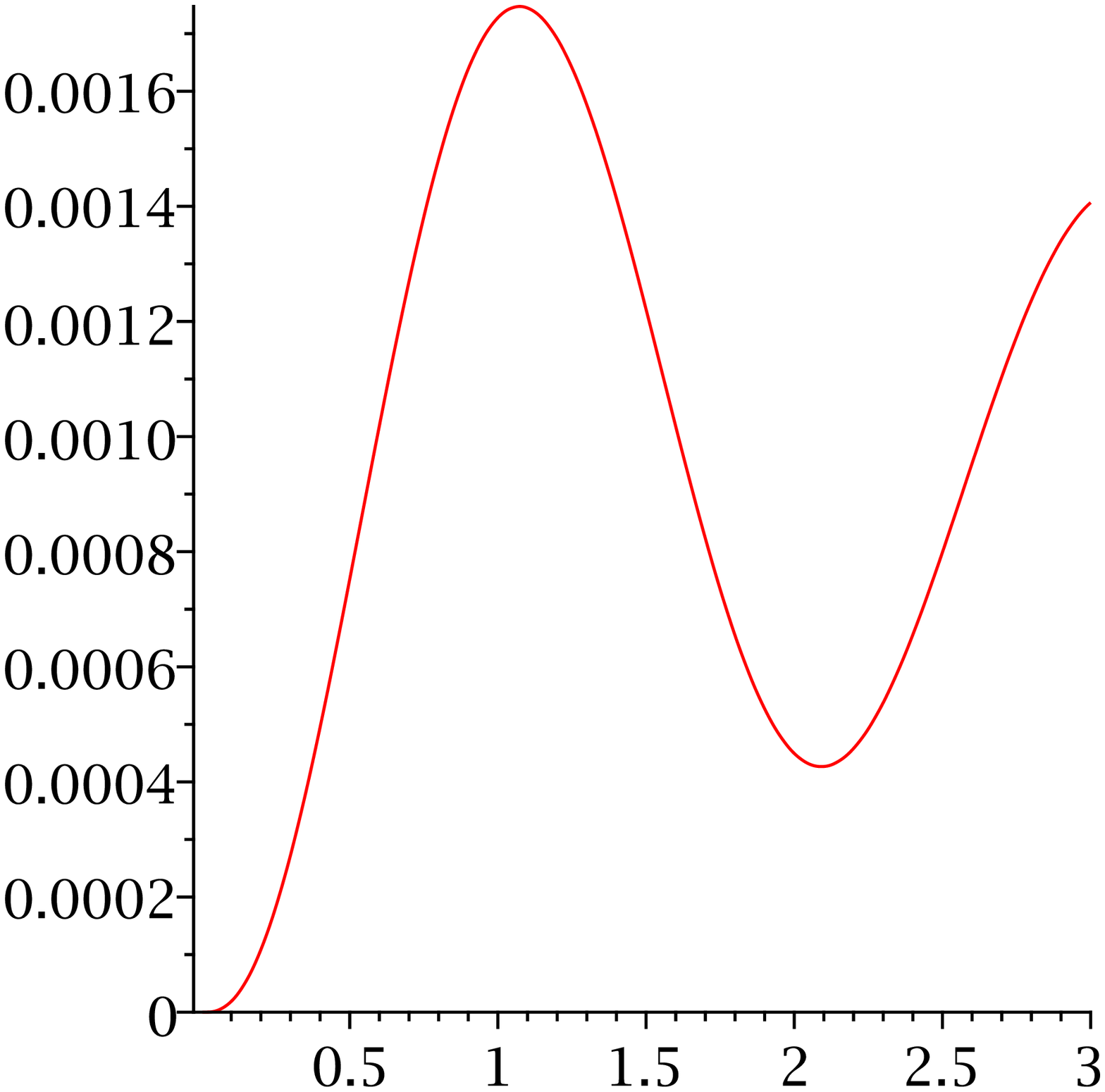}}
\quad
{\label{fig:4.6}\includegraphics[width=6cm,height=5cm]{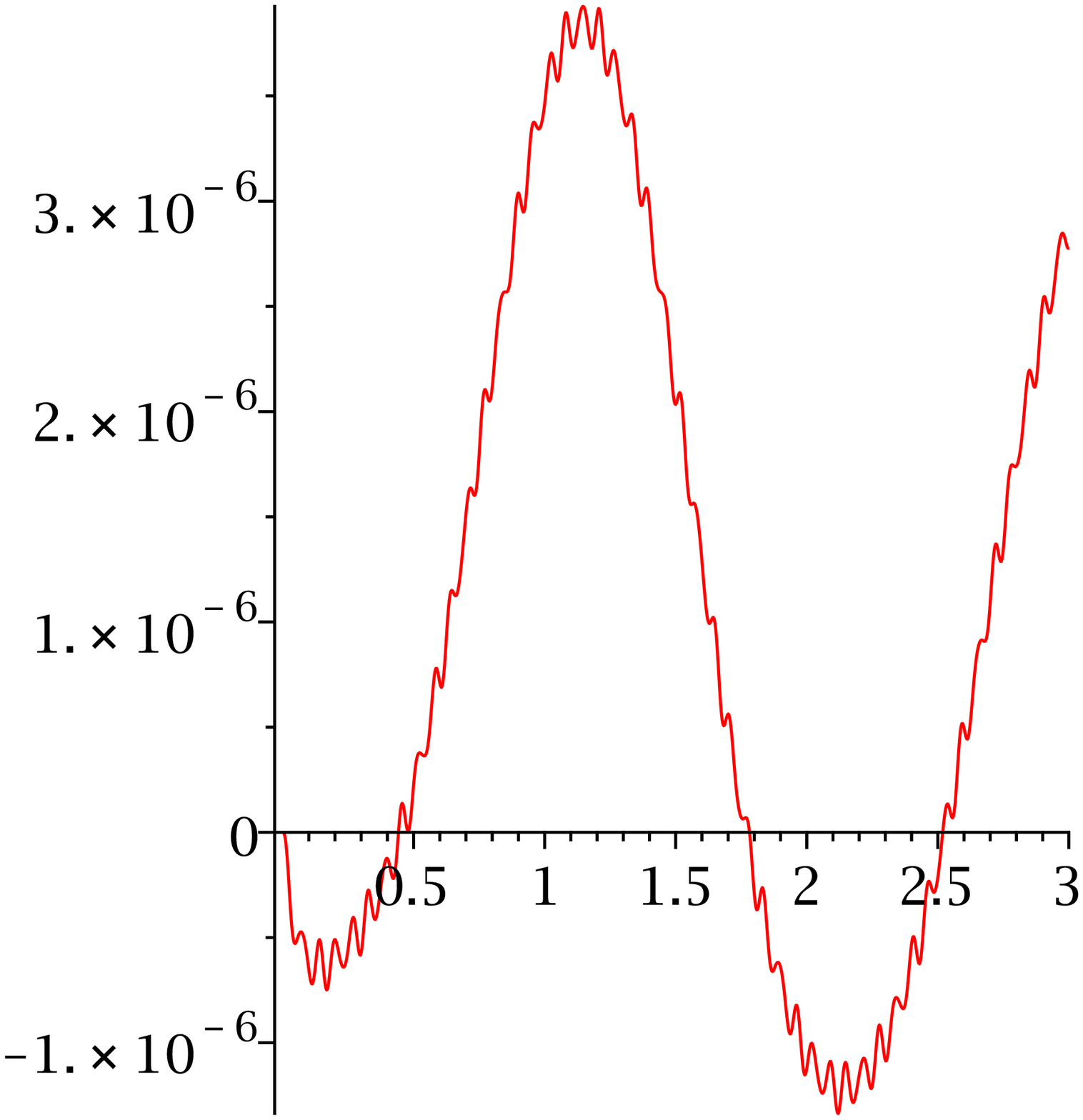}}
\\
{\label{fig:4.6}\includegraphics[width=6cm,height=5cm]{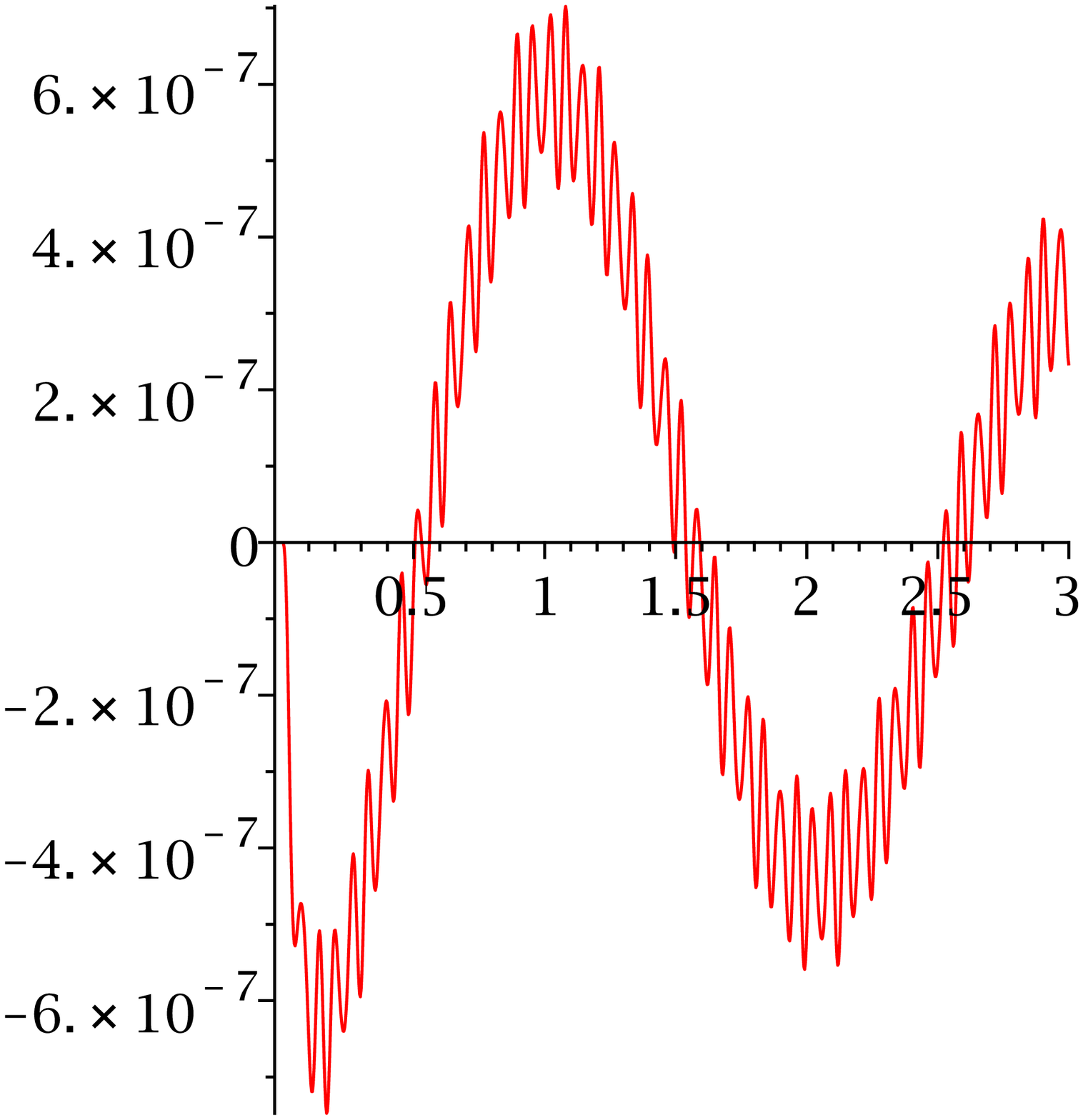}}
\quad
{\label{fig:4.6}\includegraphics[width=6cm,height=5cm]{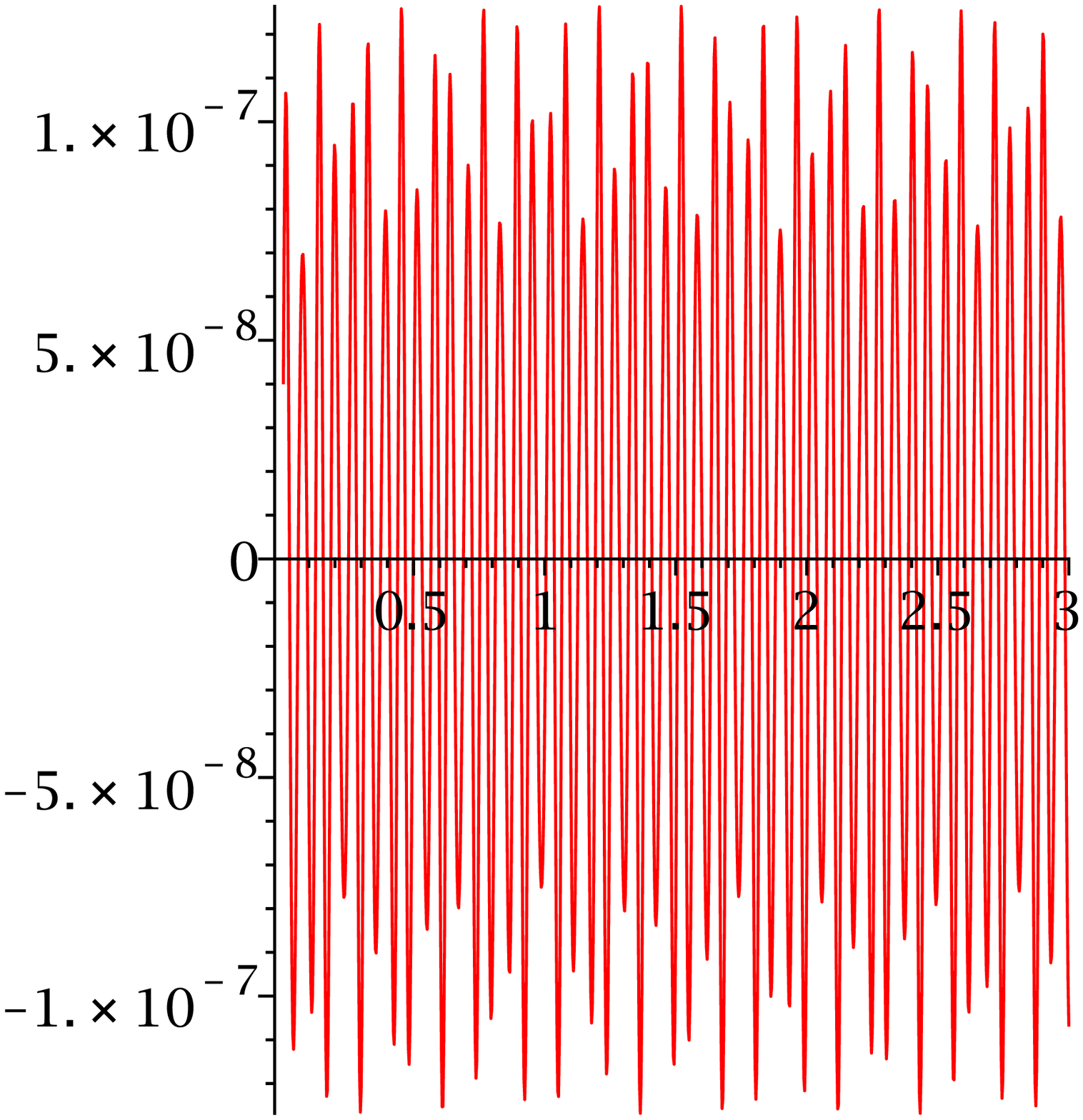}}
\\
{\label{fig:4.6}\includegraphics[width=6cm,height=5cm]{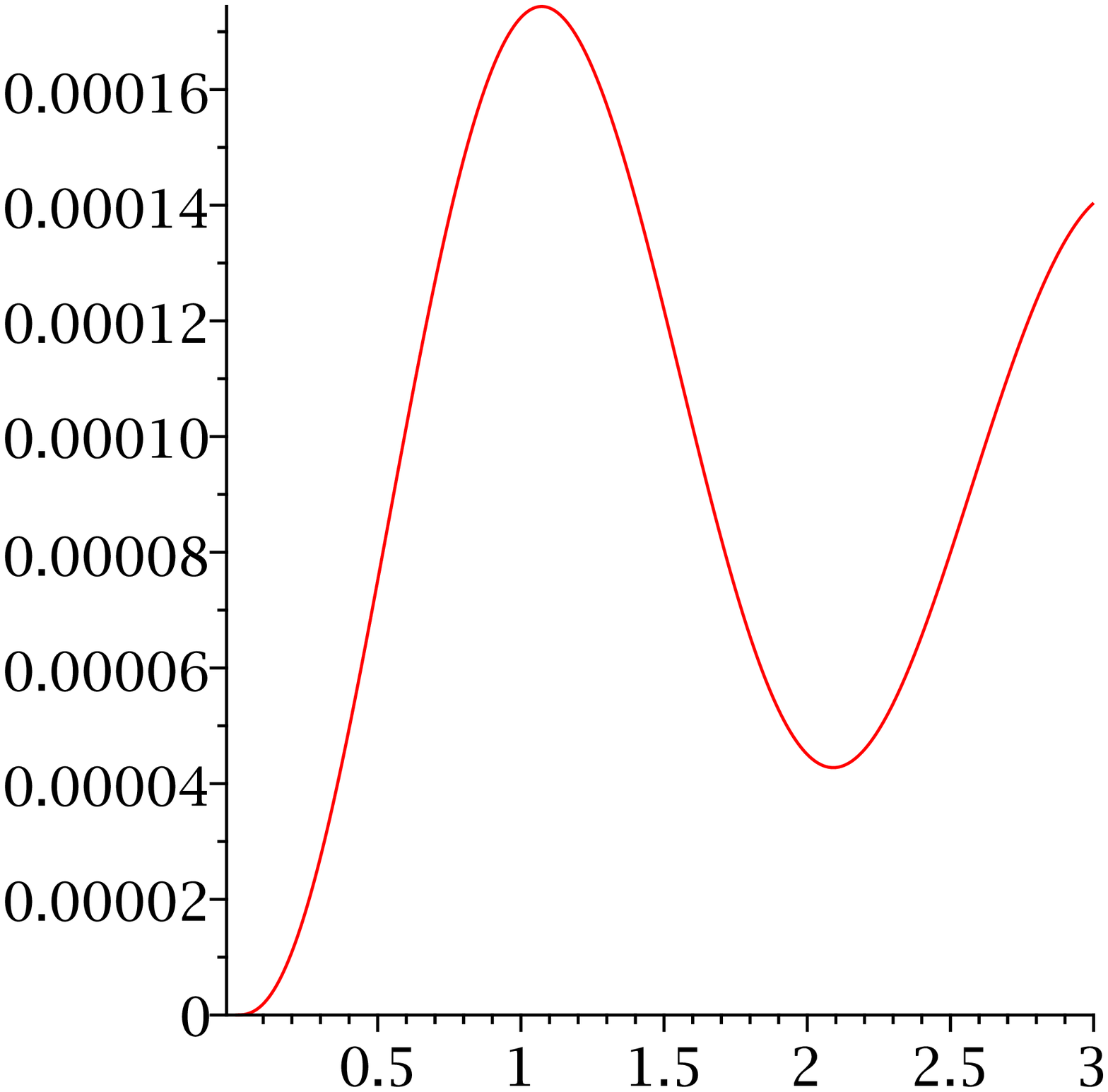}}
\quad
{\label{fig:4.6}\includegraphics[width=6cm,height=5cm]{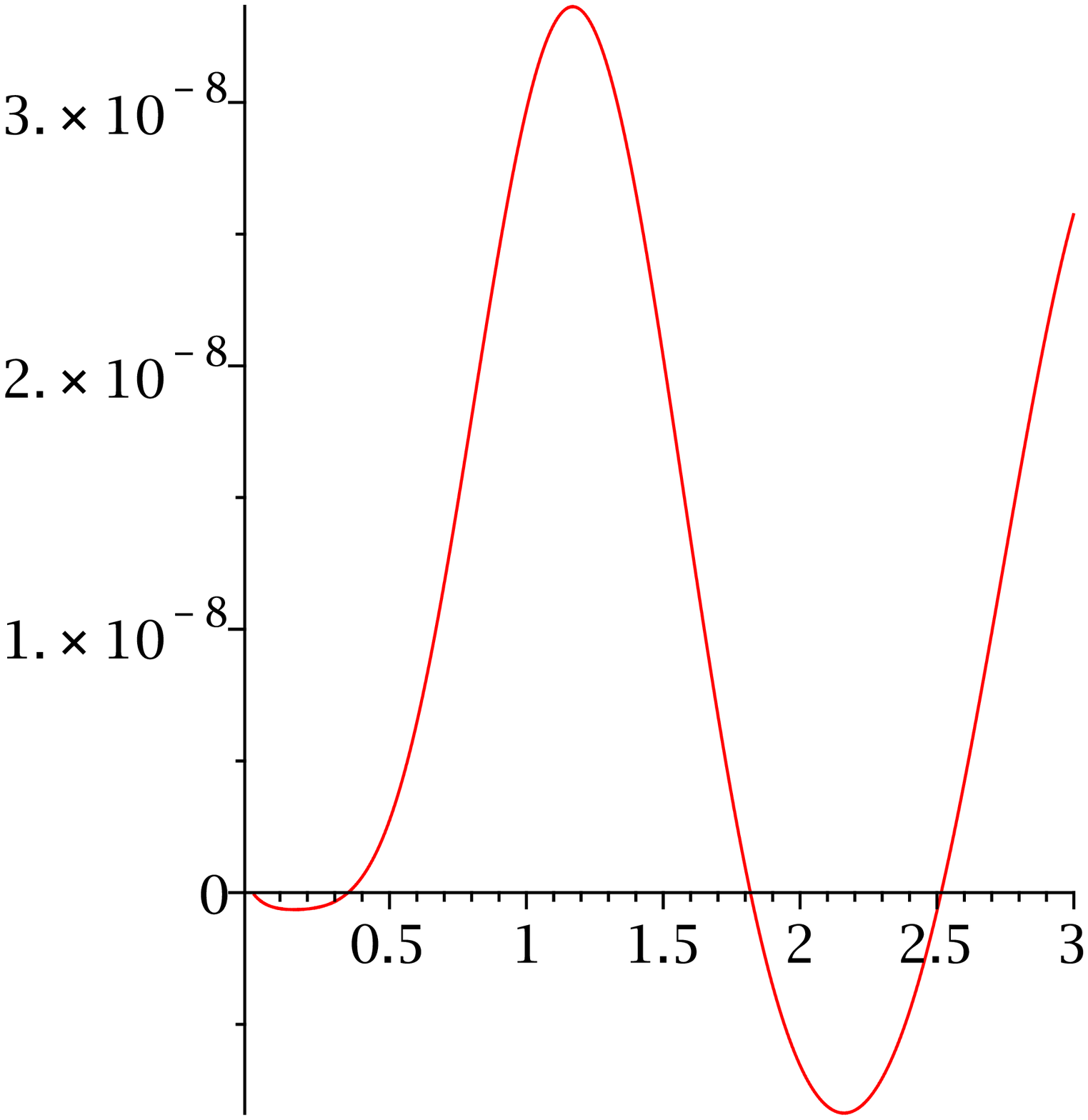}}\\
{\label{fig:4.6}\includegraphics[width=6cm,height=5cm]{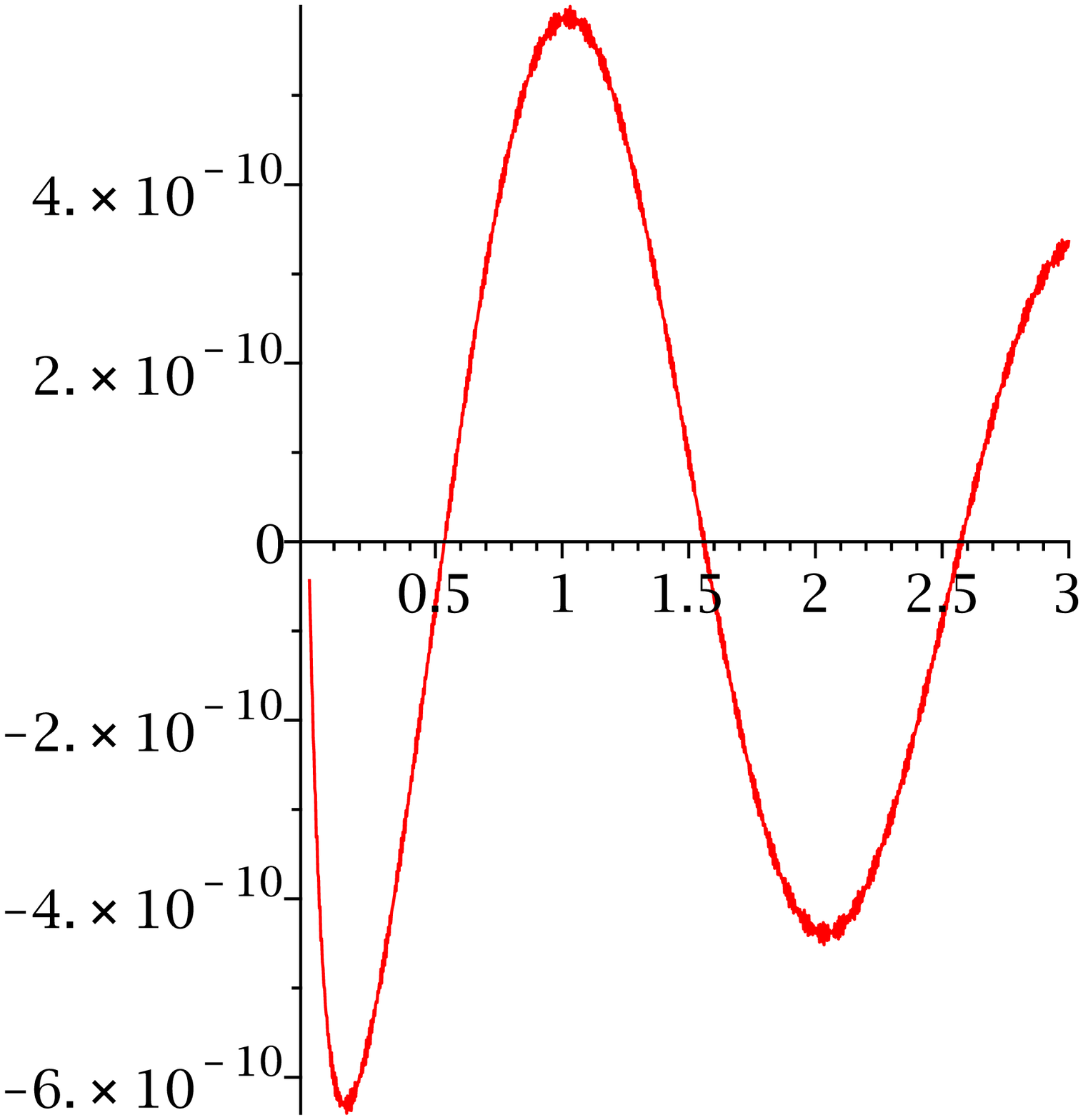}}
\quad
{\label{fig:4.6}\includegraphics[width=6cm,height=5cm]{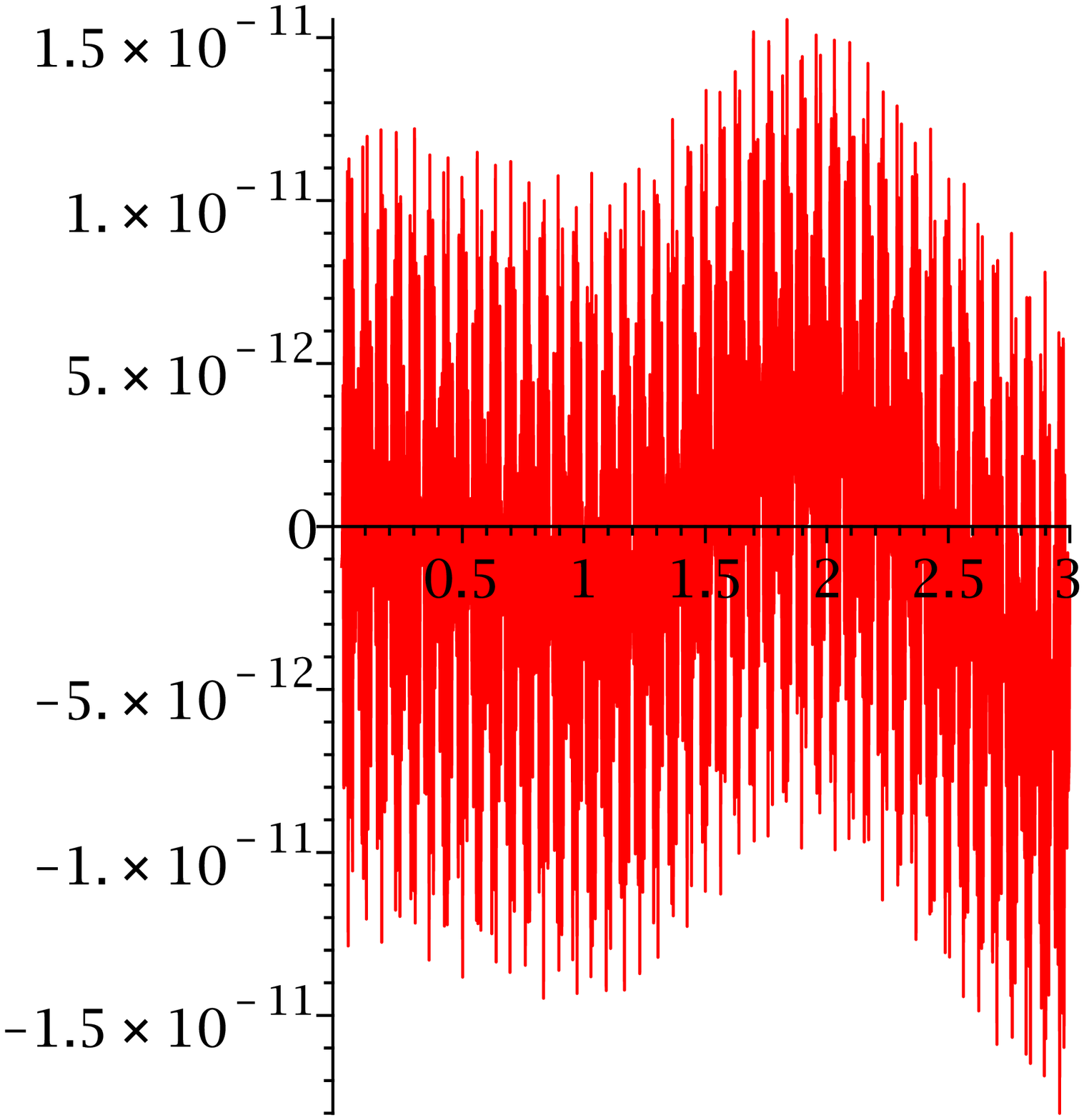}}\\
\caption{The top row: the real parts of error function with $s = 0$
(the left) and $s = 1$ (the right) for $y_1$ with $\omega = 100$.
The middle row: the real parts of error function with $s = 2$ (the
left) and $s = 3$ (the right) for $y_1$ with $\omega = 100$. The
third row: the real parts of error function with $s = 0$ (the left)
and $s = 1$ (the right) for $y_1$ with $\omega = 1000$. The fourth
row: the real parts of error function with $s = 2$ (the left) and $s
= 3$ (the right)for $y_1$ with $\omega = 1000$.}
\end{figure}

\begin{figure}[htbp]
\centering
{\label{fig:4.7}\includegraphics[width=6cm,height=5cm]{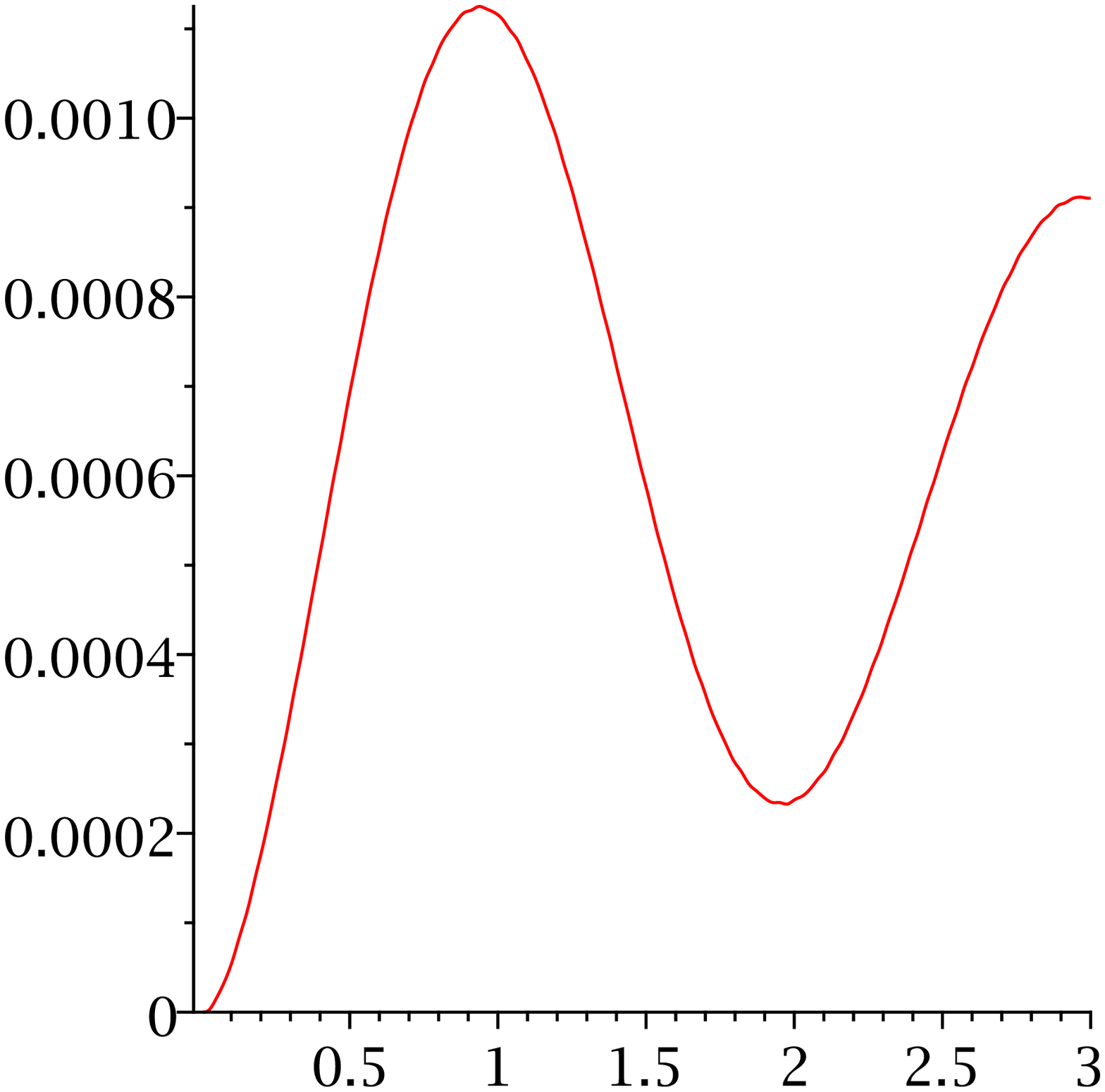}}
\quad {\label{fig:
5.3}\includegraphics[width=6cm,height=5cm]{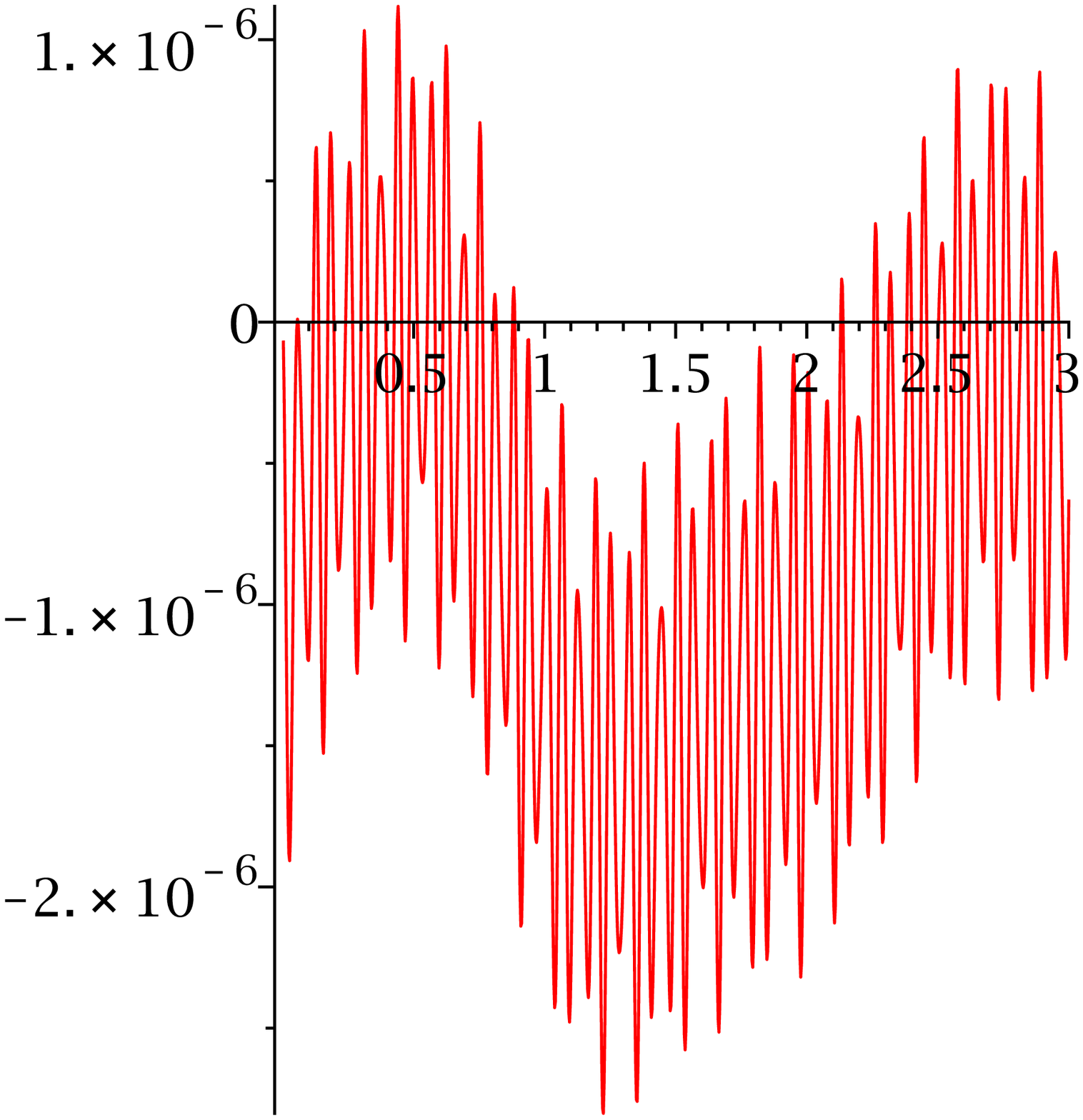}}
\\
{\label{fig:4.7}\includegraphics[width=6cm,height=5cm]{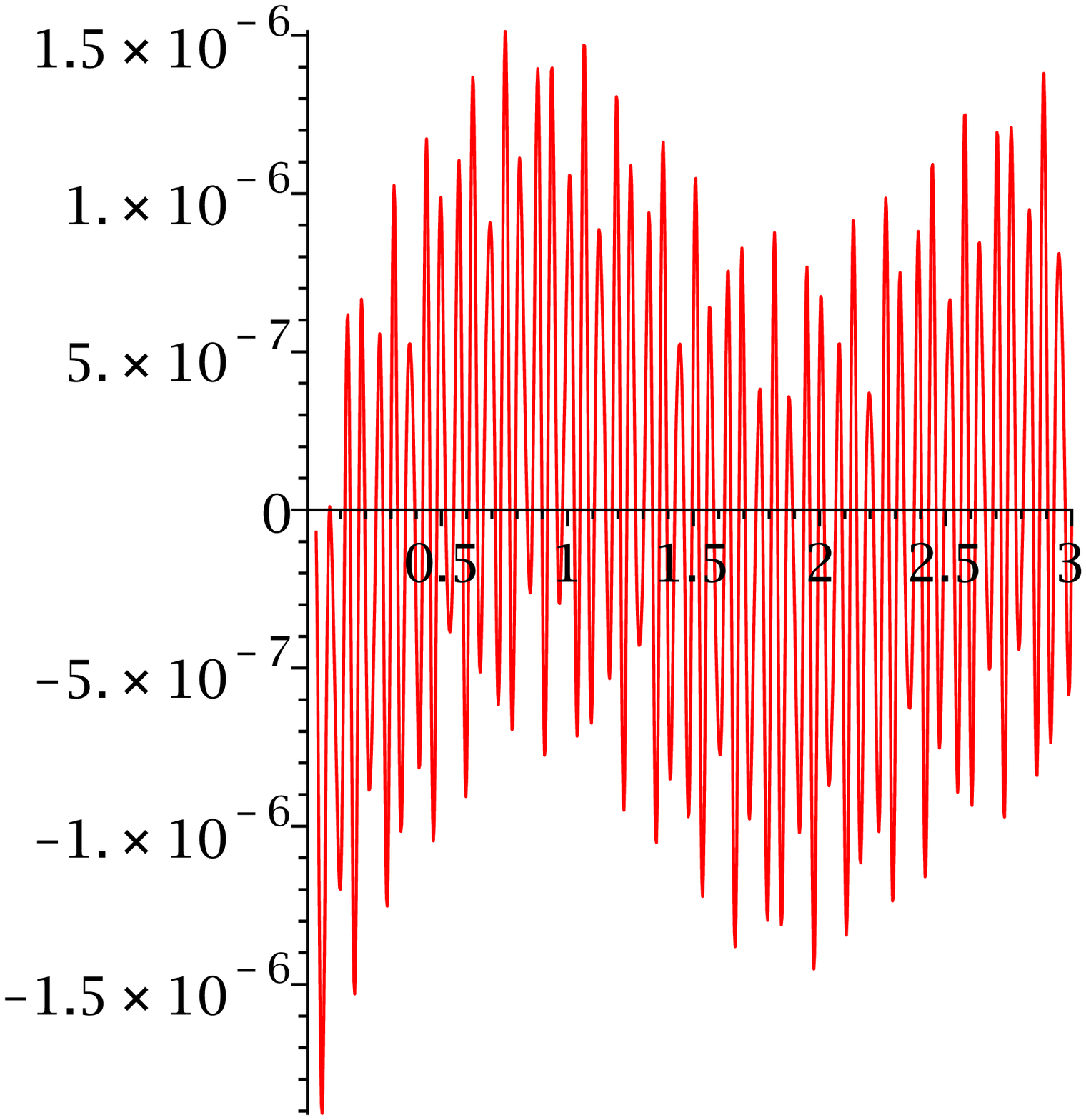}}
\quad
{\label{fig:4.7}\includegraphics[width=6cm,height=5cm]{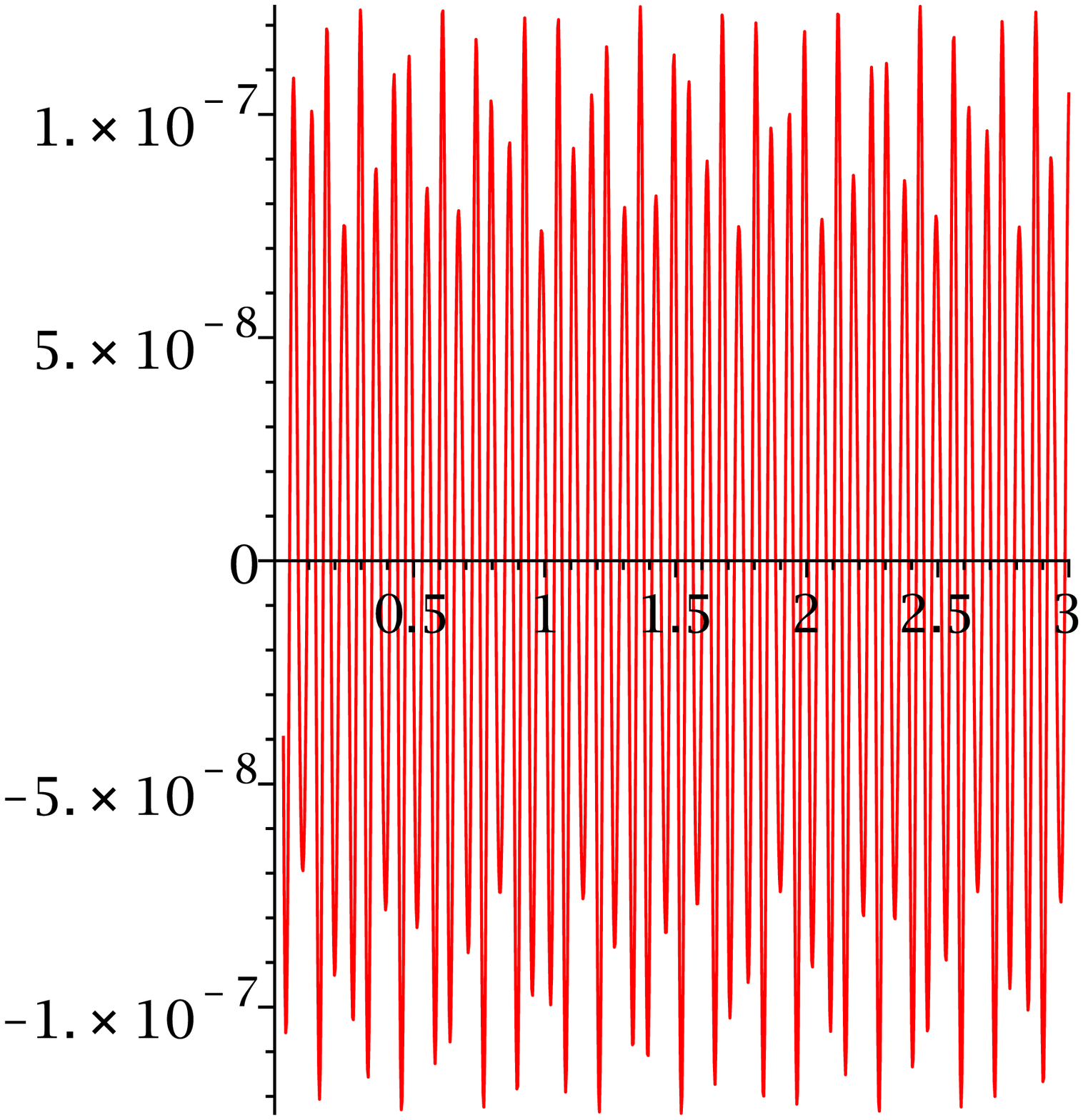}}
\\
{\label{fig:4.7}\includegraphics[width=6cm,height=5cm]{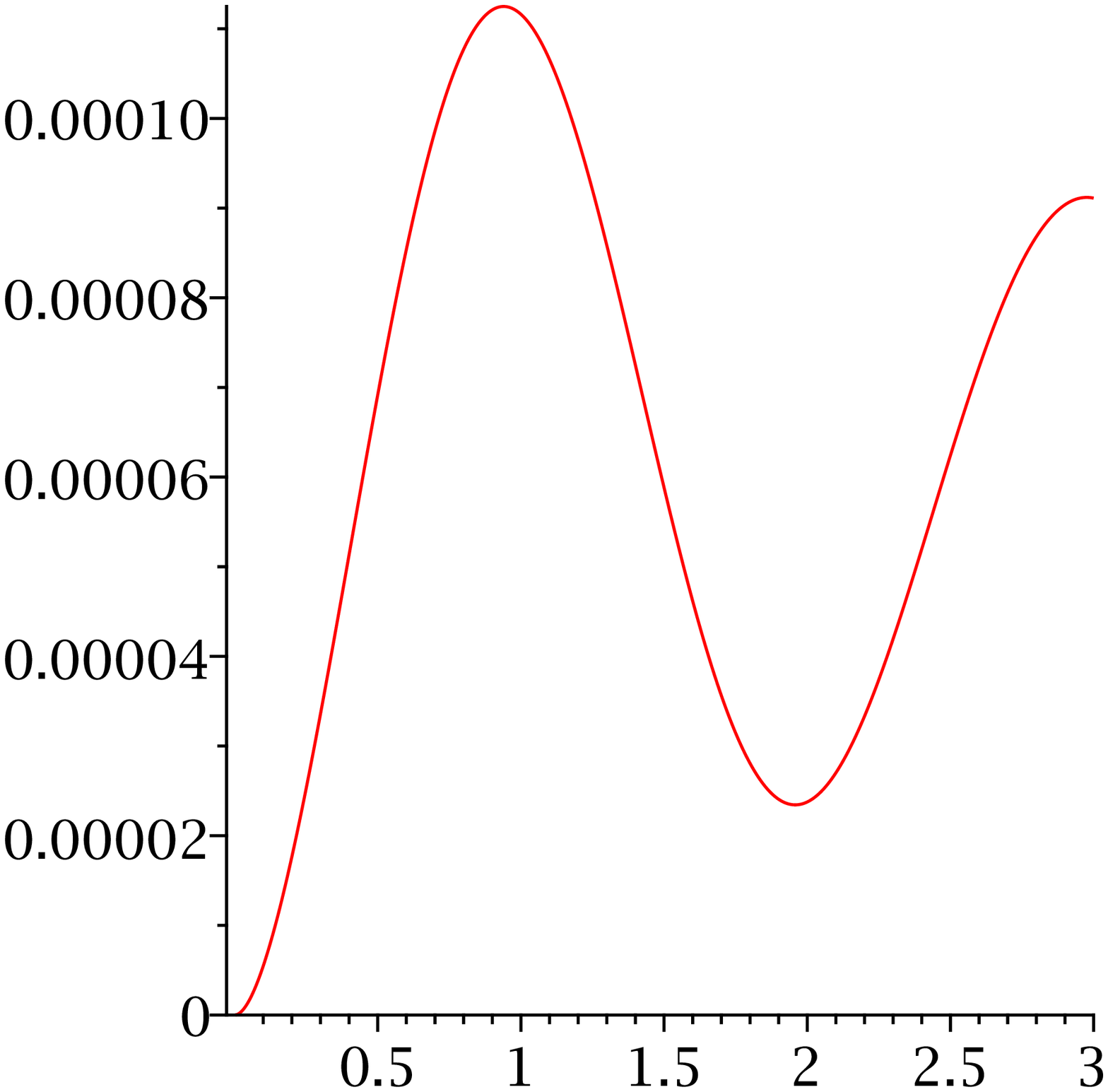}}
\quad
{\label{fig:4.7}\includegraphics[width=6cm,height=5cm]{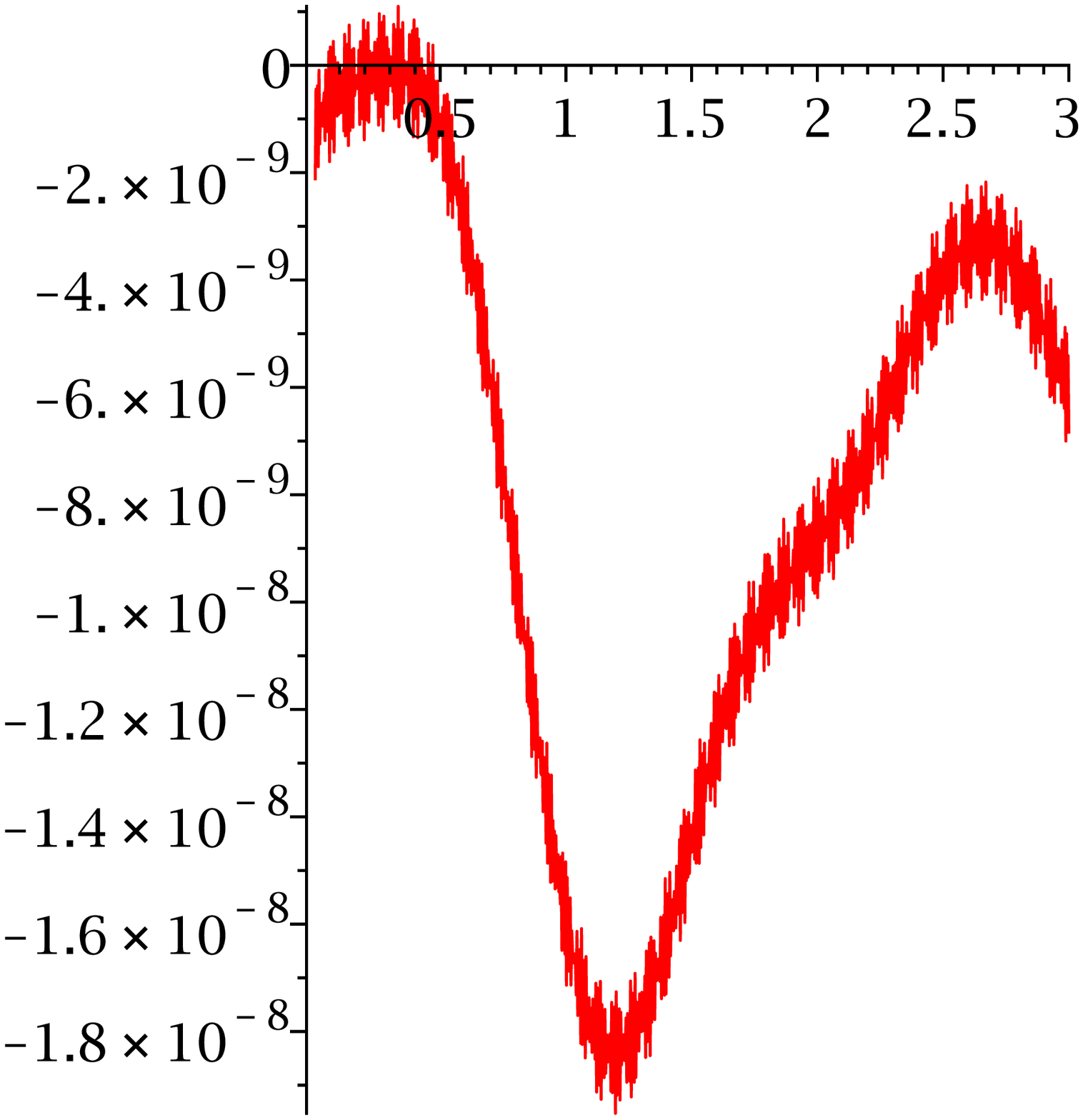}}\\
{\label{fig:4.7}\includegraphics[width=6cm,height=5cm]{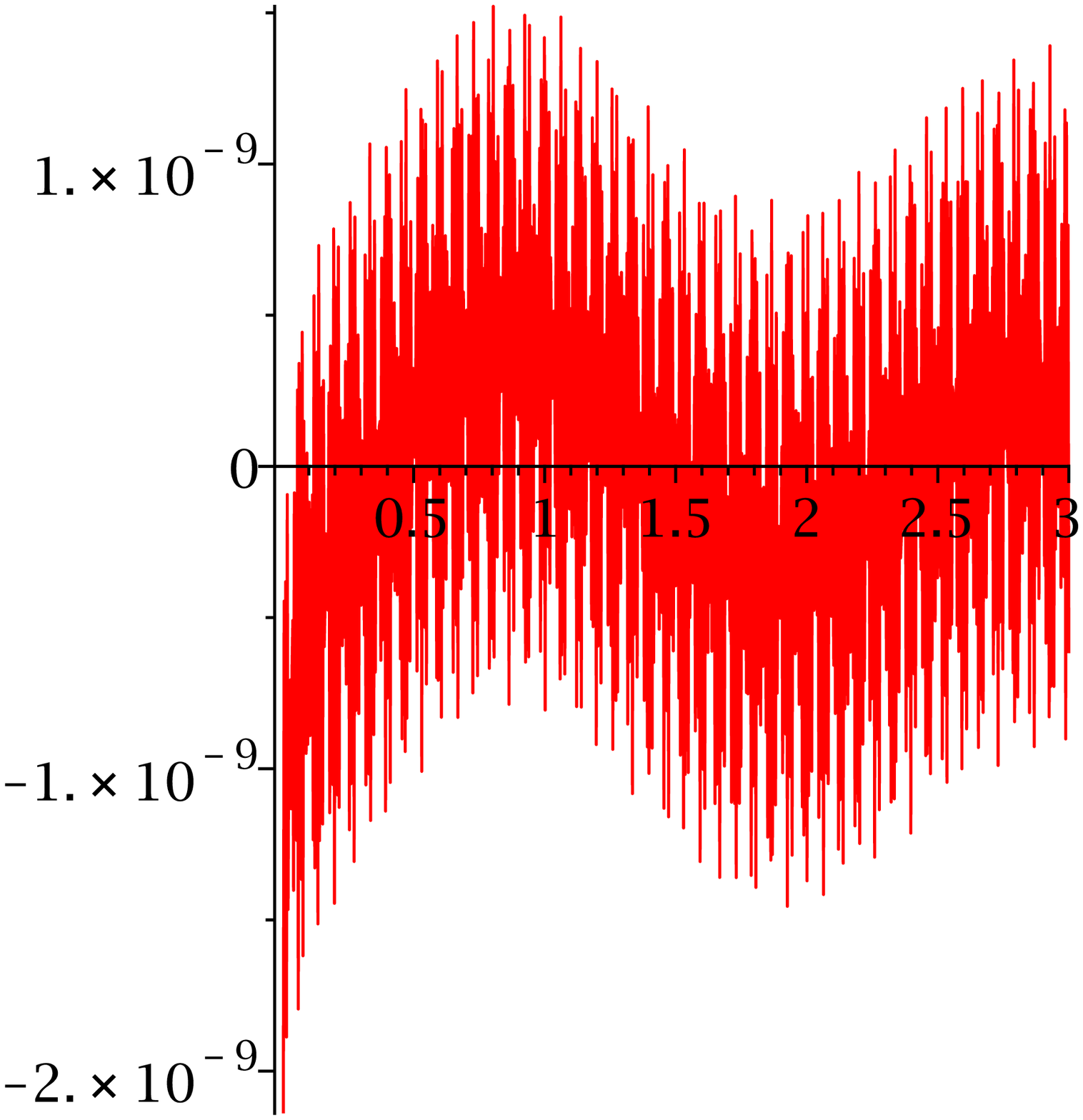}}
\quad
{\label{fig:4.7}\includegraphics[width=6cm,height=5cm]{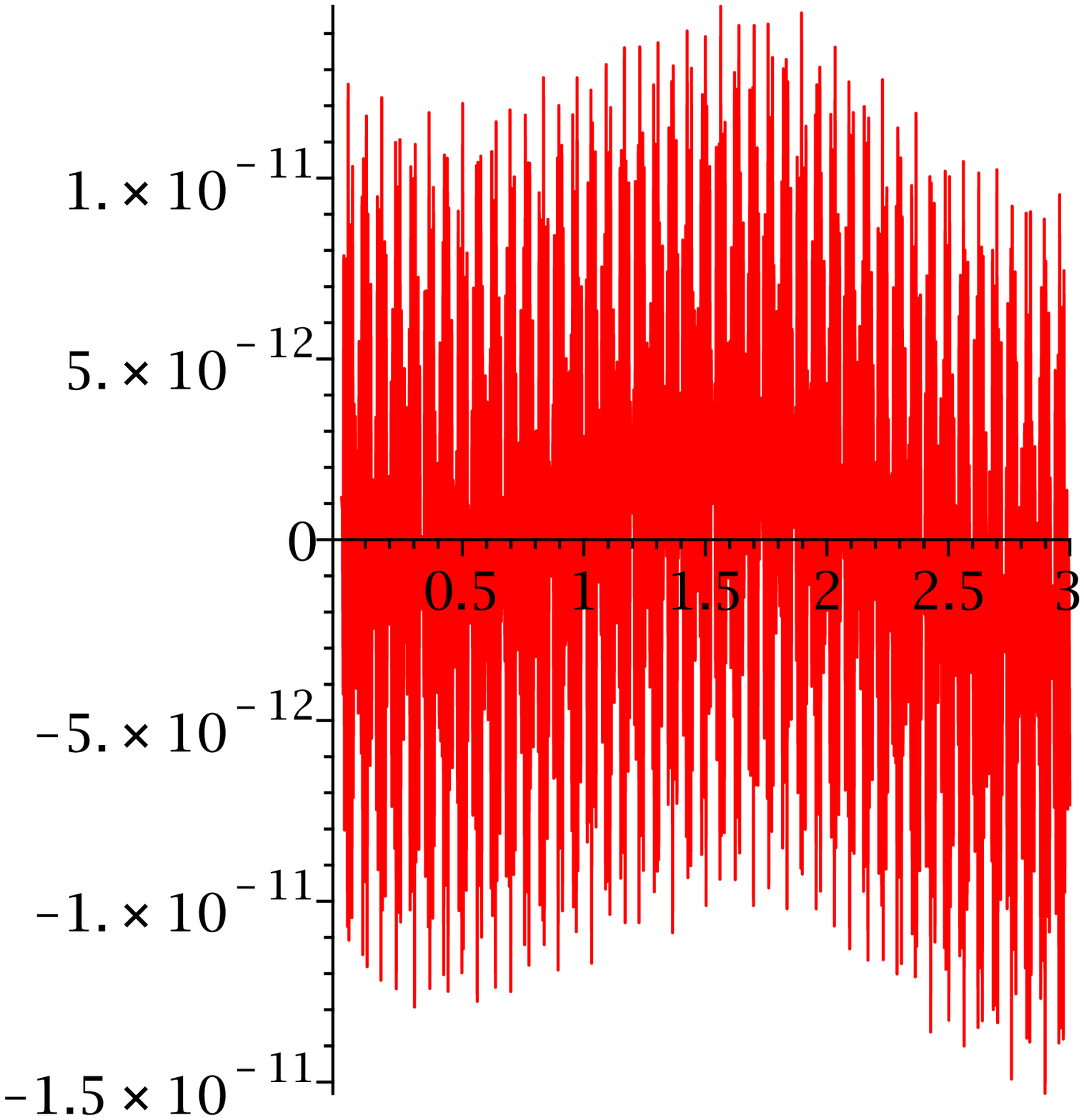}}\\
\caption{The top row: the real parts of error function with $s = 0$
(the left) and $s = 1$ (the right) for $y_2$ with $\omega = 100$.
The middle row: the real parts of error function with $s = 2$ (the
left) and $s = 3$ (the right) for $y_2$ with $\omega = 100$. The
third row: the real parts of error function with $s = 0$ (the left)
and $s = 1$ (the right) for $y_2$ with $\omega = 1000$. The fourth
row: the real parts of error function with $s = 2$ (the left) and $s
= 3$ (the right)for $y_2$ with $\omega = 1000$.}
\end{figure}

\begin{figure}[htbp]
\centering
{\label{fig:4.8}\includegraphics[width=6cm,height=5cm]{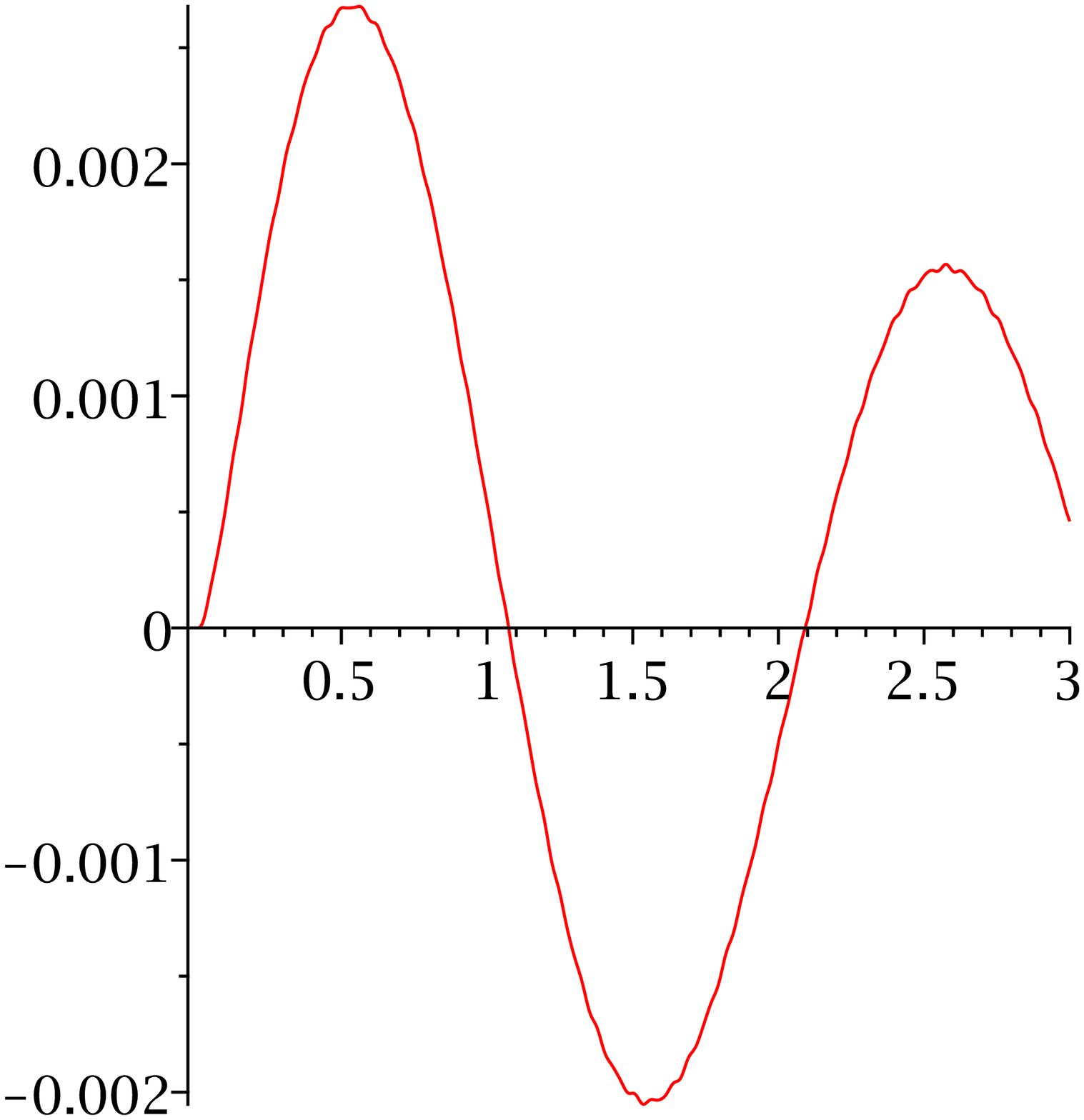}}
\quad {\label{fig:
5.4}\includegraphics[width=6cm,height=5cm]{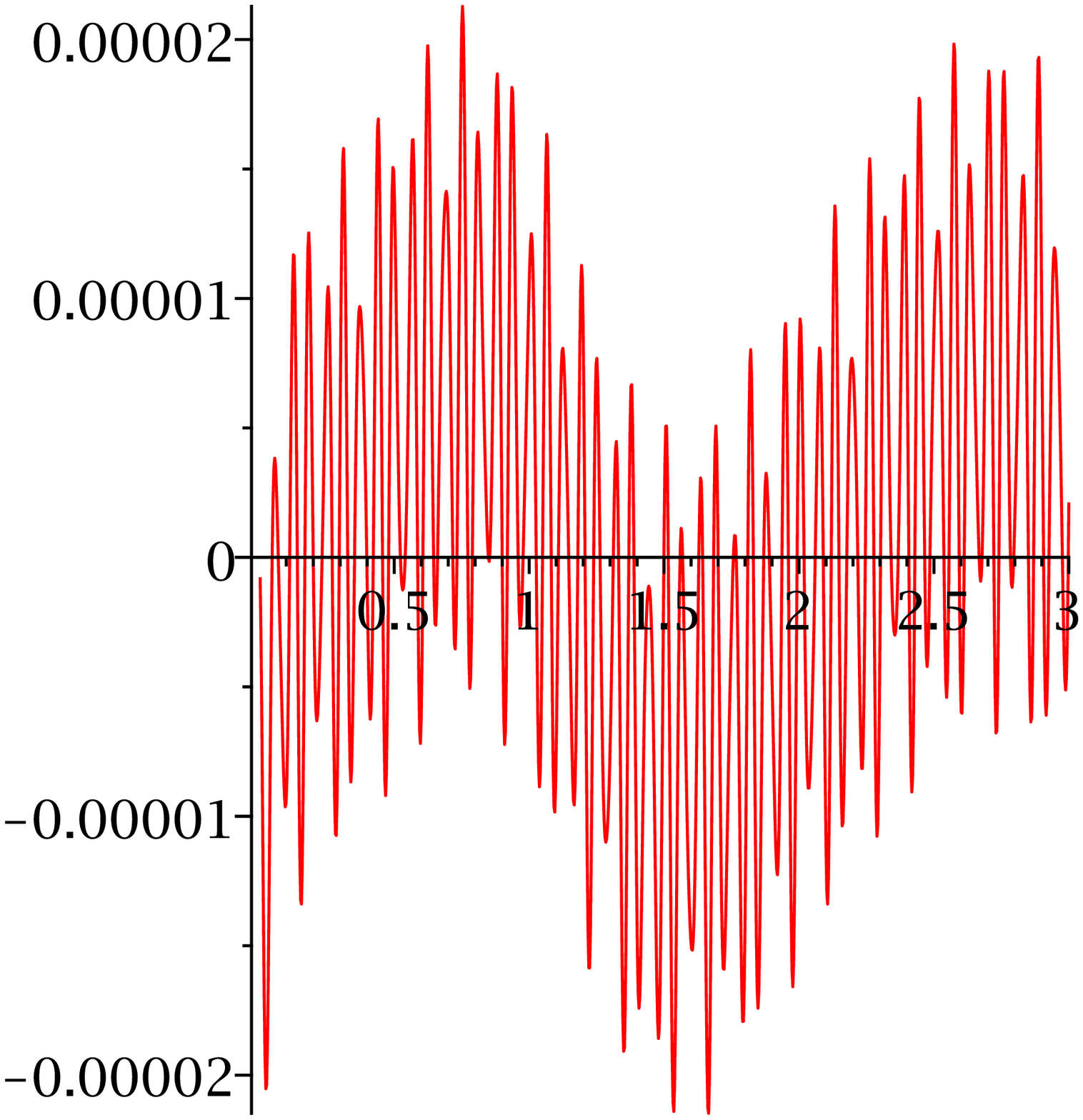}}
\\
{\label{fig:4.8}\includegraphics[width=6cm,height=5cm]{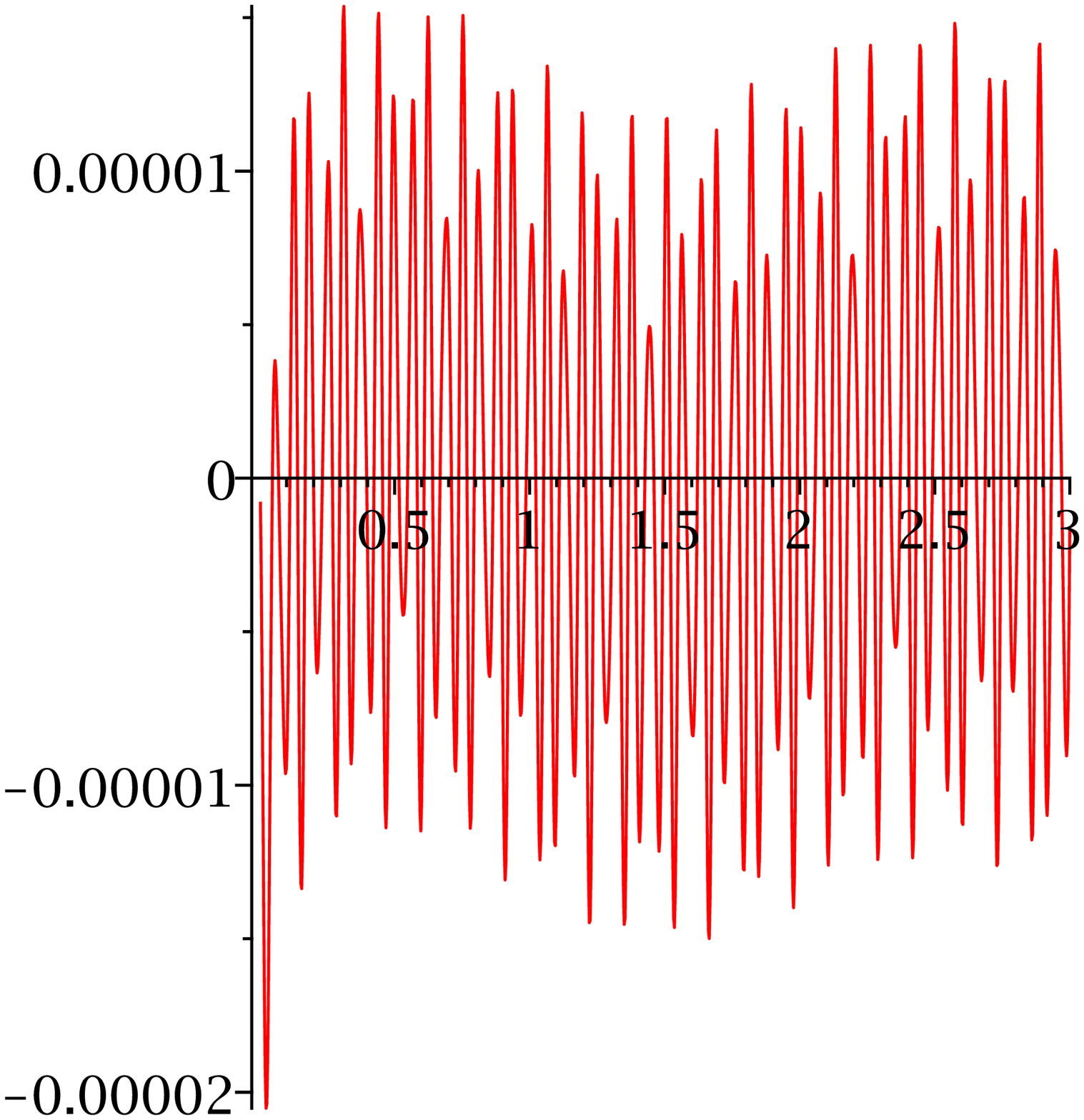}}
\quad
{\label{fig:4.8}\includegraphics[width=6cm,height=5cm]{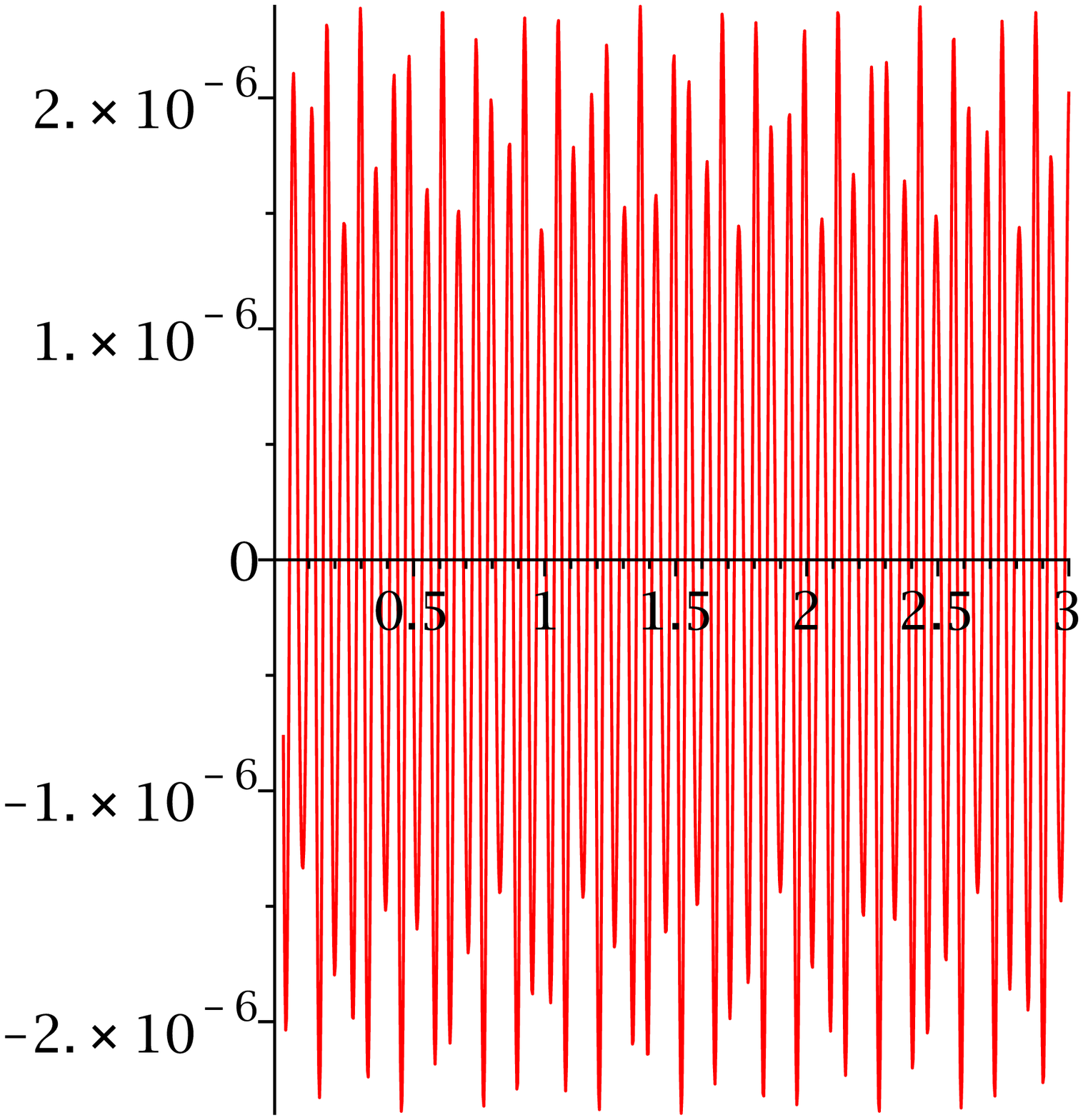}}
\\
{\label{fig:4.8}\includegraphics[width=6cm,height=5cm]{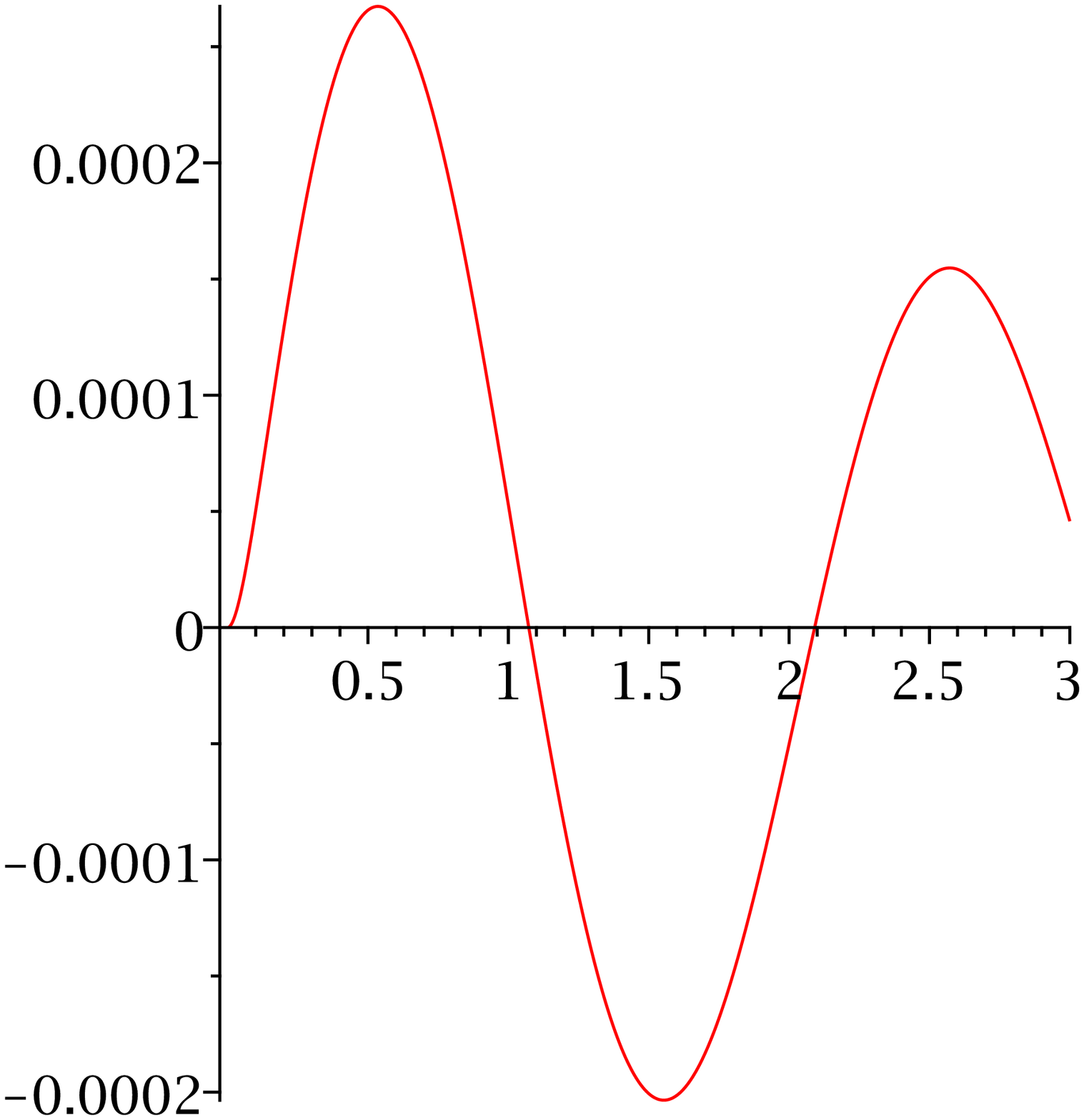}}
\quad
{\label{fig:4.8}\includegraphics[width=6cm,height=5cm]{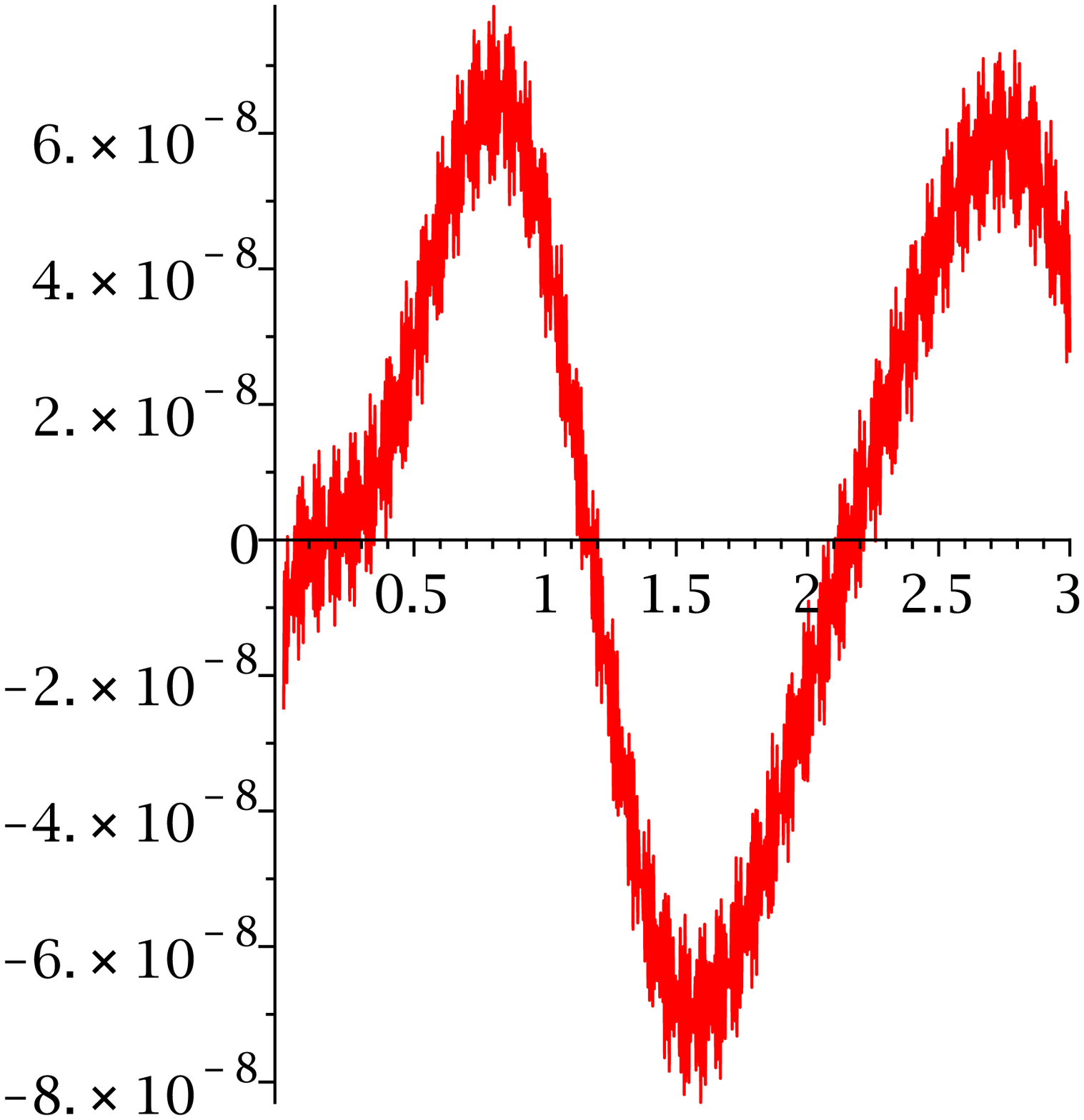}}\\
{\label{fig:4.8}\includegraphics[width=6cm,height=5cm]{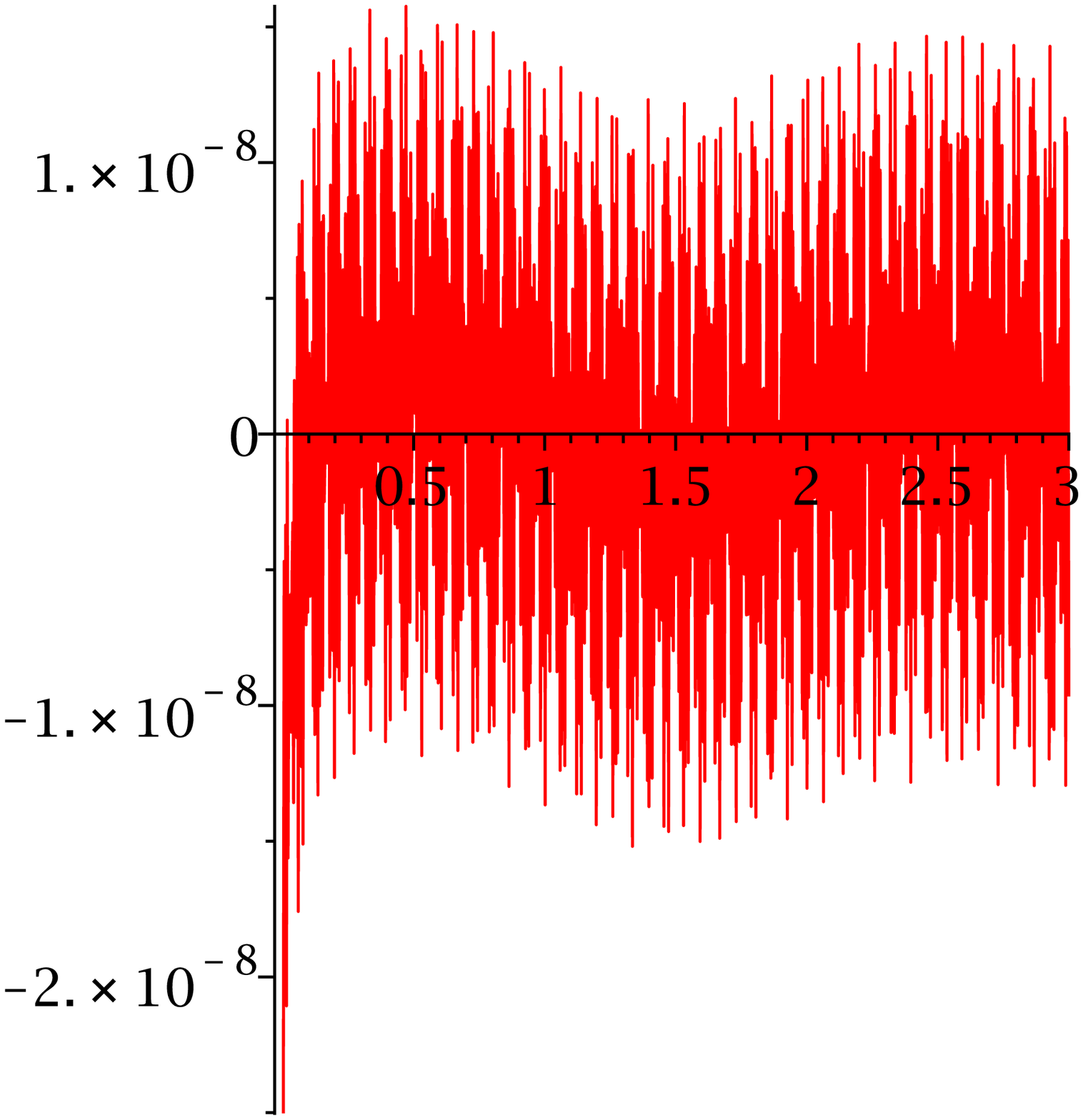}}
\quad
{\label{fig:4.8}\includegraphics[width=6cm,height=5cm]{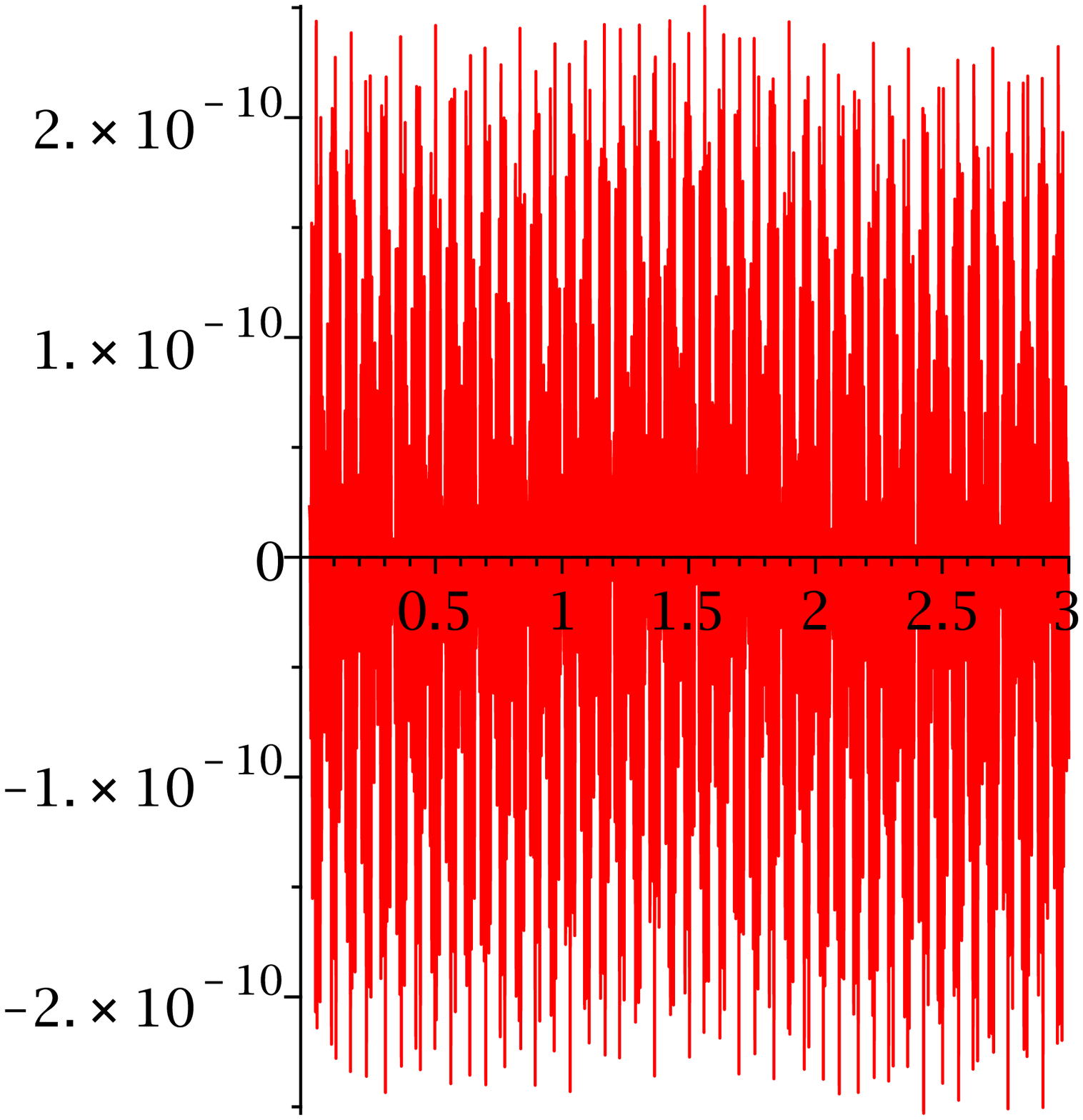}}\\
\caption{The top row: the real parts of error function with $s = 0$
(the left) and $s = 1$ (the right) for $y_3$ with $\omega = 100$.
The middle row: the real parts of error function with $s = 2$ (the
left) and $s = 3$ (the right) for $y_3$ with $\omega = 100$. The
third row: the real parts of error function with $s = 0$ (the left)
and $s = 1$ (the right) for $y_3$ with $\omega = 1000$. The fourth
row: the real parts of error function with $s = 2$ (the left) and $s
= 3$ (the right)for $y_3$ with $\omega = 1000$.}
\end{figure}

\begin{figure}[htbp]
\centering
{\label{fig:4.9}\includegraphics[width=6cm,height=5cm]{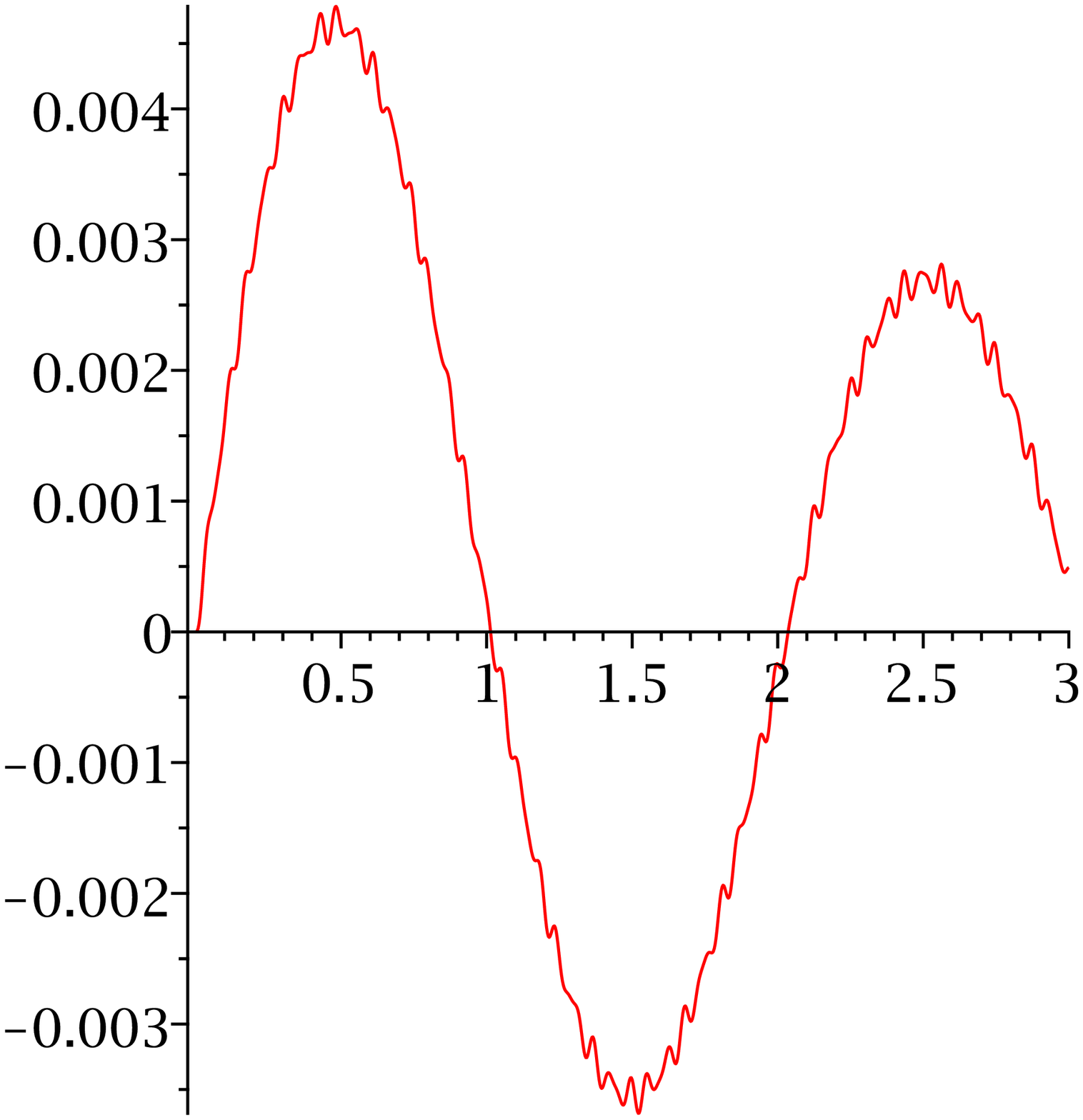}}
\quad
{\label{fig:4.9}\includegraphics[width=6cm,height=5cm]{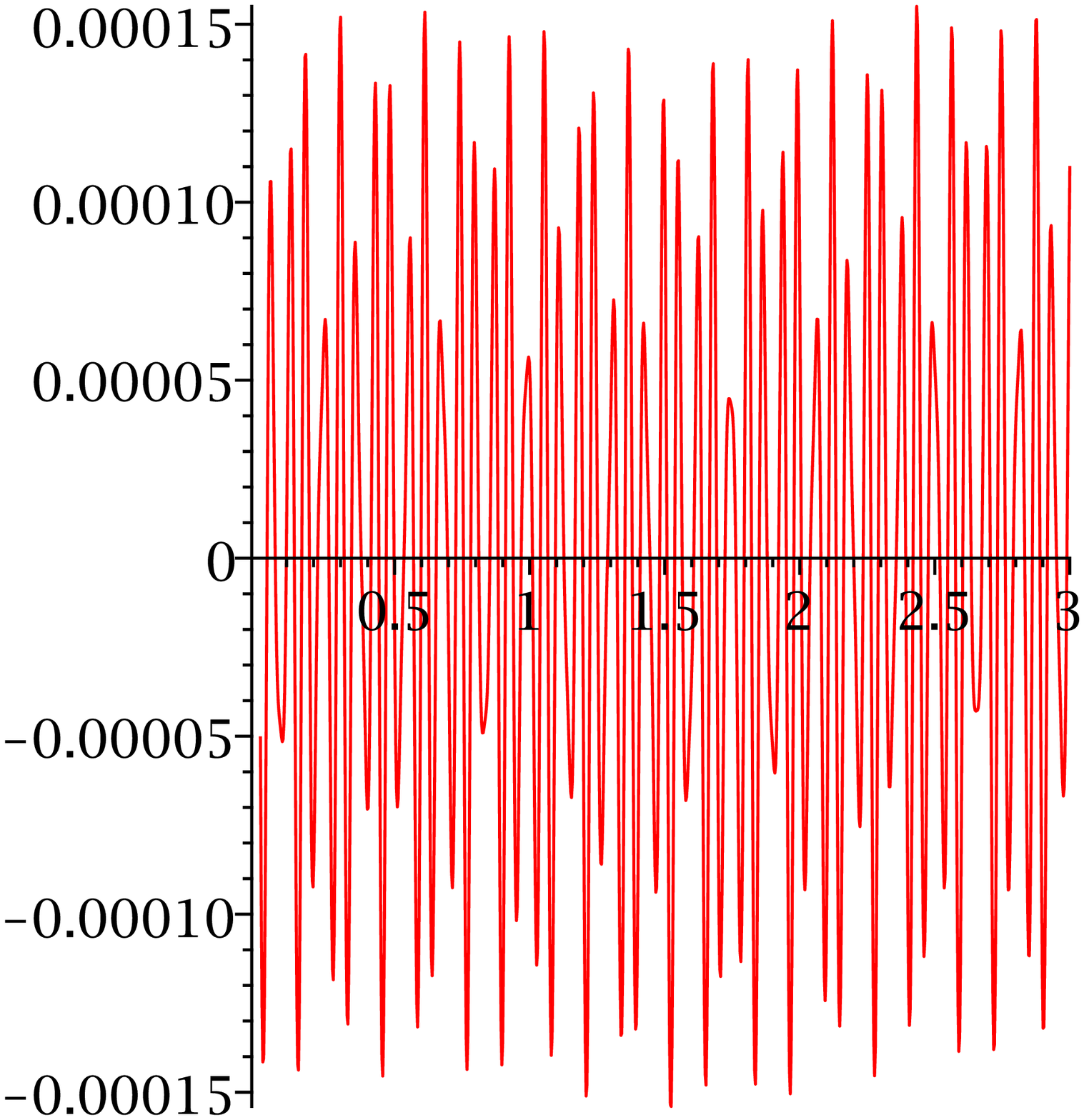}}
\\
{\label{fig:4.9}\includegraphics[width=6cm,height=5cm]{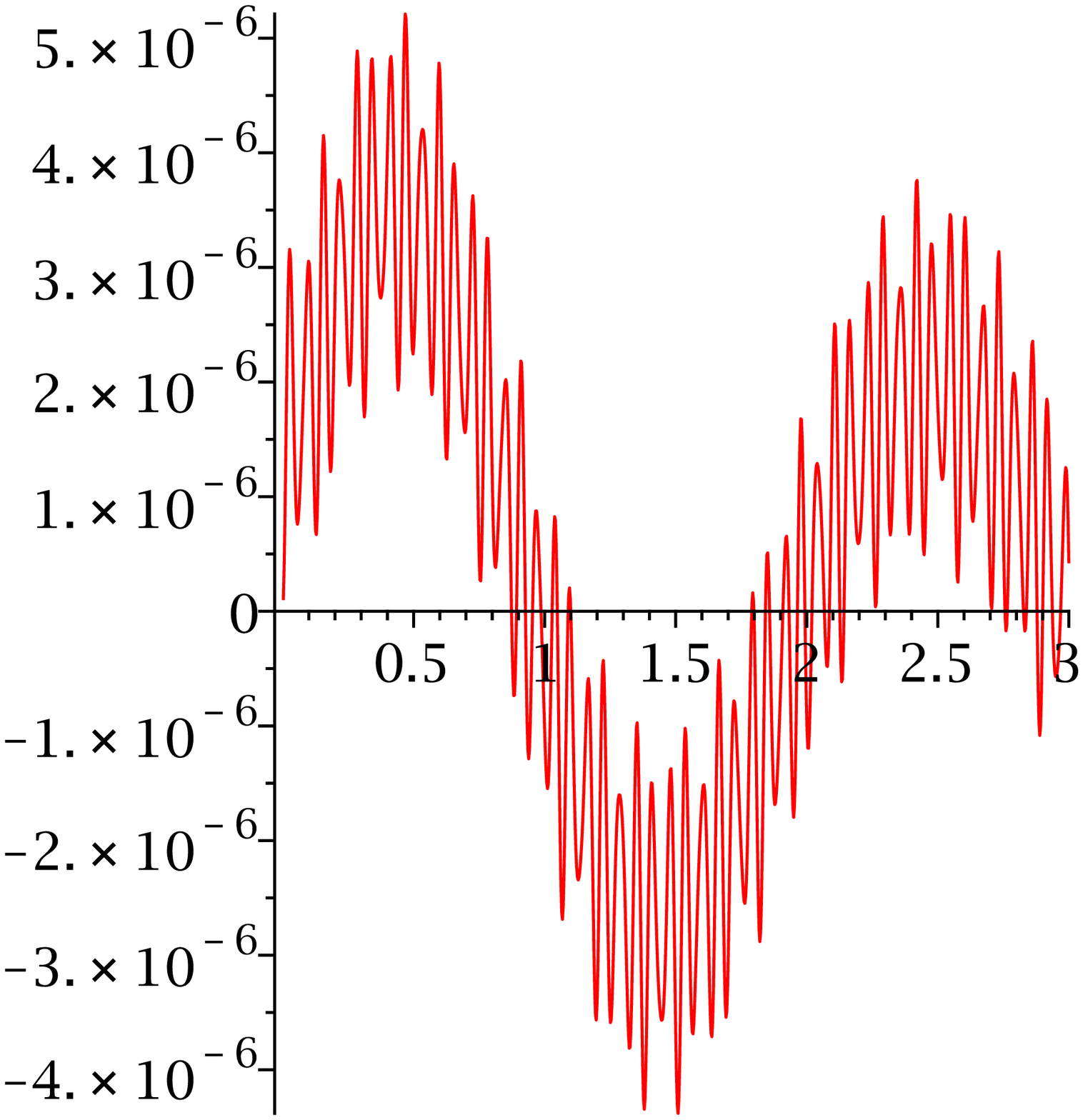}}
\quad
{\label{fig:4.9}\includegraphics[width=6cm,height=5cm]{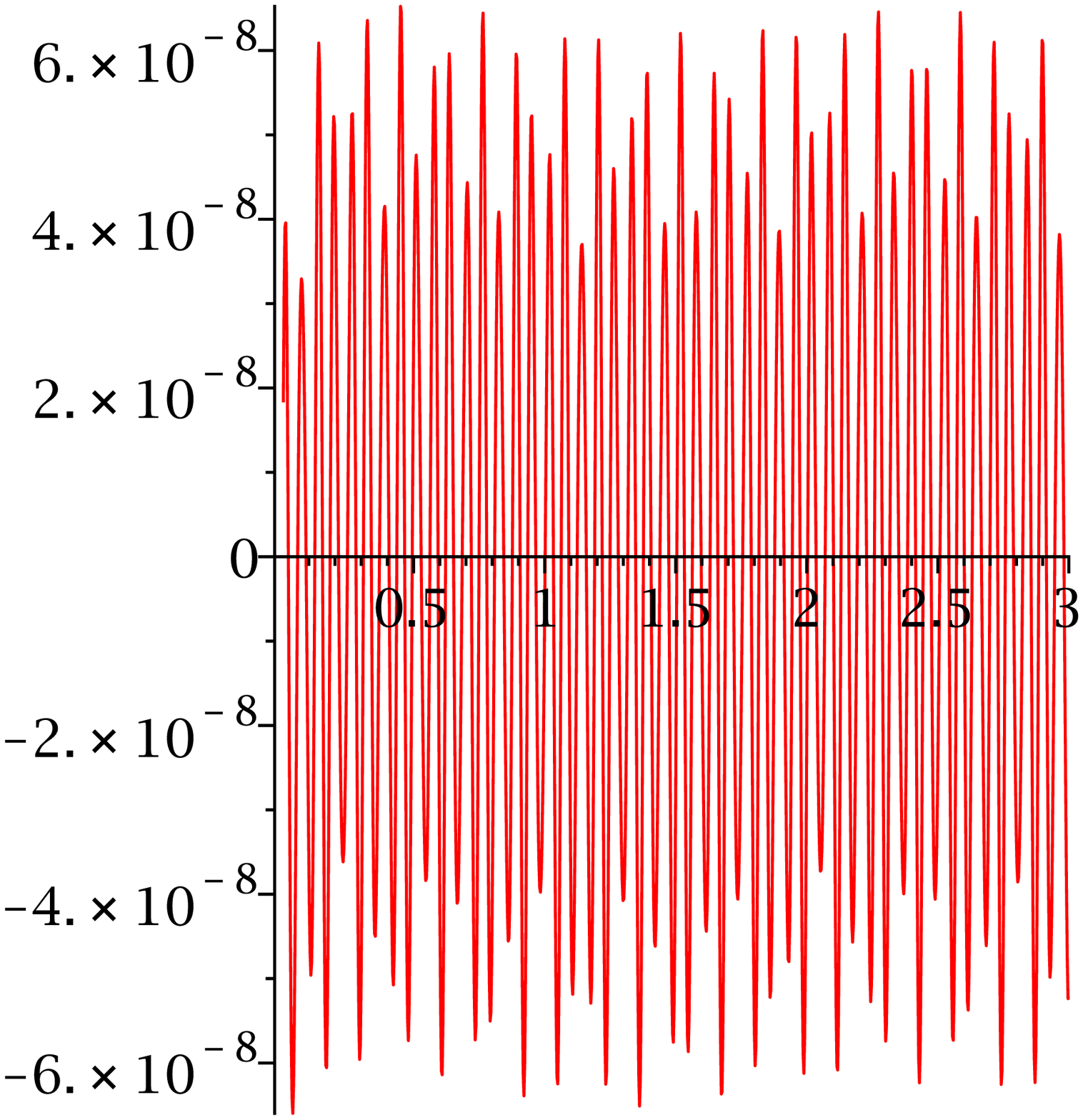}}
\\
{\label{fig:4.9}\includegraphics[width=6cm,height=5cm]{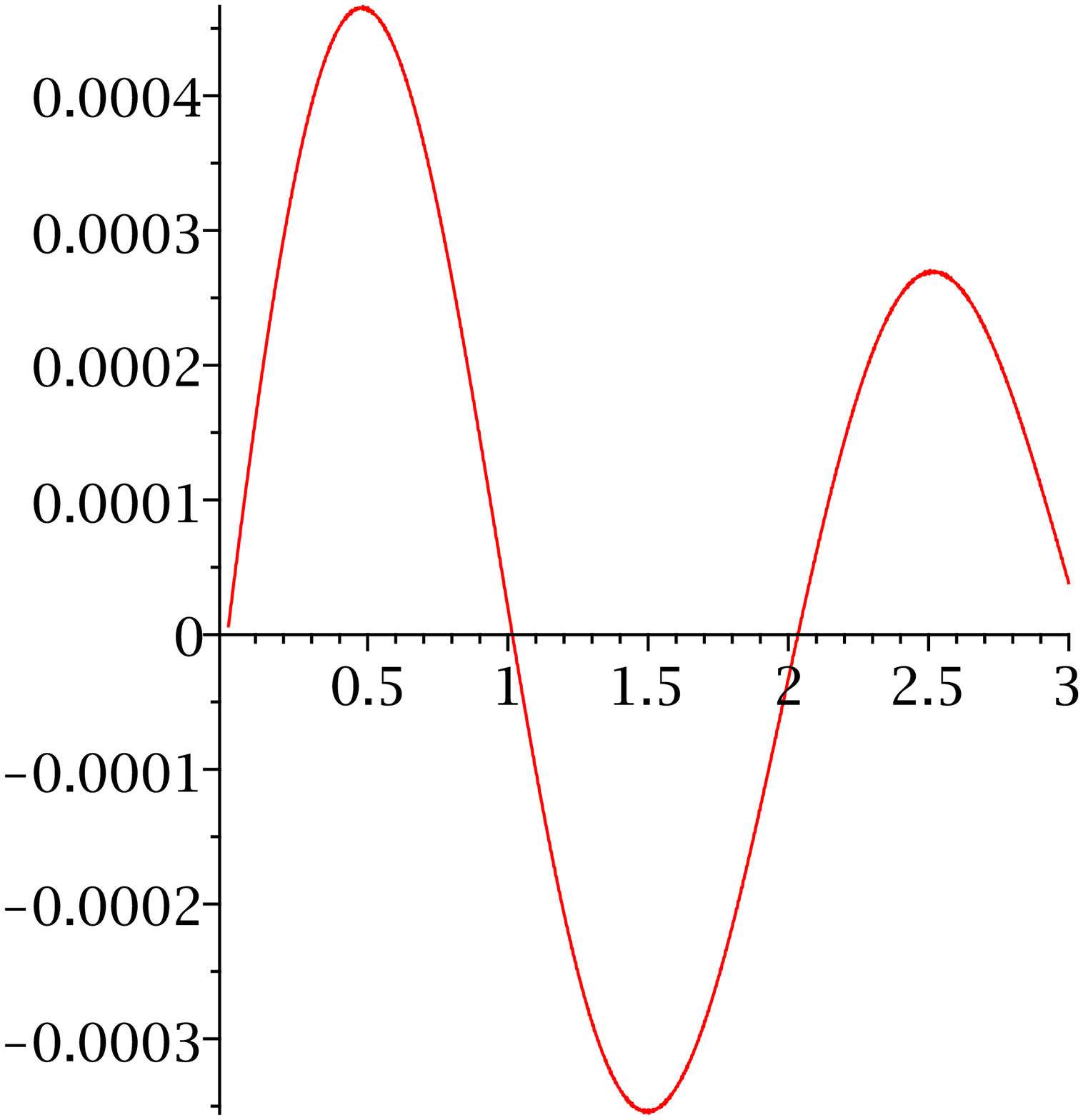}}
\quad
{\label{fig:4.9}\includegraphics[width=6cm,height=5cm]{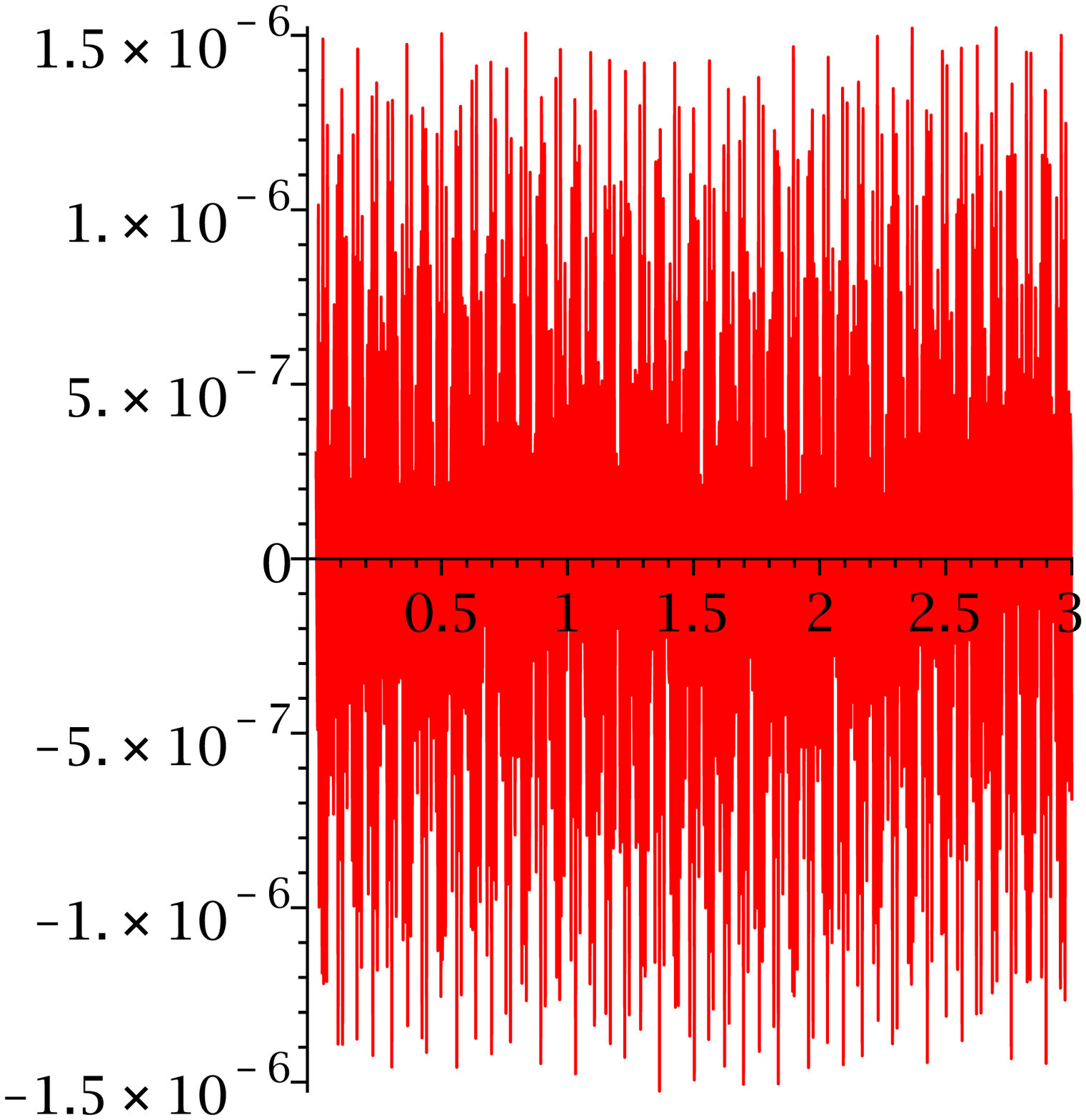}}\\
{\label{fig:4.9}\includegraphics[width=6cm,height=5cm]{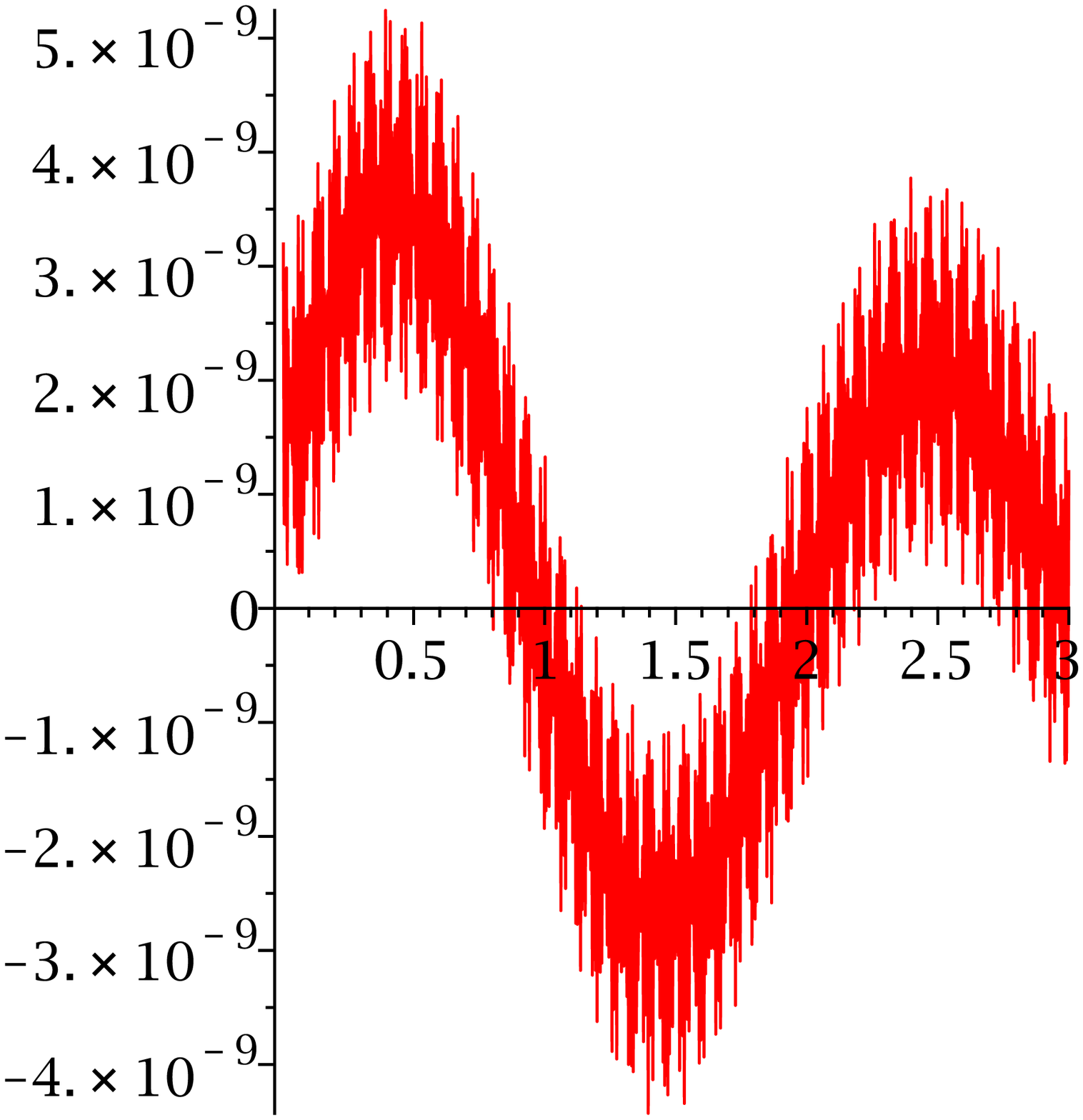}}
\quad
{\label{fig:4.9}\includegraphics[width=6cm,height=5cm]{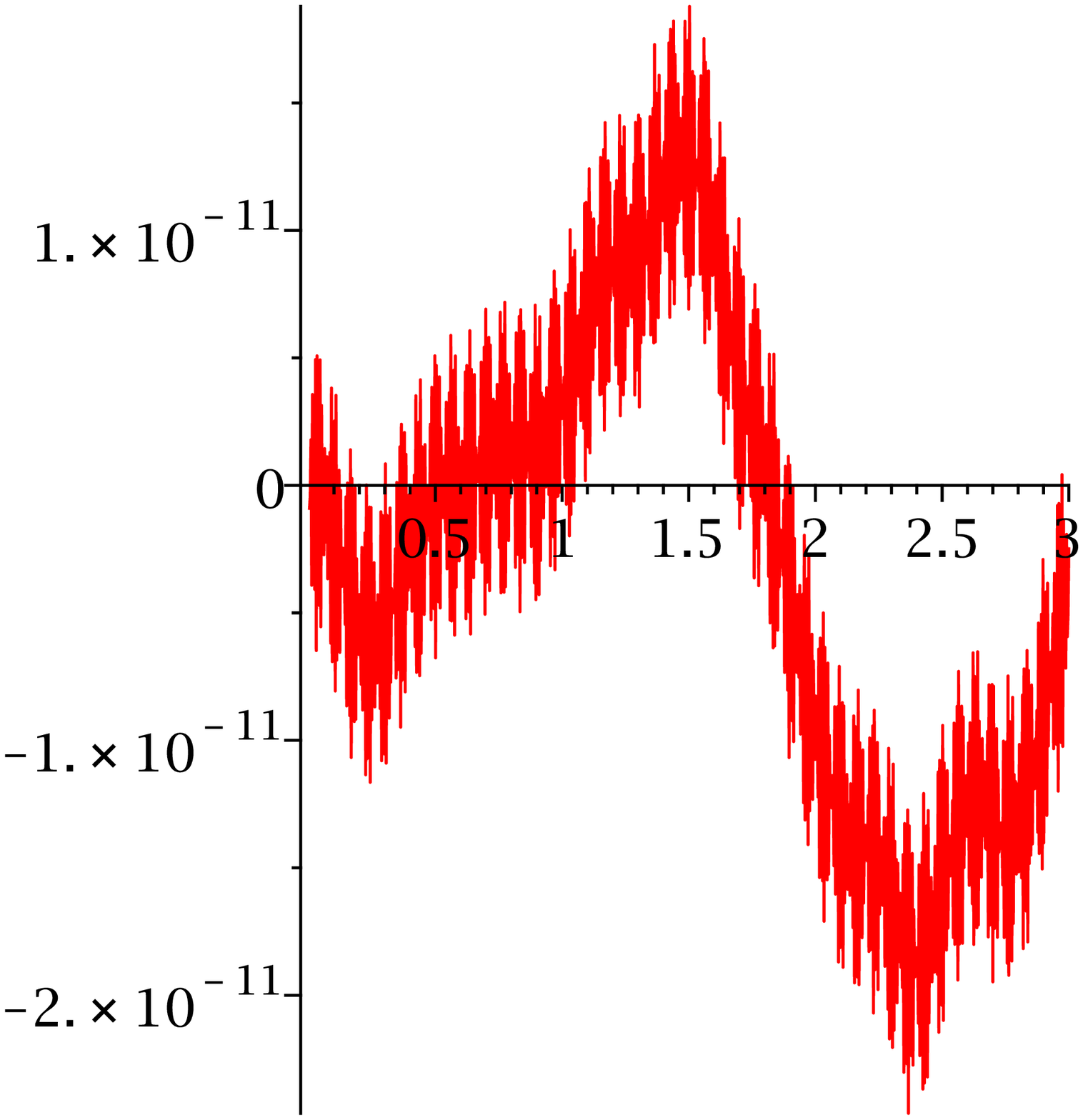}}\\
\caption{The top row: the real parts of error function with $s = 0$
(the left) and $s = 1$ (the right) for $y_4$ with $\omega = 100$.
The middle row: the real parts of error function with $s = 2$ (the
left) and $s = 3$ (the right) for $y_4$ with $\omega = 100$. The
third row: the real parts of error function with $s = 0$ (the left)
and $s = 1$ (the right) for $y_4$ with $\omega = 1000$. The fourth
row: the real parts of error function with $s = 2$ (the left) and $s
= 3$ (the right)for $y_4$ with $\omega = 1000$.}
\end{figure}

\begin{figure}[htbp]
\centering
{\label{fig:4.10}\includegraphics[width=6cm,height=5cm]{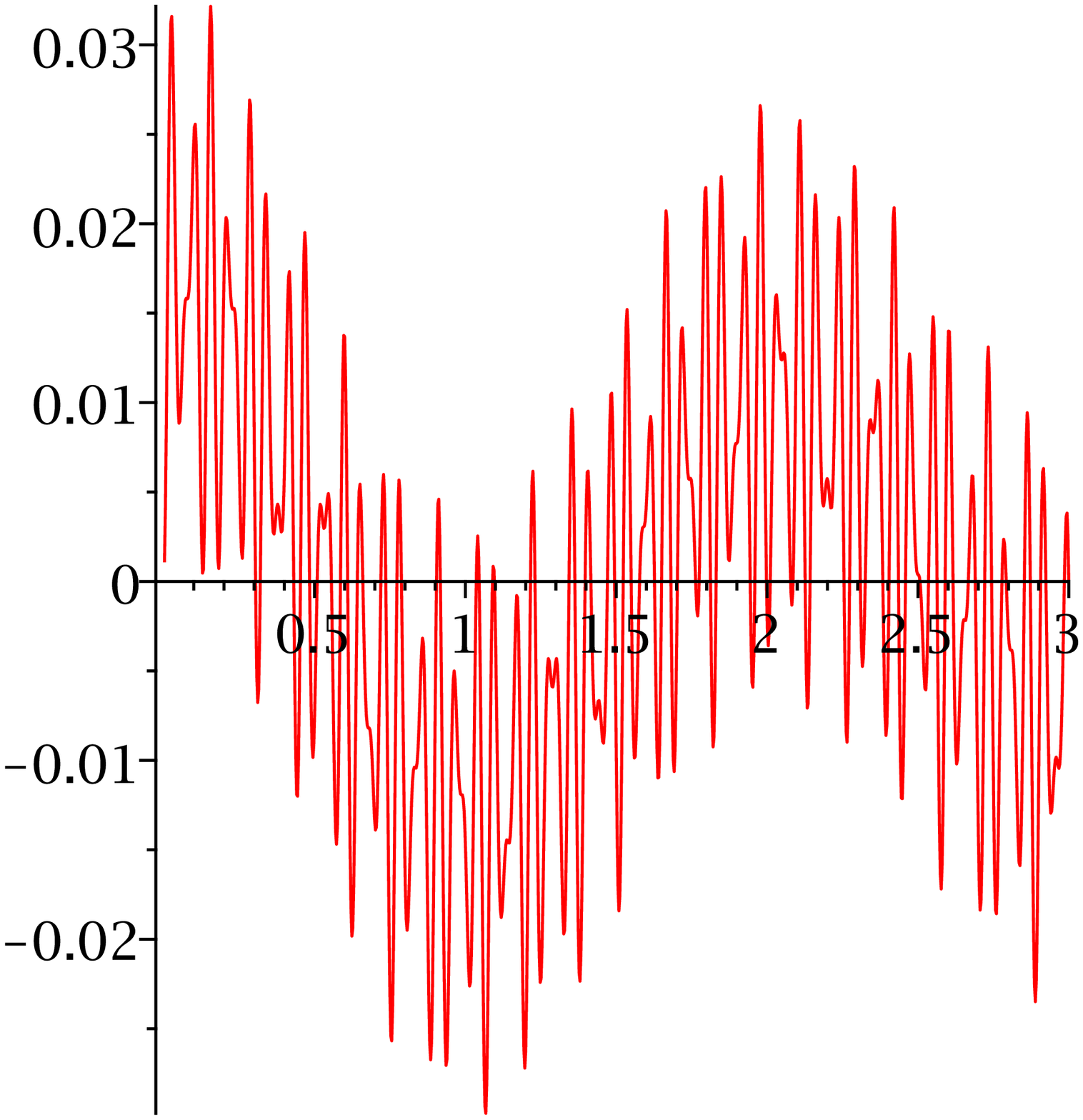}}
\quad
{\label{fig:4.10}\includegraphics[width=6cm,height=5cm]{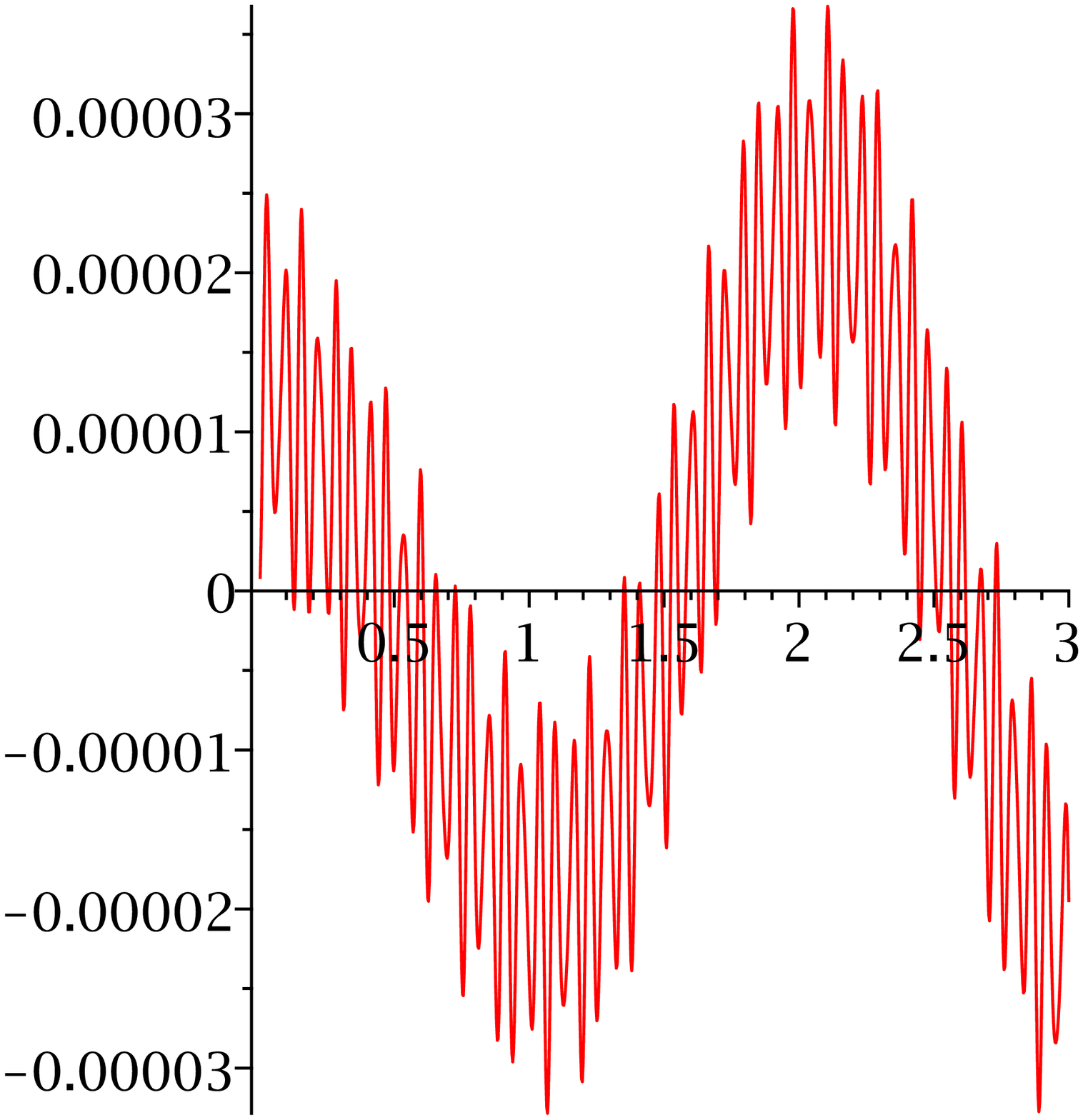}}
\\
{\label{fig:4.10}\includegraphics[width=6cm,height=5cm]{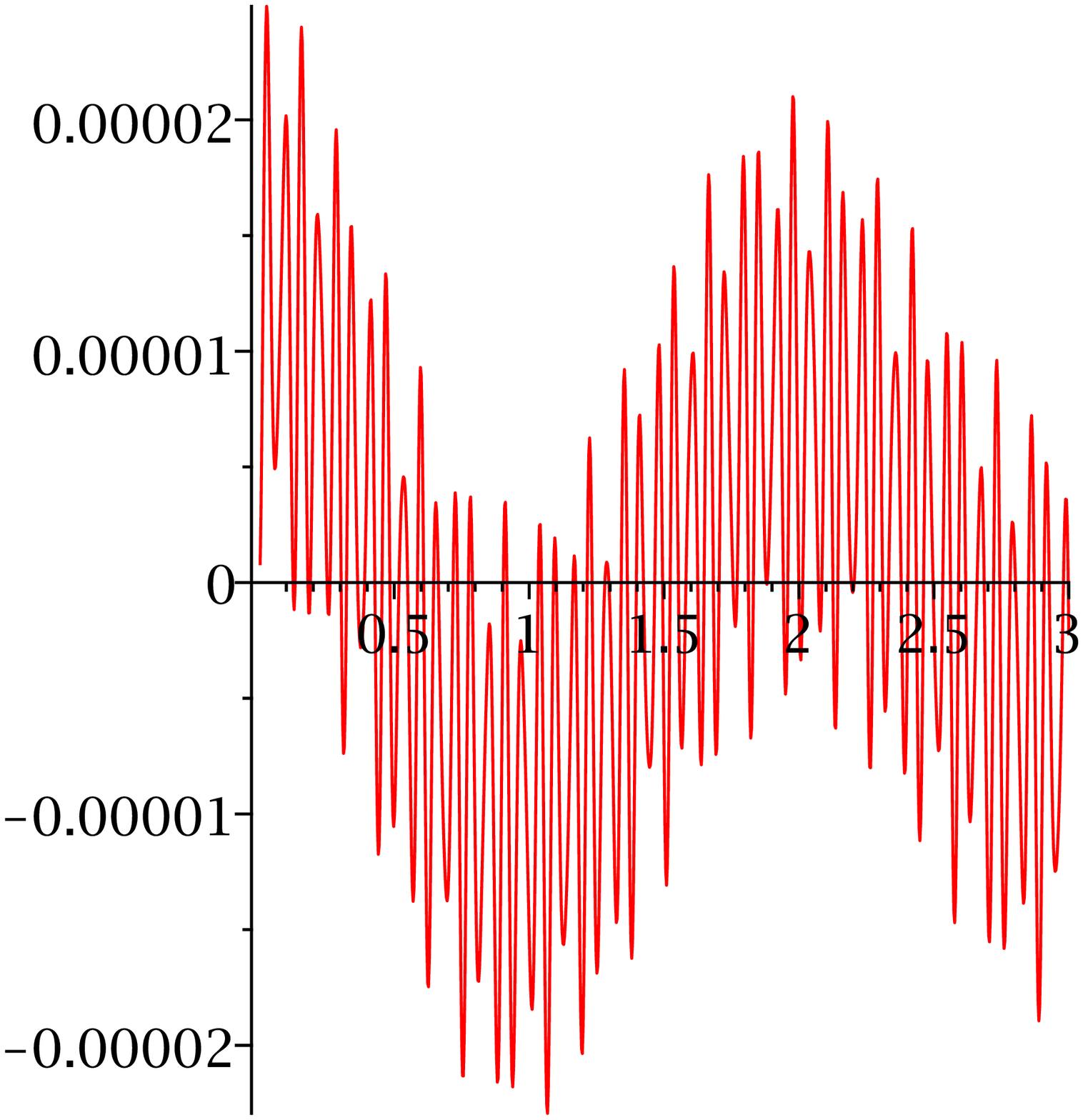}}
\quad
{\label{fig:4.10}\includegraphics[width=6cm,height=5cm]{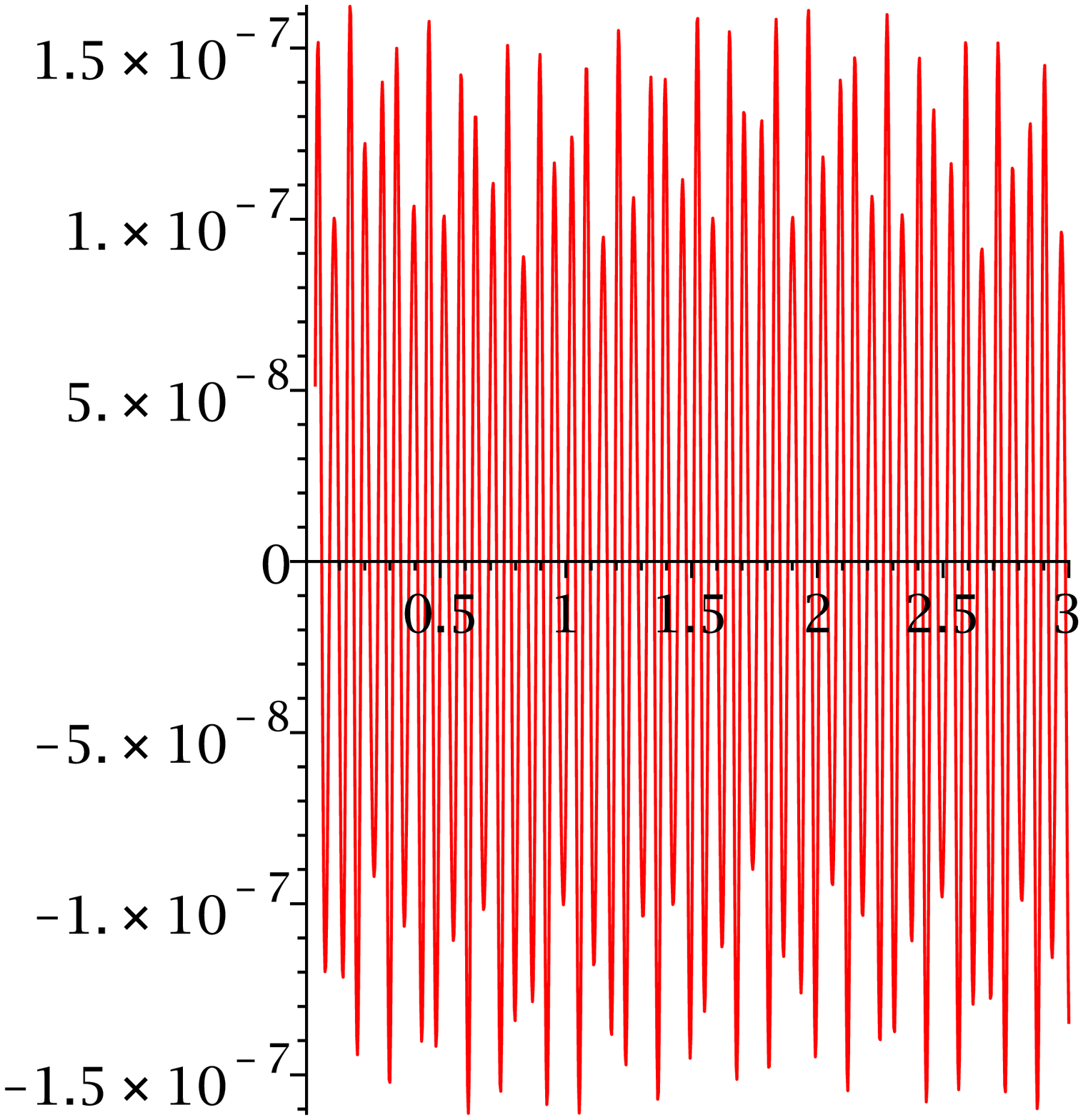}}
\\
{\label{fig:4.10}\includegraphics[width=6cm,height=5cm]{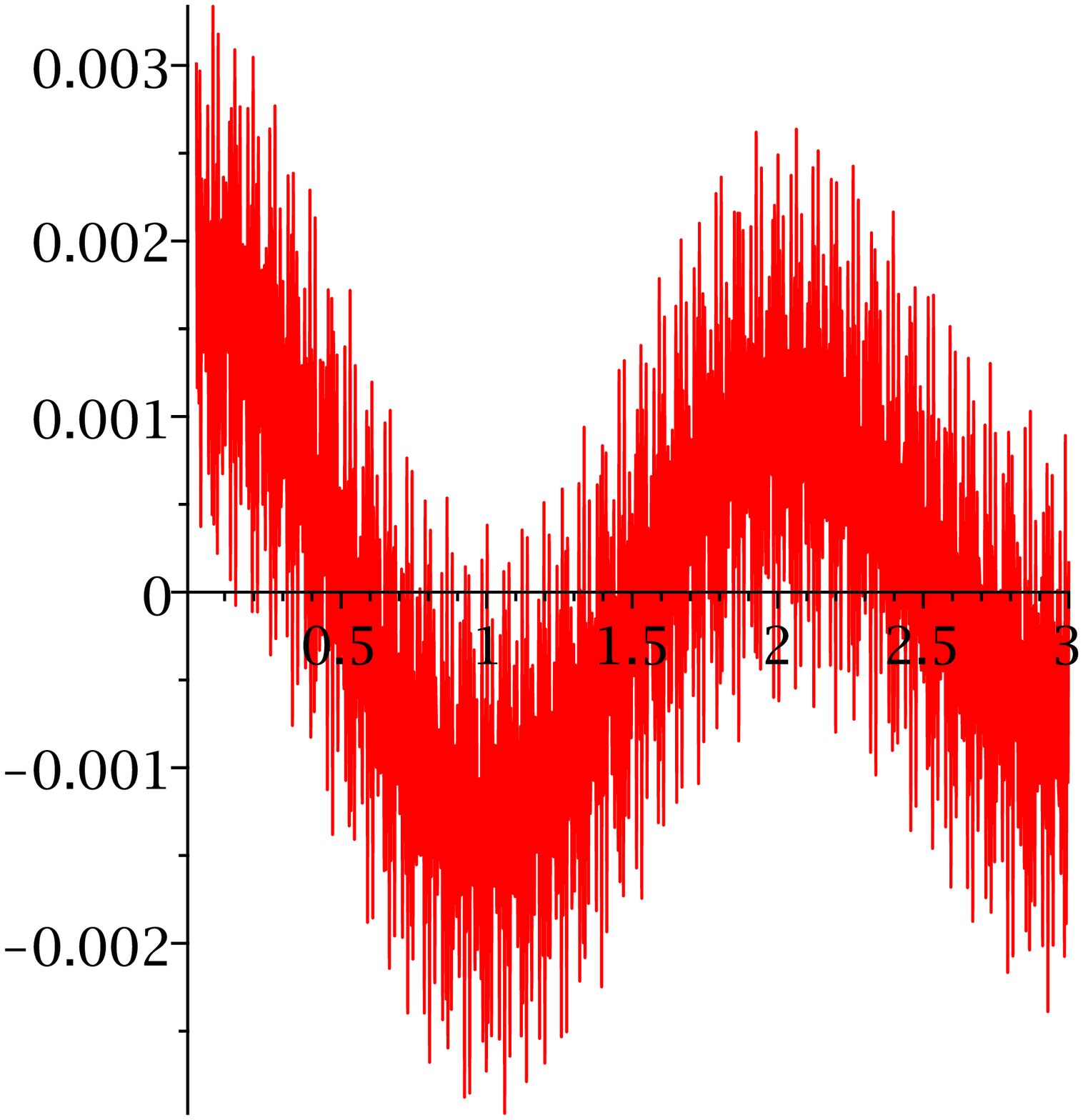}}
\quad
{\label{fig:4.10}\includegraphics[width=6cm,height=5cm]{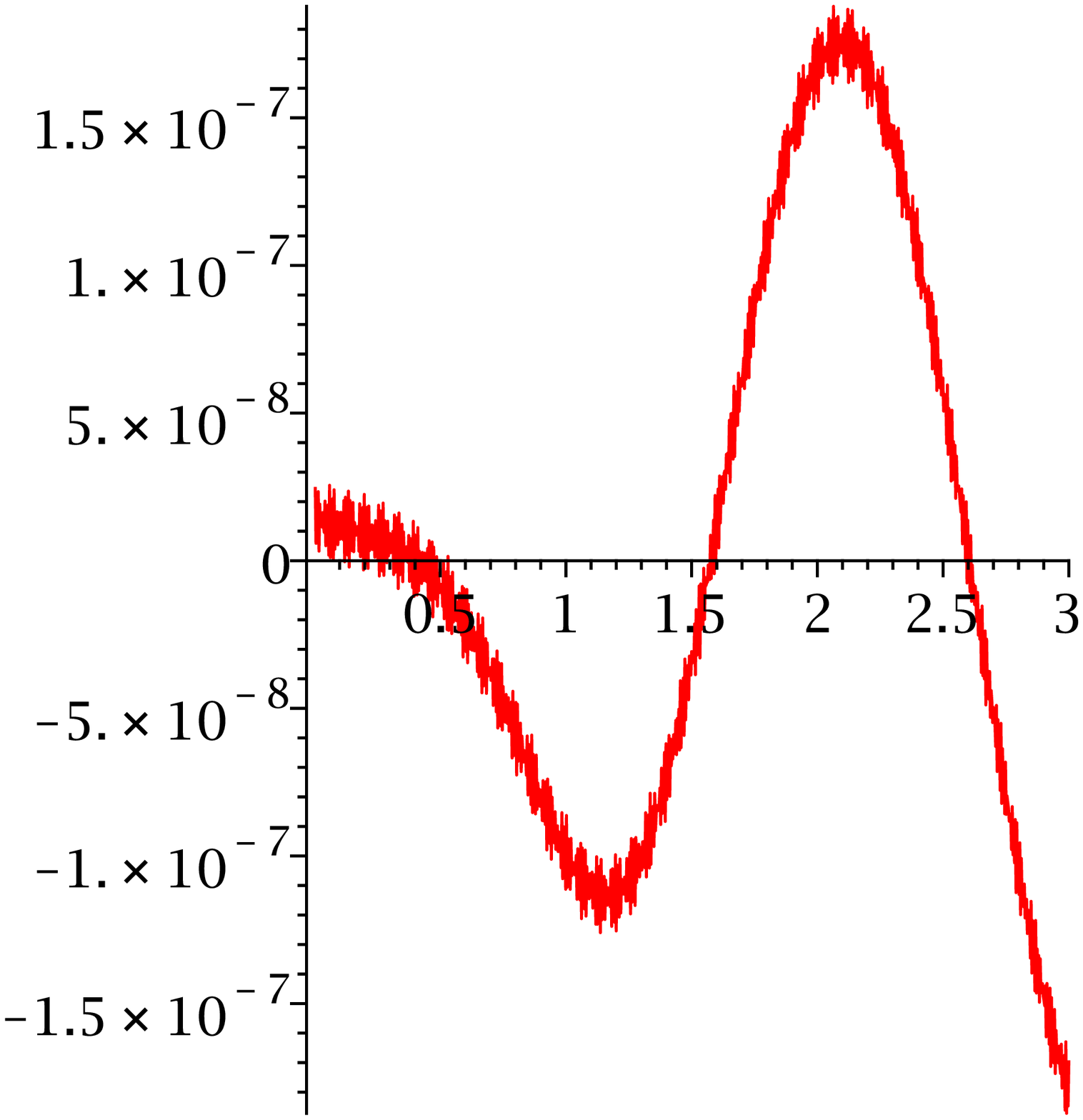}}\\
{\label{fig:4.10}\includegraphics[width=6cm,height=5cm]{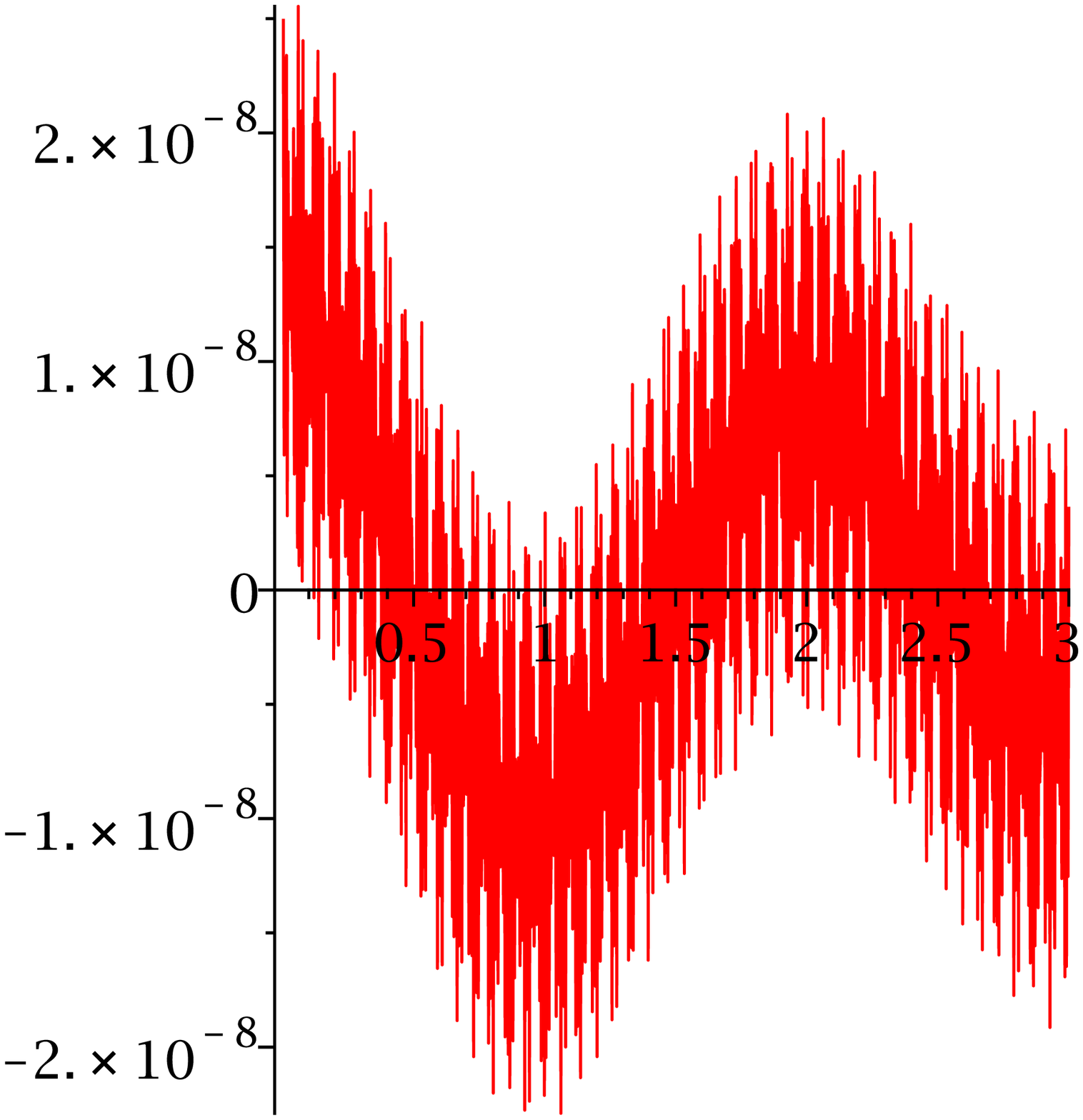}}
\quad
{\label{fig:4.10}\includegraphics[width=6cm,height=5cm]{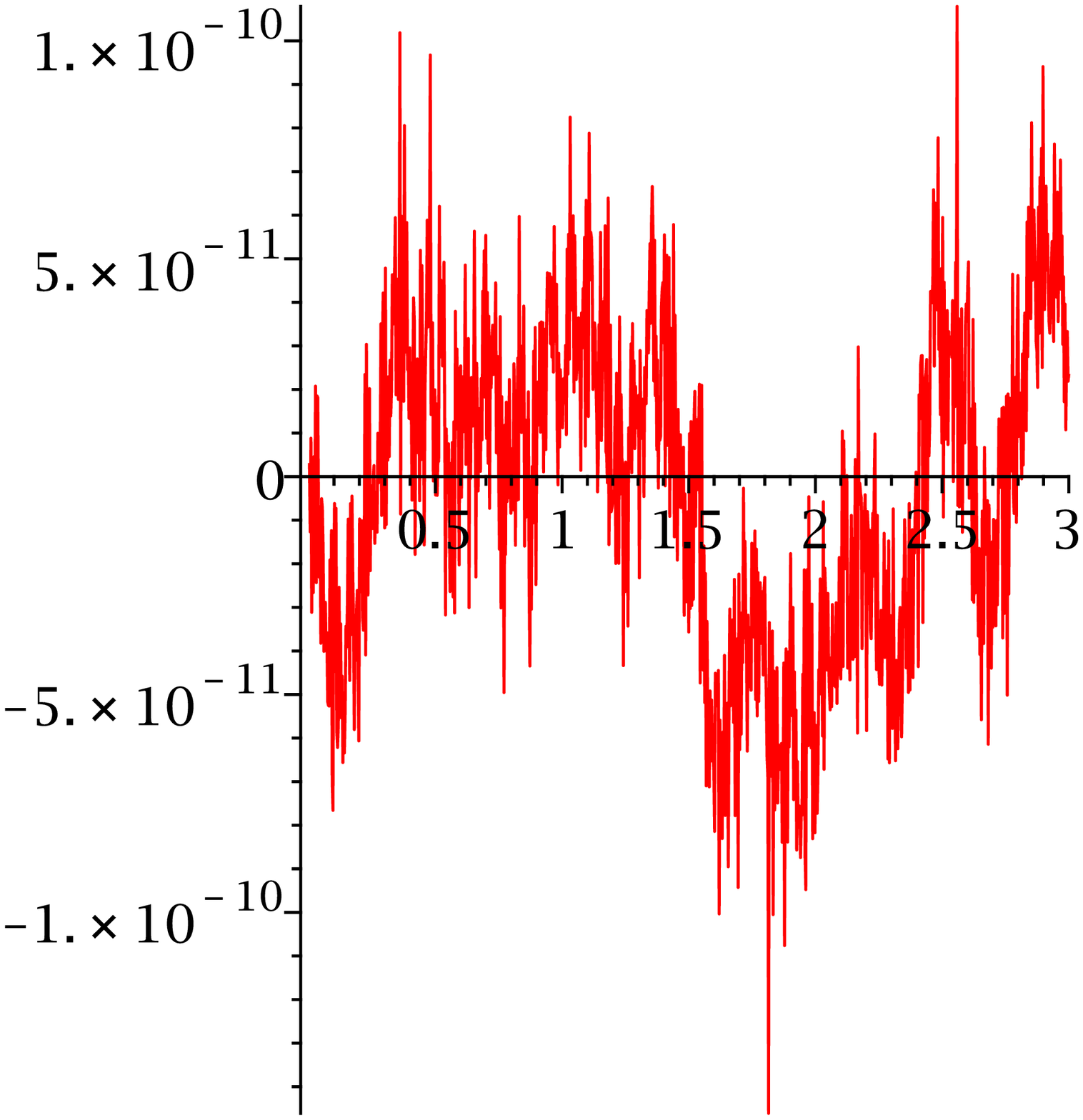}}\\
\caption{The top row: the real parts of error function with $s = 0$
(the left) and $s = 1$ (the right) for $y_5$ with $\omega = 100$.
The middle row: the real parts of error function with $s = 2$ (the
left) and $s = 3$ (the right) for $y_5$ with $\omega = 100$. The
third row: the real parts of error function with $s = 0$ (the left)
and $s = 1$ (the right) for $y_5$ with $\omega = 1000$. The fourth
row: the real parts of error function with $s = 2$ (the left) and $s
= 3$ (the right)for $y_5$ with $\omega = 1000$.}
\end{figure}

In addition, the CPU time is compared to the Runge-Kutta method
(rkf45) whose tolerance equals to $10^{-10}$. It takes $275$ seconds
for $\omega = 500$ and about $2992$ seconds for $\omega = 5000$,
respectively. The CPU time for the asymptotic method is about $11$
seconds for $\omega = 500$ and $13$ seconds for $\omega = 5000$. As
evident from our previous theoretical analysis, the computational
cost is about the same for the asymptotic method regardless of the
value of the oscillatory parameter.

\bigskip

\baselineskip=0.9\normalbaselineskip


{\small
\begin{thebibliography}{99}

\bibitem{Besicovitch:32}
{\rm Besicovitch, A. S. (1932)}, \textit{Almost Periodic Functions}, Cambridge Univ. Press,
Cambridge.

\bibitem{BoCheng:11}
{\rm Bao B. C., Shi G. D., Xu J. P., Liu Z. and Pan S. H. (2011)}. Dynamic
analysis of chaotic circuit with two memristors. \textit{Science China}, 54,
2180-2187.

\bibitem{Chartier:10}
{\rm Chartier, P., Murua, A. and Sanz-Serna, J. M. (2010)}, Higher-order averaging,
formal series and numerical integration I: B-series, \textit{Found. Comp. Maths.}
10, 695-727.

\bibitem{Chartier:12}
{\rm Chartier, P., Murua, A. and Sanz-Serna, J. M. (2012)}, Higher-order averaging,
formal series and numerical integration II: the quasi-periodic case, \textit{Found. Comp. Maths.}
to appear.

\bibitem{Chedjou:01}
{\rm Chedjou, J.C., Fotsin, H.B.,Woafo, P., Domngang, S.(2001)}, Analog simulation of the dynamics of a Van der
Pol oscillator coupled to a Duffing oscillator, \textit{IEEE Trans. Circ. Syst. I: Fundam. Theory Appl.}
48, 748¨C757.


\bibitem{Condon:10}
{\rm Condon, M., Dea\~{n}o, A. and Iserles, A. (2010)}, On systems of differential
equations with extrinsic oscillation, \textit{Discr. and Cont. Dynamical Sys.}
28, 1345-1367.

\bibitem{E:03}
{\rm E, W. and Engquist, B. (2003)}, The heterogeneous multiscale methods, \textit{Commun. Math. Sci.}
1, 87-132.

\bibitem{Fodjouong:07}
{\rm Fodjouong, G.J., Fotsin, H.B.,Woafo, P.(2007)}, Synchronizing modified van der Pol-Duffing oscillators with
offset terms using observer design: application to secure communications, \textit{Phys. Scr.}
75, 638¨C644.


\bibitem{Giannini:04}
{\rm Giannini, F. and Leuzzi, G. (2004)}, \textit{Nonlinear Microwave Circuit Design}, Wiley,
Chichester.

\bibitem{Iserles:06}
{\rm Iserles, A., N{\o}rsett, S. P. and Olver, S. (2006)}, Highly oscillatory quadrature:
The story so far, in A. Bermudez, ed., `Proceedings of ENuMath', Springer
Verlag, Berlin, pp. 97-118.


\bibitem{Ram:10}
{\rm Ram\'{\i}rez, F., Su\'{a}arez, A., Lizarraga, I. and Collantes, J.-M. (2010)}, Stability
analysis of nonlinear circuits driven with modulated signals, \textit{IEEE Trans.
Microwave Theory Tech.}, 58, 929-940.

\bibitem{Sanz:09}
{\rm Sanz-Serna, J. M. (2009)}, Modulated Fourier expansions and heterogeneous
multiscale methods, \textit{IMA J. Numer. Anal.}, 29, 595-605.

\bibitem{Slight:08}
{\rm Slight, T.J., et al.(2008)}, A Lienard oscillator resonant tunnelling diode-laser diode hybrid integrated circuit:
model and experiment, \textit{IEEE J. Quantum Electron.}, 44, 1158¨C1163.

\bibitem{Verhulst:90}
{\rm Verhulst, F. (1990)}, \textit{Nonlinear Differential Equations and Dynamical Systems},
Springer Verlag, Heidelberg.

\end{thebibliography}
}
\end{document}